\newif\ifpreprint
\preprinttrue  
\ifpreprint
\documentclass[12pt]{article}
\pdfoutput=1

\usepackage{fullpage}

\usepackage{amsmath,amssymb,amsthm}

\usepackage{algorithm}
\usepackage[noend]{algpseudocode} 

\newcommand{\citep}{\cite}

\else
\documentclass[referee,sn-mathphys]{sn-jnl}

\normalbaroutside

\jyear{2023}%

\theoremstyle{thmstyleone}%
\newtheorem{theorem}{Theorem}
\newtheorem{proposition}[theorem]{Proposition}%
\newtheorem{corollary}[theorem]{Corollary}%
\newtheorem{lemma}[theorem]{Lemma}%

\theoremstyle{thmstyletwo}%
\newtheorem{example}{Example}%
\newtheorem{remark}{Remark}%

\theoremstyle{thmstylethree}%
\newtheorem{definition}{Definition}%
\newtheorem{assumption}{Assumption}%

\fi

\usepackage{subfig}
\usepackage{verbatim}
\usepackage{graphicx}

\ifpreprint

\newcommand{\reviewChanges}{}

\fi
\usepackage{mathtools}
\usepackage{enumitem}
\usepackage{multirow}

\usepackage{csvsimple}	
\usepackage{booktabs}   
\usepackage{adjustbox}  
\usepackage{appendix}
\usepackage{multicol}
\usepackage{hyperref}

\usepackage[acronyms]{glossaries}   
\makeglossaries

\newcommand{\uncset}{{\mathcal{U}}}

\newcommand{\eg}{{\it e.g.}}
\newcommand{\ie}{{\it i.e.}}

\newcommand{\ones}{\mathbf 1}
\newcommand{\reals}{{\mbox{\bf R}}}

\newcommand{\prob}{{\mathbf P}}

\newcommand{\Tr}{\mathop{\bf tr}}

\newcommand{\Expect}{{\mbox{\bf E}}}

\newcommand{\dom}{\mathop{\bf dom}_u} 

\newcommand{\dif}{\text{d}}

\newcommand{\define}{=}
\newcommand{\defn}{\define}
\newcommand{\norm}[1]{\| #1 \|}
\renewcommand{\Re}{\reals}
\renewcommand{\tfrac}[2]{#1/#2}

\newcommand{\indicator}{\mathcal{I}}
\newcommand{\CVaR}{\mathop{\bf CVaR}}
\newcommand{\VaR}{\mathop{\bf VaR}}
\newcommand{\supp}{S}

\ifpreprint

\newtheorem{theorem}{Theorem}[section]  
\newtheorem{lemma}{Lemma}[section]  
\newtheorem{corollary}{Corollary}[theorem]  

\newtheorem{assumption}{Assumption}[section]

\newtheorem{remark}{Remark}[section]

\fi

\ifpreprint
\newcommand{\set}[2]{\{ #1 \mid #2 \}}

\else
\renewcommand{\set}[2]{\{ #1 \mid #2 \}}

\fi

\newacronym{SAA}{SAA}{sample average approximation}
\newacronym{LO}{LO}{linear optimization}
\newacronym{QO}{QO}{quadratic optimization}
\newacronym{MIQO}{MIQO}{mixed-integer quadratic optimization}
\newacronym{MIO}{MIO}{mixed-integer optimization}
\newacronym{MILO}{MILO}{mixed-integer linear optimization}
\newacronym{MINLO}{MINLO}{mixed-integer nonlinear optimization}
\newacronym{sBB}{sBB}{spacial branch and bound}
\newacronym{NLO}{NLO}{nonlinear optimization}
\newacronym{PWA}{PWA}{piecewise affine}
\newacronym{SVM}{SVM}{support vector machines}
\newacronym{ReLU}{ReLU}{rectified linear unit}
\newacronym{CPU}{CPU}{central processing unit}
\newacronym{GPU}{GPU}{graphics processing unit}
\newacronym{MPC}{MPC}{model predictive control}
\newacronym{ADMM}{ADMM}{alternating direction method of multipliers}
\newacronym{ADP}{ADP}{approximate dynamic programming}
\newacronym{FPGA}{FPGA}{field-programmable gate array}
\newacronym{DRO}{DRO}{distributionally robust optimization}
\newacronym{RO}{RO}{robust optimization}
\newacronym{MRO}{MRO}{mean robust optimization}

\begin{document}

\ifpreprint
\title{Mean Robust Optimization}
\author{Irina Wang, Cole Becker, Bart Van Parys, and Bartolomeo Stellato}
\maketitle
\begin{abstract}
Robust optimization is a tractable and expressive technique for decision-making under uncertainty, but it can lead to overly conservative decisions when pessimistic assumptions are made on the uncertain parameters.
Wasserstein distributionally robust optimization can reduce conservatism by being data-driven, but it often leads to very large problems with prohibitive solution times.
We introduce mean robust optimization, a general framework that combines the best of both worlds by providing a trade-off between computational effort and conservatism.
We propose uncertainty sets constructed based on clustered data rather than on observed data points directly thereby significantly reducing problem size.
By varying the number of clusters, our method bridges between robust and Wasserstein distributionally robust optimization.
We show finite-sample performance guarantees and explicitly control the potential additional pessimism introduced by any clustering procedure.
In addition, we prove conditions for which, when the uncertainty enters linearly in the constraints, clustering does not affect the optimal solution.
We illustrate the efficiency and performance preservation of our method on several numerical examples, obtaining multiple orders of magnitude speedups in solution time with little-to-no effect on the solution quality.


\end{abstract}

\fi

\section{Introduction}

\Gls{RO} and \gls{DRO} are popular tools for decision-making under uncertainty due to their high expressiveness and versatility.
The main idea of \gls{RO} is to define an uncertainty set and to minimize the worst-case cost across possible uncertainty realizations in that set.
However, while \gls{RO} often leads to tractable formulations, it can be overly-conservative~\citep{reduce_conserve}.
To reduce conservatism, \gls{DRO} takes a probabilistic approach, by modeling the uncertainty as a random variable following a probability distribution known only to belong to an uncertainty set (also called ambiguity set) of distributions.
%
In both \gls{RO} and \gls{DRO}, the choice of the uncertainty or ambiguity set can greatly influence the quality of the solution for both paradigms.
Good-quality uncertainty sets can lead to excellent practical performance while ill chosen sets can lead to overly-conservative actions and intractable computations.

Traditional approaches design uncertainty sets based on theoretical assumptions on the uncertainty distributions~\citep{ben-tal_robust_2000,bandi_tractable_2012,ben-tal_robust_2009,bertsimas_price_2004}.
While these methods have been quite successful, they rely on a priori assumptions that are difficult to verify in practice.
On the other hand, the last decade has seen an explosion in the availability of data.
This change has brought a shift in focus from a priori assumptions on the probability distributions to data-driven methods in operations research and decision sciences.
In \gls{RO} and \gls{DRO}, this new paradigm has fostered data-driven methods where uncertainty sets are shaped directly from data~\citep{bertsimas_data-driven_2018}.
In data-driven \gls{DRO}, a popular choice of the ambiguity set is the ball of distributions whose Wasserstein distance to a nominal distribution is at most $\epsilon > 0$~\cite{mohajerin_esfahani_data-driven_2018,kuhn2019wasserstein,DBLP:gao2020,gao2016distributionally}.
When the reference distribution is an empirical distribution, the associated Wasserstein \gls{DRO} can be formulated as a convex minimization problem where the number of constraints grows linearly with the number of datapoints~\citep{mohajerin_esfahani_data-driven_2018}.
While less conservative than \gls{RO}, data-driven \gls{DRO} can lead to very large formulations that are intractable, especially in~\gls{MIO}.

A common idea to reduce the dimensionality of data-driven decision-making problems  is to use clustering techniques from machine learning.
While clustering has recently appeared in various works within the stochastic programming literature~\citep{SAA1, SAA2, SAA3, SAA4}, the focus has been on the improvement of and comparisons to the \gls{SAA} approach and not in a distributionally robust sense.
In contrast, recent approaches in the \gls{DRO} literature cluster data into partitions and either build moment-based uncertainty sets for each partition~\citep{Rsome, perakis2022robust}, or enrich Wasserstein \gls{DRO} formulations with partition-specific information (\eg, relative weights)~\cite{partition}.
While these approaches are promising, clustering is still used as a pre-processing heuristic on the data-sets in~\gls{DRO}, without a clear understanding of how it affects the conservatism of the optimal solutions.
In particular, choosing the right clustering parameters to carefully balance computational tractability and out-of-sample performance is still an unsolved challenge.




\subsection{Our contributions}%
\label{sub:our_contributions}
In this work, we present~\gls{MRO}, a data-driven method that, via machine learning clustering, bridges between \gls{RO} and  Wasserstein \gls{DRO}.
\begin{itemize}
	\item We design the uncertainty set for \gls{RO} as a ball around clustered data. Without clustering, our formulation corresponds to the finite convex reformulation in Wasserstein~\gls{DRO}.
With just one cluster, our formulation corresponds to the classical \gls{RO} approach.
The number of clusters is a tunable parameter that provides a tradeoff between the worst-case objective value and computational efficiency, which includes both speed and memory usage.
	\item We provide probabilistic guarantees of constraint satisfaction for our method, based on the quality of the clustering procedure.
	\item We derive bounds on the effect of clustering in case of constraints with concave \reviewChanges{and maximum-of-concave} dependency on the uncertainty. In addition, we show that, when constraints are linearly affected by the uncertainty, clustering does not affect the solution nor the probabilistic guarantees.
	\item We show on various numerical examples that, thanks to our clustering procedure, our approach provides multiple orders of magnitude speedups over classical approaches while guaranteeing the same probability of constraint satisfaction. The code to reproduce our results is available at \url{https://github.com/stellatogrp/mro_experiments}.
\end{itemize}

\subsection{Related work}%
\label{sub:related_work}
\paragraph{Robust optimization.}
\gls{RO} deals with decision-making problems where some of the parameters are subject to uncertainty.
The idea is to restrict data perturbations to be within a deterministic uncertainty set, then optimize the worst-case performance across all realizations of this uncertainty.
For a detailed overview of \gls{RO}, we refer to the survey papers by Ben-Tal and Nemirovski~\cite{robustconvexopt} and Bertsimas et al.~\cite{bertsimassurvey}, as well as the books by Ben-Tal et al.~\cite{ben-tal_robust_2009} and Bertsimas and den Hertog~\cite{robustadaptopt}.
These approaches, while powerful, may be overly-conservative, and there have been approaches that provide a tradeoff between conservatism and constraint violation~\citep{reduce_conserve}.

\paragraph{Distributionally robust optimization.}
\gls{DRO} minimizes the worst-case expected loss over a probabilistic ambiguity set characterized by certain known properties of the true data-generating distribution.
Based on the type of ambiguity set considered existing literature on \gls{DRO} can roughly be defined in two.
Ambiguity sets of the first type contain all distributions that satisfy certain moment constraints~\citep{dro_1,dro_2,dro_3,dro_4}.
In many cases such ambiguity sets possess a tractable formulation, but have also been criticized for yielding overly conservative solutions~\citep{dro_5}.
Ambiguity sets of the second type enjoy the interpretation of a ball of distributions around a nominal distribution, often the empirical distribution on the observed samples.
Wasserstein uncertainty sets are one particular example \citep{mohajerin_esfahani_data-driven_2018,kuhn2019wasserstein,DBLP:gao2020,gao2016distributionally} and enjoy both a tractable primal as well as a tractable dual formulation.
We refer to the work by Chen and Paschalidis~\cite{OPT-026} for a thorough overview of \gls{DRO}, and to the work by Zhen et al.~\cite{zhen2021mathematical} for a general theory on convex dual reformulations.
When the ambiguity set is well chosen, \gls{DRO} formulations enjoy strong out-of-sample statistical performance guarantees.
As these statistical guarantees are typically not very sharp, in practice the radius of the uncertainty set is typically chosen through time consuming cross-validation \citep{DBLP:gao2020}.
At the same time, \gls{DRO} has the downside of being more computationally expensive than traditional robust approaches.
We observe for instance that the number of constraints in Wasserstein~\gls{DRO} formulations scale linearly with the number of samples, which can become practically prohibitive especially when integer variables are involved.
Our proposed method addresses this problem by reducing the number of constraints through clustering.
While many works have recently emerged on the construction of \gls{DRO} ambiguity sets through the partitioning of data~\cite{Rsome, partition,perakis2022robust}, or the discretization of the underlying distribution~\cite{droapprox}, there still exists a gap in the literature. In particular, theoretical bounds on the change in problem performance as affected by the number of clusters, as well as by the quality of the cluster assignment, remain largely unexplored.
In this work, we fill the gap by providing such insights.

\paragraph{Data-driven robust optimization.}
Data-driven optimization has been well-studied, with various techniques to learn the unknown data-generating distribution before formulating the uncertainty set.
Bertsimas et al.~\cite{bertsimas_data-driven_2018} construct the ambiguity set as a confidence region for the unknown data-generating distribution $\prob$ using several statistical hypothesis tests.
By pairing a priori assumptions on $\prob$ with different statistical tests, they obtain various data-driven uncertainty sets, each with its own geometric shape, computational properties, and modeling power.
We, however, use machine learning in the form of clustering algorithms to preserve the geometric shape of the dataset, without explicitly learning and parametrizing the unknown distribution.

\paragraph{Distributionally robust optimization as a robust program.}
Gao and Kleywegt~\cite{gao2016distributionally} consider a robust formulation of Wasserstein \gls{DRO} similar to our mean robust optimization, but without the idea of dataset reduction.
\reviewChanges{
Given $N$ samples and a positive integer $K$, they introduce an approximation of Wasserstein \gls{DRO} by defining a new ambiguity set as a subset of the standard Wasserstein \gls{DRO} set, containing all distributions supported on $NK$ points with equal probability $1/(NK)$, as opposed to the standard set supported on $N$ points.}
In this work, however, we study how to reduce, instead of increase, the number of variables and constraints to make the Wasserstein \gls{DRO} problem more tractable by linking it to robust optimization.

\reviewChanges{\paragraph{Robust optimization as a distributionally robust optimization program.} Xu et al.~\cite{RO_DRO} take inspiration from sample-based optimization problems to investigate probabilistic interpretations of \gls{RO}. 
They generalize the ideas of Delage and Ye~\cite{dro_3}, that the solution to a robust optimization problem is the solution to a special Distributionally Robust Stochastic Program (DRSP), where the distributional set contains all distributions whose support is contained
in the uncertainty set.
In a related vein, Bertsimas et al.~\cite{wass-p-1} show that, under a particular construction of the uncertainty sets, multi-stage stochastic linear optimization can be interpreted as Wasserstein-$\infty$ \gls{DRO}. 
We establish a similar equivalence between \gls{RO} and \gls{DRO}, focusing especially on Wasserstein-$p$ ambiguity sets for all $p$. 
We develop an easily interpretable construction of the primal constraints and uncertainty sets, and prove, in view of both the primal and dual problems, that $p = \infty$ is a limiting case of $p \geq 1$. 
This provides a natural extension of the equivalence proved in~\cite[Proposition 3]{wass-p-1}.}

\paragraph{Probabilistic guarantees in robust and distributionally optimization.}
Bertsimas et al.~\cite{bertsimas_probabilistic} propose a disciplined methodology for deriving probabilistic guarantees for solutions of robust optimization problems with specific uncertainty sets and objective functions.
They derive a posteriori guarantee to compensate for the conservatism of a priori uncertainty bounds.
Esfahani and Kuhn~\cite{mohajerin_esfahani_data-driven_2018} obtain finite-sample guarantees for Wasserstein \gls{DRO} for selecting the radius $\epsilon$ of order $N^{-1/\max\{2,m\}}$, where $N$ is the number of samples and $m$ is the dimension of the problem data, while Gao~\cite{DBLP:gao2020} derives finite-sample guarantees for Wasserstein \gls{DRO} for selecting $\epsilon$ of order $N^{-1/2}$ under specific assumptions.
We provide theoretical results of a similar vein, with a slightly increased $\epsilon$ to compensate for information lost through clustering and achieve the same probabilistic guarantees.
\reviewChanges{Our theoretical guarantees hold for Wasserstein-$p$ distance for all $p\geq1$ and $p = \infty$, and are independent of the uncertain function to minimize. 
These bounds, however, following the literature, are theoretical in nature and not tight in practice, typically resulting in overly-conservative $\epsilon$.
The final $\epsilon$ values are usually chosen through empirical calibration --- in which case, our formulation, by being lower dimensional, is overall much faster to solve.}

\paragraph{Clustering in stochastic optimization.}
Clustering in stochastic optimization is closely related to the idea of {\it scenario reduction}.
First introduced by Dupa{\v c}ov{\'a} et al.~\cite{scen_1}, scenario reduction seeks to approximate, with respect to a probability metric, an $N$-point distribution with a distribution with a smaller number of points.
Rujeerapaiboon et al.~\cite{scen_k} analyze the worst-case bounds on the approximation error for scenario reduction with respect to the Wasserstein metric, for initial distributions constrained to a unit ball.
They provide constant-factor approximation algorithms for performing scenario reduction through $K$-medians and $K$-means clustering~\citep{kmeans}.
Later, Bertsimas and Mundru~\cite{scen_3} apply this idea to two-stage stochastic optimization problems, and provide an alternating-minimization method for finding optimal reduced scenarios under the modified objective.
They also provide performance bounds on the stochastic optimization problem for different scenarios.
Jacobson et al.~\citep{SAA1}, Emelogu et al.~\cite{SAA2}, Beraldi et al.~\cite{SAA3}, and Chen \cite{SAA4} apply a similar idea of clustering to reduce the sample/scenario size, then compare the results against the classical \gls{SAA} approach where the sample size is not reduced.
In \gls{MRO}, we adapt and extend the scenario reduction approach to Wasserstein \gls{DRO}, where upon fixing the reduced scenario points to ones found by the clustering algorithm, we allow for variation around these reduced points.
We then provide performance bounds on the \gls{DRO} problem depending on the number of clusters.

\paragraph{Data compression in data-driven problems.}
Fabiani and Goulart~\cite{data-comp} compress data for robust control problems by minimizing the Wasserstein-1 distance between the original and compressed datasets, and observe a slight loss in performance in exchange for reduced computation time. While related, this is orthogonal to our approach of using machine learning clustering to reduce the dataset, where we include results and theoretical bounds for a more general set of robust optimization problems with Wasserstein-$p$ distance, and demonstrate conditions under which no performance loss is necessary.

\subsection{Layout of the paper}%
\reviewChanges{
In Section~\ref{sec:mro}, we present our approach for concave uncertainty constraints, then extend the results to maximum-of-concave functions. 
In Section~\ref{sub:connectionsdro}, we present connections to distributionally robust optimization, and give theoretical guarantees on constraint satisfaction.
In Section~\ref{sec:conservatism}, we analyze the effect on clustering on the worst-case value of the~\gls{MRO} solutions for both concave and maximum-of-concave constraints}.
In Section~\ref{sec:parameter_selection}, we give guidelines for choosing hyperparameters.
In Section~\ref{sec:examples}, we provide computational verification of the speedups obtained through our methodology.
In Section~\ref{sec:conclusions}, we summarize our conclusions.

\section{Mean robust optimization}%
\label{sec:mro}
\subsection{The problem}
\label{sec:intro}
We consider an uncertain constraint of the form
\begin{equation}
	\label{eq:constr_original}
	g(u,x) \le 0,
\end{equation}
where $x \in \mathcal{X} \subseteq \reals^n$ is an optimization variable on the compact set $\mathcal{X}$, $u \in S \subseteq \reals^m$ is an uncertain parameter with convex and closed support $S$, and $-g(u, x)$ is proper, convex, and lower-semicontinuous in~$u$ for all $x$.
Throughout this paper, we assume the support $\supp$ of $u$ to live within the domain of $g$ for the variable $u$, which we will refer to as $\dom{g}$, \ie, $\supp \subseteq \dom{g}$.
We assume $\dom{g}$ is independent of $x$, and that the following assumption holds. 
\begin{assumption}
	\label{assp:dom}
	The domain $\dom{g}$ is $\reals^m$. \reviewChanges{Otherwise, $g$ is either element-wise monotonically increasing in $u$ and only has a (potentially) lower-bounded domain, or element-wise monotonically decreasing in $u$ and only has a (potentially) upper-bounded domain.}
\end{assumption}
\reviewChanges{This assumption on the domain and monotonicity of $g$ is very common in practice as it is satisfied by linear and quadratic functions, as well as other common functions (\eg, $\log(u)$, and $1/(1+u)$).}

\reviewChanges{In Section~\ref{sec:max_concave}, we extend our results for $g$ being the maximum of concave functions, each satisfying the aforementioned conditions.}

The \gls{RO} approach defines an uncertainty set~$\uncset \subseteq \reals^m$ and forms the {\it robust counterpart} as
\begin{equation*}
	g(u,x) \le 0, \quad \forall u \in \uncset,
\end{equation*}
where the uncertainty set is chosen so that for any solution $x$, the above holds with a certain probability. We define this in terms of expectation,
\begin{equation}
	\label{eq:expectation}
	\Expect^\prob (g(u,x)) \le 0,
\end{equation}
where $\prob$ is the unknown distribution of the uncertainty $u$.

\reviewChanges{
\paragraph{Risk measures.}%
\label{par:risk_measures}
Expectation constraints of the form~\eqref{eq:expectation} can represent popular risk measures, and can imply constraints commonly used in chance-constrained programming (CCP). 
In CCP, the probabilistic constraint considered is
  \begin{equation}
		\label{eq:chance}
		\prob(g(u,x)\leq 0) \geq 1-\alpha,
	 \end{equation}
which corresponds to the {\it value at risk} being nonpositive, \ie,
\begin{equation*}
\VaR(g(u,x),\alpha) = \inf\{\gamma \mid \prob(g(u,x)\leq \gamma) \geq 1-\alpha\} \leq 0.
\end{equation*}
Unfortunately, except in very special cases, the value at risk function is intractable~\citep{cvaropt}.
A tractable approximation of the value at risk is the {\it conditional value at risk}~\citep{cvaropt,cvar}, defined as
\begin{equation*}
\CVaR(g(u,x),\alpha) = \inf_{\tau}\{\Expect(\tau + (1/\alpha)(g(u,x)- \tau)_+)\},
\end{equation*}
where $(a)_+ = \max\{a,0\}$. 
This expression can be modeled through our approach, by writing $\CVaR(g(u,x),\alpha) = \inf_{\tau}\{\Expect(\hat{g}(u,x,\tau))\}$, where ${\hat{g}(u,x,\tau) = \tau + (1/\alpha)(g(u,x)- \tau)_+}$ is the maximum of concave functions, which we study in Sections~\ref{sec:max_concave},~\ref{sub:sparse_portfolio_optimization}, and~\ref{sub:facility}.
It is well known from~\cite{cvaropt} that the relationship between these probabilistic guarantees of constraint satisfaction is
\begin{equation*}
\CVaR(g(u,x),\alpha) \leq 0 ~~ \Longrightarrow~~ \VaR(g(u,x),\alpha)\leq 0 ~~\Longleftrightarrow ~~\prob(g(u,x)\leq 0) \geq 1 - \alpha.
\end{equation*}
Therefore, our expectation constraint implies common chance constraints.
}

\paragraph{Finite-sample guarantees.}%
\label{par:finite_sample_guarantees}
In data-driven optimization, while $\prob$ is unknown, it is partially observable through a finite set of $N$ independent samples of the random vector $u$.
We denote the training dataset of these samples by ${\mathcal{D}}_N \define \{d_i\}_{i\leq N} \subseteq \supp$, and note that this dataset is governed by $\prob^N$, the product distribution supported on $\supp^{N}$.
A data-driven solution of a robust optimization problem is a feasible decision $\hat{x}_N \in \reals^n$ found using the data-driven uncertainty set $\uncset$, which in turn is constructed by the training dataset ${\mathcal{D}}_N $.
Specifically, the feasible decision and data-driven uncertainty set $\uncset$ we construct must imply the probabilistic guarantee
\begin{equation}
	\label{eq:prob_guarantees1}
	\prob^N\left( \Expect^{\prob}(g(u, \hat{x}_N)) \le 0\right) \ge 1 - \beta,
\end{equation}
where $\beta > 0$ is the specified probability of constraint violation. From now on, when we refer to probabilistic guarantees of constraint satisfaction, it will be a reference to~\eqref{eq:prob_guarantees1}.

\subsection{Our approach}
\begin{figure}[h]%
	\centering
	\includegraphics[width=0.6\textwidth]{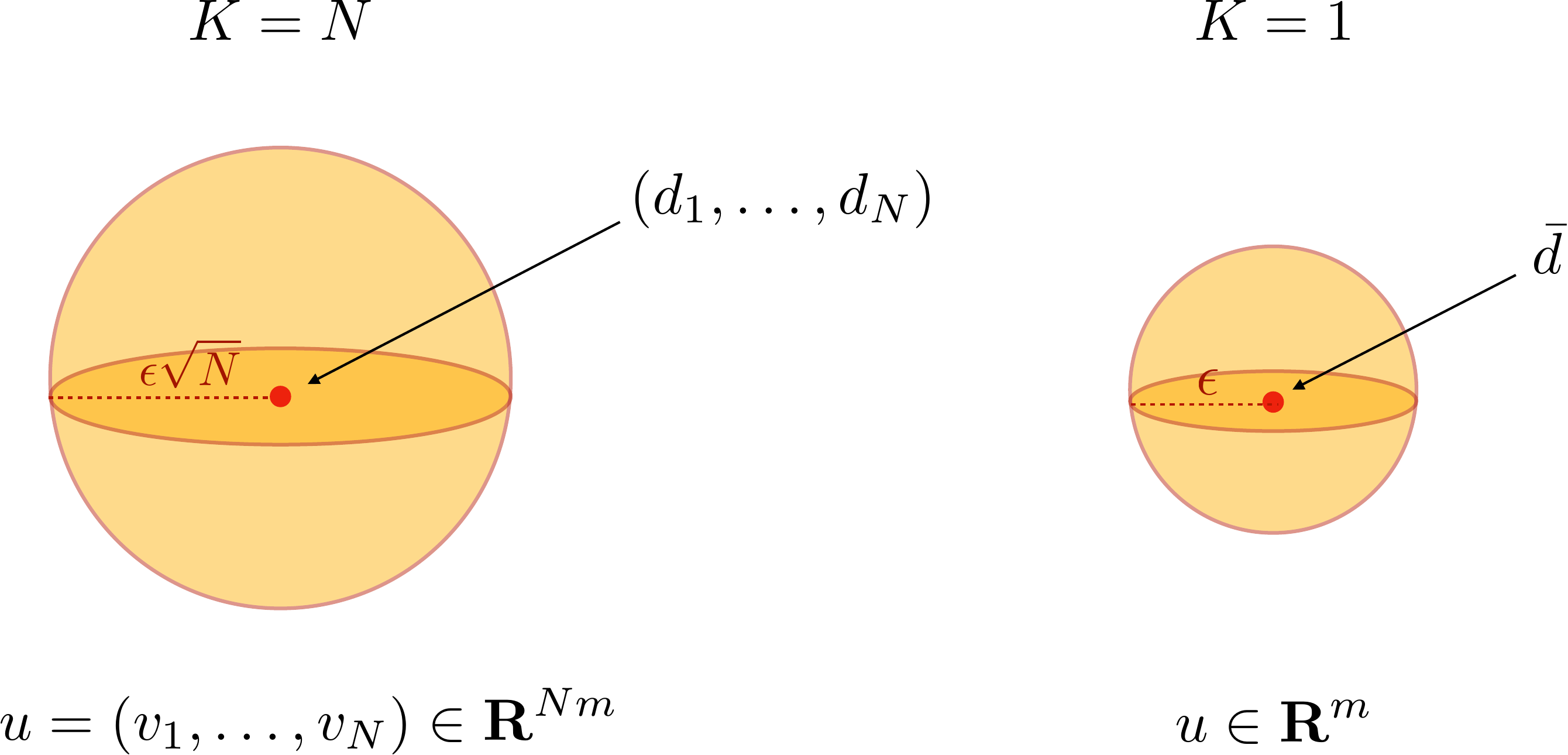}
	\caption{Visualizing the uncertainty set $\uncset(N, \epsilon)$ and $\uncset(1,\epsilon)$ as high dimension balls when $p = 2$.}%
	\label{fig:unset}%
\end{figure}

To meet the probabilistic guarantees outlined above, we propose to construct $\hat{x}_N$ to satisfy particular constraints, with respect to a particular uncertainty set. 
\paragraph{Case $p \geq 1$.}
In the case where $p \geq 1$, the set we consider takes the form
\begin{equation*}
	\label{eq:uncset_p}
	\uncset(K, \epsilon) \define \left\{ u = (v_1,\dots,v_K) \in \supp^{K} \quad\middle|\quad  \sum_{k=1} ^K w_k \| v_k - \bar{d}_k \|^p \le \epsilon^p  \right\},
\end{equation*}
where we partition ${\mathcal{D}}_N$ into $K$ disjoint subsets $C_k$, and $\bar{d}_k$ is the centroid of the $k$-th subset, for $k=1,\dots,K$.
The weight $w_k>0$ of each subset is equivalent to the proportion of points in the subset, \ie, $w_k = |C_k|/N$.
We choose $p$ to be an integer exponent, and $\epsilon > 0$ will be chosen depending on the other parameters to ensure satisfaction of the probability guarantee~\eqref{eq:prob_guarantees1}.
When $p = 2$ and $\supp = \reals^m$, the set can be visualized as an ellipsoid in $\mathbf{R}^{Km}$ with the center formed by stacking together all $\bar{d}_k$ into a single vector of dimension $\mathbf{R}^{Km}$.
When we additionally have $K = N$ or $K = 1$, this ellipsoid becomes a ball of dimension $\reals^{Nm}$ or $\reals^{m}$ respectively, as shown in Figure~\ref{fig:unset}.

\begin{figure}[h]%
	\centering
	\includegraphics[width=0.6\textwidth]{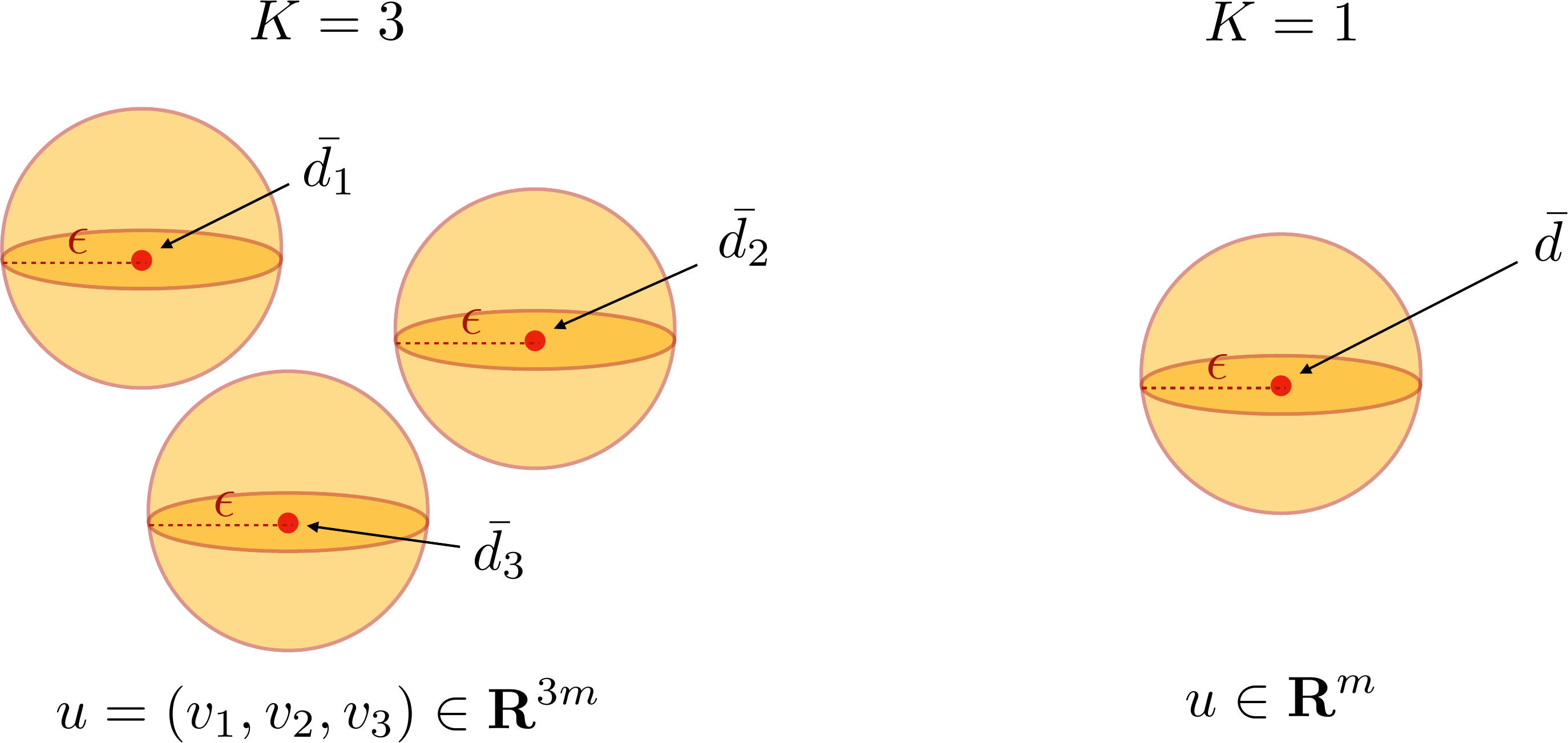}
	\caption{Visualizing the decoupled uncertainty set $\uncset(K, \epsilon)$ with $p = \infty$.}%
	\label{fig:unset-p}%
\end{figure}

\paragraph{Case $p = \infty$.}
In the case where $p = \infty$, the set we consider takes a more specific form,
\begin{equation*}
	\label{eq:uncset_inf}
	\uncset(K, \epsilon) \define \left\{ u = (v_1,\dots,v_K) \in \supp^{K} \quad\middle|\quad \max_{k=1,\dots,K} \| v_k - \bar{d}_k \| \le \epsilon\right\},
\end{equation*}
where the constraints for individual $v_k$ become decoupled.
See Figure~\ref{fig:unset-p} for an example when $K = 3$ and $K = 1$.
This decoupling follows the result for the Wasserstein type $p = \infty$ metric \citep[Equation 2]{wass-p}, as our uncertainty set is analogous to the set of all distributions within Wasserstein-$\infty$ distance of $\bar{d}$.
We note that, if any of the decoupled constraints are violated, then $\lim_{p \rightarrow \infty} \sum_{k=1} ^K w_k \| v_k - \bar{d}_k \|^p \geq \epsilon^p$, and the summation constraint will be violated.

For both cases, $p \geq 1$ and $p = \infty$, when~$K = 1$, we have a simple uncertainty set: a ball of radius $\epsilon$ around the empirical mean of the entire dataset, $\uncset(1, \epsilon) \define \left\{ v \in \supp \mid \| v - \bar{d} \| \le \epsilon  \right\}.$
This is equivalent to the uncertainty set of traditional \gls{RO}, as it is of the same dimension $m$ as the uncertain parameter.
When $K=N$ and $w_k = 1/N$, both cases closely resemble the ambiguity sets of Wasserstein-$p$ \gls{DRO}.

Having defined the uncertainty set, we now introduce constraints of the form
\begin{equation}
	\label{eq:constraints}
\bar{g}(u, x) \define  \sum_{k=1} ^K w_k g(v_k, x),
\end{equation}
where $g$ is defined in the original constraint~\eqref{eq:constr_original}.
The weights $w_k$ correspond to the ones defined in the uncertainty set.
Putting everything together, $\hat{x}_N$ is the solution to the robust optimization problem
\begin{equation}
	\label{eq:robustopt}
	\tag{MRO}
	\begin{array}{ll}
		\mbox{minimize} & f(x)\\
		\mbox{subject to} & \bar{g}(u, x)  \le 0  \quad \forall u \in \uncset(K, \epsilon),
	\end{array}
	\end{equation}
	where $f$ is the objective function. We call this problem the \glsentryfull{MRO} problem.

\begin{figure}[h]%
		\centering
		\subfloat{{\includegraphics[width=0.9\linewidth]{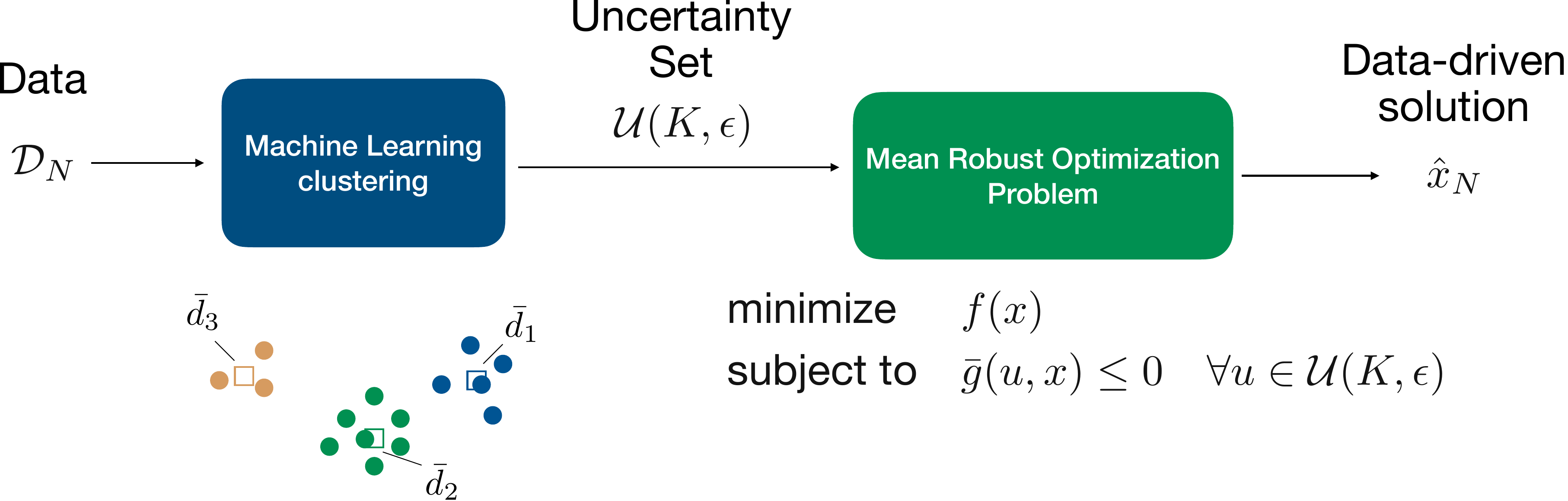} }}%
		\caption{Mean robust optimization procedure.}%
		\label{fig:mro_proc}%
	\end{figure}

\paragraph{Data-driven procedure.}%
\label{par:data_driven_procedure_}
Given the problem data, we formulate the uncertainty set from clustered data using machine learning, with the choice of $K$ and $\epsilon$ chosen experimentally.
Then, we solve the \gls{MRO} problem to arrive at a data-driven solution $\hat{x}_N$ which satisfies the probabilistic guarantee~\eqref{eq:prob_guarantees1}, see Figure~\ref{fig:mro_proc}.

\subsection{Solving the robust problem}
We now outline two ways to solve the~\gls{MRO} problem, using a direct convex reformulation and using a cutting plane algorithm. The reformulation and cutting plane procedure follow usual techniques for \gls{RO} problems in existing literature~\citep{ben-tal_robust_2009,kuhn2019wasserstein,robustadaptopt,cuttingplane}, with adaptations made for the \gls{MRO} setup, as well as a reformulation derived for the case $p = \infty$. We include simple examples and pseudocode for completeness.

\subsubsection{Direct convex reformulation for $p \geq 1$}
\label{sub:direct}
In the case where $p \geq 1$, the~\gls{MRO} can be rewritten as the optimization problem
\begin{equation}
	\label{eq:unif-const}
	\begin{array}{ll}
		\mbox{minimize} & f(x)\\
		\mbox{subject to} & \begin{dcases}
			\begin{rcases}
		\underset{v_1, \dots, v_K \in \supp}{\text{maximize}} &\quad \sum_{k=1} ^K w_k g(v_k, x) \\
		\mbox{subject to} & \quad \sum_{k=1} ^K w_k \| v_k - \bar{d}_k \|^p \le \epsilon^p
			\end{rcases}
			\end{dcases} \le 0,
	\end{array}
\end{equation}
which, by dualizing the inner maximization problem, has the following reformulation for any $\epsilon > 0$:
\begin{equation}
	\label{eq:robustopt_p}
	\begin{array}{ll}
		\mbox{minimize} & f(x)\\
		\mbox{subject to} & \sum_{k=1}^{K} w_k s_k  \le 0\\
		& [-g]^*(z_{k} - y_{k}, x) + \sigma_{{\supp}}(y_{k}) - z_{k}^T\bar{d}_k  + \phi(q)\lambda \left\|z_{k}/\lambda \right\|^q_*  +\lambda \epsilon^p \leq s_k\\
		&\hspace{6cm} \quad k = 1,\dots, K\\
		& \lambda \geq 0,
	\end{array}
\end{equation}
with variables $\lambda \in \reals$, $s_k \in \reals$, $z_{k} \in \reals^{m}$, and $y_{k} \in \reals^m$.
Here, $[-g]^*(z, x) \define \sup_{u\in \dom g} z^Tu - [-g(u, x)] $ is the conjugate of $-g$, $\sigma_ \supp(z) \define\sup_{u \in  \supp}  z^Tu$ is the support function of $ \supp \subseteq \mathbf{R}^m$, $\| \cdot \|_*$ is the dual norm of  $\| \cdot \|$, and $\phi(q) = (q-1)^{(q-1)}/q^q$ for $q > 1$~\citep[Theorem 8]{kuhn2019wasserstein}.
Note that $q$ satisfies $1/p + 1/q = 1$, \ie, $q = p/(p - 1)$.
\reviewChanges{When $p=1$ and $q =\infty$, we note the formulation in~\eqref{eq:robustp1}.}
The support function $\sigma_ \supp$ is also the conjugate of $\chi_ \supp$, which is defined $\chi_{\supp}(u) = 0$ if $u \in  \supp$, and $\infty$ otherwise.
The proof of the derivation and strong duality of the constraint is delayed to Appendix~\ref{app:dual_form}.
Since the dual of the constraint becomes a minimization problem, any feasible solution that with objective less than or equal to $0$ will satisfy the constraint, so we can remove the minimization to arrive at the above form.
While traditionally we take the supremum instead of maximizing, here the supremum is always achieved as we assume $g$ to be upper-semicontinuous.
For specific examples of the conjugate forms of different $g$, see Bertsimas and den Hertog~\cite[Section 2.5]{robustadaptopt} and Beck~\cite[Chapter 4]{amirbeck}.

When $K$ is set to be $N$, $w_k$ is $1/N$, and this is of an analogous form to the convex reduction of the worst case problem for Wasserstein \gls{DRO}, which we will introduce in Section~\ref{sub:connectionsdro}.

\reviewChanges{
We note the special case when $p = 1$. We observe from~\citep[Section 2.2 Remark 1]{kuhn2019wasserstein} that
\begin{equation*}
	 \lim_{q \rightarrow \infty}\phi(q)\lambda \left\|z_{k}/\lambda \right\|^q_* =
	\begin{dcases} 0 & \text{if}\; \|z_{k} \| \leq \lambda \\
		\infty & \text{otherwise}.
	\end{dcases}
	\end{equation*}
Therefore, the above formulation becomes
\begin{equation}
	\label{eq:robustp1}
	\begin{array}{ll}
		\mbox{minimize} & f(x)\\
		\mbox{subject to} & \sum_{k=1}^{K} w_k s_k  \le 0\\
		& [-g]^*(z_{k} - y_{k}, x) + \sigma_{{\supp}}(y_{k}) - z_{k}^T\bar{d}_k  + \lambda\epsilon^p\leq s_k\\
		&\hspace{6cm} \quad k = 1,\dots, K\\
		& \|z_{k} \| \leq \lambda, \quad k = 1,\dots, K\\
		& \lambda \geq 0.
	\end{array}
\end{equation}}

\paragraph{Example with affine constraints.}
Consider a single affine constraint of the form
\begin{equation}
	\label{eq:affineg}
	(a + Pu)^Tx \le b,
\end{equation}
where $a \in \reals^n$, $P \in \reals^{n \times m}$, and $b \in \reals$.
In other words, $g(u, x) = (a + Pu)^Tx - b$, and the support set is $\supp = \reals^m$.
Note that, in this case, $y_k$ must be $0$ for the support function $\sigma_{\supp}(y_k)$ to be finite.
We compute the conjugate as
\begin{equation}
\label{eq:conj_affine}
[-g]^*(z, x) =  \sup_u z^Tu +b - (a + Pu)^Tx =
\begin{dcases} a^Tx - b & \text{if}\; z + P^Tx = 0\\
	\infty & \text{otherwise}.
\end{dcases}
\end{equation}
To substitute $\sigma_{\supp}(y_k)$ and $[-g]^*(z_k - y_k, x)$ into~\eqref{eq:robustopt_p}, we note that $y_k = 0$ and $z_k = -P^Tx$, \ie, $z_k$ is independent from $k$.
By combining the $K$ constraints in~\eqref{eq:robustopt_p}, we arrive at the form
\begin{equation}
	\label{eq:example-affine-p}
	\begin{array}{ll}
		\mbox{minimize} & f(x)\\
		\mbox{subject to} & a^Tx - b + \phi(q)\lambda \left\|P^Tx/\lambda \right\|^q_*  +\lambda \epsilon^p + (P^Tx)^T\sum_{k=1}^{K} w_k \bar{d}_k \le 0\\
		& \lambda \ge 0,\\
	\end{array}
\end{equation}
where the number of variables or constraints does not depend on $K$.
Since vector $\sum_{k=1}^{K} w_k \bar{d}_k$ is the average of the datapoints in $\mathcal{D}_N$ for any $K \in \{1,\dots,N\}$, this formulation corresponds to always choosing $K=1$.

\subsubsection{Direct convex reformulation for $p = \infty$ }
\label{sub:direct-p}
In the case where ${p = \infty}$, the \gls{MRO} can be rewritten as the optimization problem
\begin{equation}
	\label{eq:unif-const-p}
	\begin{array}{ll}
		\mbox{minimize} & f(x)\\
		\mbox{subject to} & \begin{dcases}
			\begin{rcases}
		\underset{v_1, \dots, v_K \in\supp}{\text{maximize}} &\quad \sum_{k=1} ^K w_k g(v_k, x) \\
		\mbox{subject to} & \quad  \| v_k - \bar{d}_k \| \le \epsilon, \quad k = 1,\dots, K
			\end{rcases}
			\end{dcases} \le 0,
	\end{array}
\end{equation}
which has the following reformulation where the constraint above is dualized, for all $\epsilon > 0$:
	\begin{equation}
	\label{eq:unif-dual-p}
	\begin{array}{ll}
		\text{minimize} & f(x)\\
 \text{subject to} \quad & \sum_{k=1}^{K} w_k s_k \leq 0\\
&[-g]^*(z_{k} - y_{k},x) + \sigma_{{\supp}}(y_{k}) - z_{k}^T\bar{d}_k  + \lambda_k \epsilon \leq s_k\\
& \hspace{3cm} \quad k = 1,\dots, K\\
&  \left\|z_{k} \right\|_* \leq \lambda_k \quad k = 1,\dots, K,
\end{array}
\end{equation}
with new variables $s_k \in \reals$, $z_{k} \in \reals^m$, and $y_{k} \in \reals^m$. The proof is delayed to Appendix~\ref{app:dual_form-p}.

\begin{remark}[Case ${p = \infty}$ is the limit of case ${p \geq 1}$]
	\label{limit-proof}
	In terms of the primal problem,~\eqref{eq:unif-const-p} is the limiting case of~\eqref{eq:unif-const} as $p \rightarrow \infty$.
	In terms of the reformulated problem with dualized constraints, problem~\eqref{eq:unif-dual-p} is the limiting case of~\eqref{eq:robustopt_p} as $p \rightarrow \infty$ .
	The proofs are delayed to Appendix~\ref{app:primal_limit} and Appendix~\ref{app:dual_limit} respectively. These proofs extend the ideas stated in~\cite[Proposition 3]{wass-p-1}.
\end{remark}

\paragraph{Example with affine constraints.}%
\label{par:example}
Consider again the case of affine constraint as in~\eqref{eq:affineg} with support set $\supp = \reals^m$, now with $p=\infty$.
Following a similar derivation as~\eqref{eq:example-affine-p}, we substitute the conjugate function $[-g]^*$~\eqref{eq:conj_affine} in problem~\eqref{eq:unif-dual-p}, we can obtain
\begin{equation}
	\label{eq:example-affine-inf}
	\begin{array}{ll}
		\mbox{minimize} & f(x)\\
		\mbox{subject to} & a^Tx - b + \epsilon\left\|P^Tx \right\|_* + (P^Tx)^T\sum_{k=1}^{K} w_k \bar{d}_k \le 0,\\
	\end{array}
\end{equation}
where the number of constraints and variables does not depend on $K$.
Similarly to problem~\eqref{eq:example-affine-p}, the term $\sum_{k=1}^{K} w_k \bar{d}_k$ is the average of the datapoints in $\mathcal{D}_N$ for any $K \in \{1,\dots,N\}$. Therefore, the choice of $K$ does not affect this formulation.
\reviewChanges{This can be viewed as the robust counterpart when the uncertainty set is a norm ball of radius $\epsilon$ centered at $(1/N)\sum_{i=1}^N d_i$}

Note that, if $\bar{d} = 0$ the constraint can be simplified even further, obtaining ${a^Tx + \epsilon\|P^Tx\|_* \le b}$, which corresponds to the robust counterpart in \gls{RO} with norm uncertainty sets~\citep[Section 2.3]{robustadaptopt},~\citep[Chapter 2]{ben-tal_robust_2009}.

\begin{remark}
	When $g$ is affine and \reviewChanges{$\supp = \reals^m$}, for any $\epsilon$ and norm, the convex reformulations for $p = 1$ and $p = \infty$ are identical. The proof appears in Appendix~\ref{app:equivalence}.
\end{remark}

\subsubsection{Cutting plane algorithm}
The second approach to solve problem~\eqref{eq:robustopt} is to use a cutting plane procedure, in which we consider the minimization problem where $x$ is the variable and $S$ a finite set of values for the uncertainty,
\begin{equation}
	\label{eq:mro_min}
	\begin{array}{ll}
		\mbox{minimize} & f(x)\\
		\mbox{subject to} & \bar{g}(u, x)  \le 0, \quad \forall u \in
		\hat{S},
		\end{array}
\end{equation}
and the maximization problem over $u$ with $x^k$ fixed,
\begin{equation}
	\label{eq:mro_max}
	\begin{array}{ll}
		\mbox{maximize} & \bar{g}(u, x^k)\\
		\mbox{subject to} & u \in \uncset(K, \epsilon)
	\end{array}
\end{equation}
The procedure works as follows.
We first solve~\eqref{eq:mro_min} with a set $\hat{S} = \{\bar{u}\}$, where $\bar{u}$ is nominal value of the uncertainty, obtaining $x^k$.
Then, we solve~\eqref{eq:mro_max}, obtaining $u^k$.
If $\bar{g}(u^k, x^k) > 0$, then we add $u^k$ to the set $\hat{S}$.
Otherwise, we terminate. This procedure is summarized in Algorithm~\ref{alg:cutting_plane}.
As demonstrated by Bertsimas et al.~\cite{cuttingplane}, the cutting plane and convex reformulation methods are comparable in terms of performance, thus both are viable.
\begin{algorithm}
	\caption{Cutting plane algorithm to solve~\eqref{eq:robustopt}}
	\label{alg:cutting_plane}
	\begin{algorithmic}[1]
	 \State {\bf given} $\hat{S} = \{\bar{u}\}$
	 \For{$k=1,\dots,k_{\rm max}$}
	 \State $x^k \gets $ solve minimization problem~\eqref{eq:mro_min} over $x$
	 \State $u^k \gets $ solve maximization problem~\eqref{eq:mro_max} over $u$
         \If{$\bar{g}(u^{k}, x^{k}) > 0$}
	 \State $\hat{S} \gets \hat{S} \cup \{u^k\}$
	 \Else
	 \State {\bf return} $x^k$
	 \EndIf
	 \EndFor
	\end{algorithmic}
\end{algorithm}

\subsection{\reviewChanges{Maximum-of-concave constraint function}}
\label{sec:max_concave}
\reviewChanges{We now consider a more general maximum-of-concave function 
$$g(u, x) = \max_{j\leq J}g_j(u, x),$$ with each $-g_j$ being proper, convex, and lower-semicontinuous in~$u$ for all $x$. When we take $J=1$, we arrive back at the formulations given in Section~\ref{sec:mro}.
Note that any problem with multiple uncertain constraints ${g_j(u,x),\; j = 1,\hdots,J}$, where we assume the usual conditions on $g_j$, can be combined to create a joint constraint of this maximum-of-concave form. 
As mentioned in Section~\ref{sec:intro}, this can also be used to model $\CVaR$ constraints, which has a maximum-of-concave analytical form.}

\paragraph{\reviewChanges{Problem parametrization.}}
\reviewChanges{We now consider constraints of the form
\begin{equation}
	\label{eq:constraints_max}
\bar{g}(u, x) \define  \sum_{k=1}^K\sum_{j=1}^J \alpha_{jk} g_j(v_{jk}, x),
\end{equation}
where $\alpha\in \Gamma $, with $\Gamma = \{\alpha~|~\sum_{j=1}^J\alpha_{jk}= w_k, \alpha_{jk}\geq 0 ~\forall k,j\}$. For each constituent function $g_j$, the uncertainty set contains a set of vectors $(v_{j1},\dots,v_{jK})$, and a set of parameters $(\alpha_{j1},\dots, \alpha_{jK})$ to denote the fraction of mass assigned to that function for each $k$. The total amount of mass assigned for each cluster, $\sum_{j=1}^J\alpha_{jk}$, is the weight of the cluster, $w_k$. }

\reviewChanges{We use a summation over weighted pieces $g_j$ instead of a maximum over $g_j$, as this is a generalization of the maximum, and has a more natural dual reformulation. We take inspiration from~\citep[Section 4.2]{mohajerin_esfahani_data-driven_2018}, where $\alpha$ arises from the extremal distribution for Wasserstein \gls{DRO}.
Note that the intuitive maximization over $g_j$'s is analogous to setting $\alpha_{jk} = w_k$ for a specific $j$ for each $k$, and $\alpha_{jk} = 0$ otherwise.}

\reviewChanges{The uncertainty set is given as follows.}
\reviewChanges{
\paragraph{Case $p \geq 1$.}
In the case where $p \geq 1$, we have
\begin{equation*}
	\label{eq:uncset_p_max}
	\uncset(K, \epsilon) \define \left\{ u = (v_{11},\dots,v_{JK}) \in \supp^{K\times J} \quad\middle|\quad \sum_{k=1}^K \sum_{j=1}^J \alpha_{jk} \| v_{jk} - \bar{d}_k \|^p \le \epsilon^p,\alpha \in \Gamma \right\}.
\end{equation*}
Note that the single concave case given previously follows when we take $J = 1$. All parameters are defined as in the single concave case.}

\reviewChanges{\paragraph{Case $p = \infty$.}
In the case where $p = \infty$, the set we consider becomes
\begin{equation*}
	\label{eq:uncset_inf_max}
	\uncset(K, \epsilon) \define \left\{ u = (v_{11},\dots,v_{JK}) \in \supp^{K\times J} \quad\middle|\quad \max_{k=1,\dots,K} \sum_{j=1}^J (\alpha_{jk}/w_k) \| v_{jk} - \bar{d}_k \| \le \epsilon, \alpha \in \Gamma\right\},
\end{equation*}
where we once again introduce weight parameters $\alpha$. 
}

Following these changes, $\hat{x}_N$ is again the solution to the robust optimization problem~\eqref{eq:robustopt}, defined now with the generalized uncertainty set and constraint.

\paragraph{Solving the robust problem.}
We give the direct reformulation approach for solving the generalized problem for $p \geq 1$. The case $p = \infty$ is delayed to Appendix~\ref{app:dual_form-p_max}.
We write the~\gls{MRO} problem as the optimization problem
\begin{equation}
	\label{eq:unif-const_max}
	\begin{array}{ll}
		\mbox{minimize} & f(x)\\
		\mbox{subject to} & \begin{dcases}
			\begin{rcases}
		\underset{v_{11},\dots,v_{JK} \in \supp, \alpha \in \Gamma}{\text{maximize}} &\quad \sum_{k=1} ^K \sum_{j=1}^J \alpha_{jk} g_j(v_{jk}, x) \\
		\mbox{~~~~subject to} & \quad \sum_{k=1}^K \sum_{j=1}^J \alpha_{jk} \| v_{jk} - \bar{d}_k \|^p \le \epsilon^p\\
			\end{rcases}
			\end{dcases} \le 0,
	\end{array}
\end{equation}
and, by dualizing the inner maximization problem, arrive at the reformulation for all $\epsilon > 0$:
\begin{equation}
	\label{eq:robustopt_p_max}
	\begin{array}{ll}
		\mbox{minimize} & f(x)\\
		\mbox{subject to} & \sum_{k=1}^{K} w_k s_k  \le 0\\
		& [-g_j]^*(z_{jk} - y_{jk}, x) + \sigma_{{\supp}}(y_{jk}) - z_{jk}^T\bar{d}_k  + \phi(q)\lambda \left\|z_{jk}/\lambda \right\|^q_*  +\lambda \epsilon^p \leq s_k\\
		&\hspace{4cm} \quad k = 1,\dots, K, \quad j =1, \dots, J\\
		& \lambda \geq 0,
	\end{array}
\end{equation}
with variables $\lambda \in \reals$, $s_k \in \reals$, $z_{jk} \in \reals^{m}$, and $y_{jk} \in \reals^m$. The proof is delayed to Appendix~\ref{app:dual_form_max}.
Again, while traditionally we take the supremum instead of maximizing, here the supremum is always achieved as we assume \reviewChanges{$g_j$ to be upper-semicontinuous for all $j$.}

\reviewChanges{In addition, when $K$ is set to be $N$, and $w_k$'s are $1/N$, this is also of an analogous form to the convex reduction of the worst case problem for Wasserstein \gls{DRO}, given in Section~\ref{sub:connectionsdro}. }

\section{Links to Wasserstein \glsentrylong{DRO}}
\label{sub:connectionsdro}
\Glsentryfull{DRO} solves the problem
\begin{equation}
	\label{eq:dro_prob}
	\begin{array}{ll}
		\mbox{minimize} & f(x)\\
		\mbox{subject to} &   \sup_{\mathbf{Q}\in \mathcal{P}_N} \Expect^\mathbf{Q}(g(u,x))  \le 0,
	\end{array}
	\end{equation}
where the ambiguity set $\mathcal{P}_N$ contains, with high confidence, all distributions that could have generated the training samples $\mathcal{D}^N$, such that the probabilistic guarantee~\eqref{eq:prob_guarantees1} is satisfied.
Wasserstein \gls{DRO} constructs $\mathcal{P}_N$ as a ball of radius $\epsilon$ with respect to the Wasserstein metric around the empirical distribution ${\hat{\prob}^N \define \sum_{i=1}^N \delta_{d_i}/N}$, where $\delta_{d_i}$ denotes the Dirac distribution concentrating unit mass at $d_i \in \mathbf{R}^m$.
Specifically, we write
\begin{equation*}
	\label{eq:Wball}
	\mathcal{P}_N = \mathbf{B}_\epsilon^p (\hat{\prob}^N) \define \{\mathbf{Q} \in \mathcal{M}(\supp) \mid W_p(\hat{\prob}^N,\mathbf{Q}) \leq \epsilon \},
\end{equation*}
where $\mathcal{M}(\supp)$ is the set of probability distributions $\mathbf{Q}$ supported on $\supp$ with bounded $p$-th moments, \ie, $\int_S \|u\|^p \, \mathrm{d}\mathbf{Q}(u) < \infty$, and
$$W_p(\mathbf{Q},\mathbf{Q'}) \define \inf \left\{ \left(\int_\reviewChanges{\supp} \|u - u'\|^p \Pi (\dif u,\dif u') \right)^{1/p} \right\}.$$
Here, $p$ is any integer greater than 1, and $\Pi$ is any joint distribution of $u$ and $u'$ with marginals $\mathbf{Q}$ and $\mathbf{Q'}$.

When $K=N$, the constraint of the \gls{DRO} problem~\eqref{eq:dro_prob} is equivalent to the constraint of~\eqref{eq:robustopt}.
In particular, for case $p \geq 1$, the expression
\begin{equation}
\label{eq:dro-reduct}
\begin{aligned}
    \sup_{\mathbf{Q}\in \mathbf{B}_\epsilon^p(\hat{\prob}^N)} \Expect^\mathbf{Q}(g(u,x)),
	\end{aligned}
\end{equation}
is equivalent to the dual of the constraint of~\eqref{eq:unif-const_max}, when $K = N$, and $w_k = 1/N$. 
This is noted in~\citep{kuhn2019wasserstein,zhen2021mathematical}.
We give a proof of strong duality in Appendix~\ref{app:convex_red}.
This is the dual of the generalized max-of-concave form, which is equivalent to the dual of the single concave form~\eqref{eq:unif-const} when $J = 1$.
By the same logic, in the case where $p = \infty$, the expression is equivalent to the dual of the constraint of~\eqref{eq:unif-const-p}.
Given the above reductions, we can rewrite the Wasserstein \gls{DRO} problem in the same form as~\eqref{eq:unif-const_max}, the \gls{MRO} problem.

Our approach can then be viewed as a form of Wasserstein \gls{DRO}, with the difference that, when $K<N$, we deal with the clustered and averaged dataset.
We form $\mathcal{P}_N$ as a ball around the empirical distribution $\hat{\prob}^K$ of the centroids of our clustered data
\begin{equation*}
\hat{\prob}^K \define \sum_{k=1}^K w_k\delta_{\bar{d}_k},
\end{equation*}
where $w_k$ is the proportion of data in cluster $k$.
This formulation allows for the reduction of the sample size while preserving key properties of the sample, which translates directly to a reduction in the number of constraints and variables, while maintaining high quality solutions.

\subsection{Satisfying the probabilistic guarantees}
\label{par:gau}
As we have noted the parallels between \gls{MRO} and Wasserstein \gls{DRO}, we now show that the conditions for satisfying the probabilistic guarantees are also analogous.
\paragraph{Case $p \geq 1$.} Wasserstain \gls{DRO} satisfies~\eqref{eq:prob_guarantees1} if the data-generating distribution, supported on a convex and closed set $\supp$, satisfies a {\it light-tailed assumption}~\citep{fournier_guillin_2013,mohajerin_esfahani_data-driven_2018}: there exists an exponent $a > 0$ and $t > 0 $ such that ${A = \Expect^{\prob}(\exp (t \|u \|^a))  = \int_{\reviewChanges{\supp}} \exp (t \| u \| ^a) \prob(du) < \infty.}$
We refer to the following theorem.
\begin{theorem}[Measure concentration {\citep[Theorem 2]{fournier_guillin_2013}}]
\label{theorem_ball}
If the light-tailed assumption holds, we have
	\begin{equation*}
	\prob^N(W_p(\prob,\hat{\prob}^N) \geq \epsilon) \le \phi(p, N, \epsilon),
	\end{equation*}
where $\phi$ is an exponentially decaying function of $N$.
\end{theorem}
Theorem~\eqref{theorem_ball} estimates the probability that the unknown data-generating distribution $\prob$ lies outside the Wasserstein ball $\mathbf{B}_\epsilon^p (\hat{\prob}^N)$, which is our ambiguity set.
Thus, we can estimate the smallest radius $\epsilon$ such that the Wasserstein ball contains the true distribution with probability $1 - \beta$, for some target $\beta \in (0,1)$.
We equate the right-hand-side to $\beta$, and solve for $\epsilon_N(\beta)$ that provides us the desired guarantees for Wasserstein \gls{DRO}~\citep[Theorem 3.5]{mohajerin_esfahani_data-driven_2018}.

\paragraph{Case $p = \infty$.} When $p = \infty$, Bertsimas et al.~\cite[Section 6]{wass-p-1} note that the light-tailed assumption is no longer sufficient.
Wasserstein \gls {DRO} satisfies~\eqref{eq:prob_guarantees1} under stronger assumptions, as given in the following theorem.
\begin{theorem}[Measure concentration, $\mathbf{p = \infty}$ {\cite[Theorem 1.1]{wass-p-guarantee}}]
\label{theorem-inf}
	Let the support $\supp \subset \reals^m$ of the data-generating distribution be a bounded, connected, open set with Lipschitz boundary. Let $\prob$ be a probability measure on  $\supp$ with density $\rho : \supp \rightarrow(0,\infty)$, such that there exists $\lambda \geq 1$ for which ${1/\lambda \leq \rho(x) \leq \lambda, \hspace{2mm} \forall x \in \supp}$. Then,
	\begin{equation*}
		\prob^N(W_\infty(\prob,\hat{\prob}^N) \geq \epsilon) \le \phi(N, \epsilon),
		\end{equation*}
	where $\phi$ is an exponentially decaying function of $N$.
\end{theorem}
We can again equate the right-hand-side to $\beta$ and find  $\epsilon_N(\beta)$. We extend this result to the clustered set in \gls{MRO}.

\begin{theorem}[\gls{MRO} finite sample guarantee] \label{theorem_mro} Assume the light-tailed assumption holds when $p \geq 1$, and the corresponding assumptions hold when $p = \infty$. If ${\beta \in (0,1)}$, $\eta_N(K)$ is the average $p$-th powered distance of data-points in $\mathcal{D}_N$ from their assigned cluster centers, and $\hat{x}_N$ is the optimal solution to~\eqref{eq:robustopt} with uncertainty set $\uncset(K,\epsilon_N(\beta) + \eta_N(K)^{1/p})$, then the finite sample guarantee~\eqref{eq:prob_guarantees1} holds.
\end{theorem}

\begin{proof}
Compared with Wasserstein \gls{DRO}, \gls{MRO} has to account for the additional difference between the two empirical distributions $\hat{\prob}^N$ and $\hat{\prob}^K$.
We can write
\begin{equation*}
	\begin{array}{ll}
		&\displaystyle \hat{\prob}^N \define \sum_{i=1}^N \frac{1}{N}\delta_{d_i} = \sum_{k=1}^K \sum_{i \in C_k} \frac{|C_k|}{N}\frac{1}{|C_k|}   \delta_{d_i}, \\ 
		&\displaystyle \hat{\prob}^K = \sum_{k=1}^K w_k\delta_{\bar{d}_k} = \sum_{i=1}^K \frac{|C_k|}{N}\delta_{\bar{d}_k}. 
		\end{array}
	\end{equation*}
If we introduce a new parameter, $\eta_N(K)$, defined as
\begin{equation*}
	\eta_N(K) =  \frac{1}{N} \sum_{i = 1}^K \sum_{i \in C_k} \|d_i - \bar{d}_k\|^p
	\end{equation*}
	 the average $p$-powered distance with respect to the norm used in the Wasserstein metric, of all data-points in $\mathcal{D}_N$ from their assigned cluster centers $\bar{d}_k$, we notice that
\begin{equation*}
	\begin{aligned}
	W_p(\hat{\prob}^K,\hat{\prob}^N)^p &= \inf_\Pi \left\{ \int_\reviewChanges{\supp} \|u - u'\|^p \Pi (\dif u,\dif u') \right\} \quad \text{($\Pi$ any joint distribution of $\hat{\prob}^K$, $\hat{\prob}^N$  )}\\
&\leq \sum_{i = 1}^K \frac{|C_k|}{N}  \int_\reviewChanges{\supp} \|u - \bar{d}_k\|^p \reviewChanges{\hat{\prob}^N(u|u' =\bar{d}_k)}({\rm d}u)\\
&\leq \sum_{i = 1}^K \frac{|C_k|}{N}  \frac{1}{|C_k|} \sum_{i \in C_k} \|d_i - \bar{d}_k\|^p\\
&= \eta_N(K),
	\end{aligned}
	\end{equation*}
where we have replaced the integral with a finite sum, as the distributions are discrete. Therefore, by Theorems~\ref{theorem_ball},~\ref{theorem-inf} and the \reviewChanges{triangle inequality for the Wasserstein metric~\cite{wass_tri}},
\begin{equation*}
	\begin{aligned}
	W_p(\prob,\hat{\prob}^K) &\leq W_p(\prob,\hat{\prob}^N) + W_p(\hat{\prob}^K,\hat{\prob}^N)\\
	&\leq \epsilon_N(\beta) + \eta_N(K)^{1/p},
\end{aligned}
	\end{equation*}
with probability at least $1 - \beta$. We thus have
\begin{equation*}
		\prob(\prob \in \mathbf{B}_{\epsilon_N(\beta) + \eta_N(K)^{1/p}}^p (\hat{\prob}^K) )\geq 1 - \beta,
	\end{equation*}
which implies the uncertainty set $\uncset(K, \epsilon_N(\beta) + \eta_N(K)^{1/p})$ contains all possible realizations of uncertainty with probability $1-\beta$, so the finite sample guarantee~\eqref{eq:prob_guarantees1} holds.
\end{proof}

This result introduces a tradeoff between conservatism (affected by the radius) and computational efficiency (affected by the sample size), controlled by the value of $K$. When a large $K$ is chosen, the effect on the radius is small, but the computational savings may also be minimal. The choice of $K$ is thus paramount.

One might note that a blanket increase in the radius regardless of the optimization problem may introduce unnecessary conservatism. This brings us to the next section, where we bring the constraint functions $g$ into our considerations. 

\section{Worst-case value of the uncertain constraint}
\label{sec:conservatism}
In the previous section, we proposed a theoretical increase in $\epsilon$ to maintain the same finite sample guarantee before and after clustering. 
However, a question remains: what is the extent of the effects of clustering if we {\it don't} increase $\epsilon$?
In this section, we thus approach the analysis in a different manner: keeping $\epsilon$ constant, we quantify the change in the {\it worst-case value of the constraint function} that arises from clustering. 
In fact, for select cases, our results suggest there is no need to increase $\epsilon$ after clustering; under specific curvature conditions on the  constraint function $g$, we may obtain a more conservative, or even unchanged solution after clustering, in which case the original finite sample guarantee is retained. 

We begin with a remark on the clustering value attained. The \gls{MRO} approach is closely centered around the concept of clustering to reduce sample size while maintaining sample diversity.
We wish to cluster points that are close together, such that the objective is only minimally affected. With this goal, we would like to minimize the average distance of the points in each cluster to their data-center,
\begin{equation*}
	D(K)^\star =\text{minimize } D(K) =  \text{minimize } \frac{1}{N}\sum_{k=1}^K \sum_{d_i \in C_k} \| d_i - \bar{d}_k \|^2_2,
\end{equation*}
	where $\bar{d}_k$ is the mean of the points in cluster $C_k$.
	While the best performance is attained with $D(K)^\star$, in practice we work with the approximation $D(K)$, where $C_k$ is decided by a clustering algorithm. This value upper bounds $D(K)^\star$. A well-known algorithm is $K$-means~\citep{kmeans}, where we create $K$ clusters by iteratively solving a least-squares problem.
From here on we use only $D(K)$, and note that for the case $p=2$, we have $D(K)=\eta_N(K)$ from Theorem~\ref{theorem_mro}.

In this section, we then show the effects of clustering on the {\it worst-case value of the constraint function in~\eqref{eq:robustopt}}.
\reviewChanges{We prove two sets of results, corresponding to $g$ given as a single concave function, and as a more general maximum-of-concave function. For the latter, we also include the special case of the maximum-of-affine function.}

\subsection{\reviewChanges{Single concave function}}
\label{sub:single-con}
For the simplest case of a single concave function, we prove that when the support is large enough,
\begin{itemize}
	\item If $g$ is affine in $u$, \gls{MRO} does not increase the worst-case value, regardless of $K$.
	\item If $g$ is concave in $u$ and satisfies certain smoothness conditions, \gls{MRO} has a higher worst-case value than Wasserstein \gls{DRO} and the increase is inversely related to the number of clusters~$K$. In other words, the smaller the $K$, the higher the worst-case value.
\end{itemize}

\paragraph{Quantifying the clustering effect.}%
\label{par:quantify_clustering_effect_}
To quantify the effect of clustering, we calculate the difference between the following formulations of the worst-case value of the constraint in~\eqref{eq:robustopt}
\begin{equation}
	\tag{MRO-N}
	\begin{aligned}
	\label{eq:nocluster_one}
	\bar{g}^N(x) = \underset{v_1\dots v_N}{\text{maximize}} &\quad\frac{1}{N} \sum_{i=1} ^N g(v_i, x) \\
	\mbox{subject to} & \quad \frac{1}{N} \sum_{i=1} ^N  \| v_i - d_i \|^p \le \epsilon^p\\
	& \quad v_i\in \supp \quad i = 1,\dots,N,
	\end{aligned}
\end{equation}
\begin{equation}
	\tag{MRO-K}
	\label{eq:cluster_one}
	\begin{aligned}
		\bar{g}^K(x) = \underset{u_1\dots u_K}{\text{maximize}} &\quad  \sum_{k=1} ^K  \frac{|C_k|}{N} g(u_k, x) \\
	\mbox{subject to} & \quad\sum_{k=1} ^K \frac{|C_k|}{N}  \| u_k - \bar{d}_k \|^p \le \epsilon^p\\
	& \quad u_k\in \supp \quad k = 1,\dots,K,
	\end{aligned}
\end{equation}
and
\begin{equation}
	\tag{MRO-N*}
	\label{eq:nocluster_inf}
	\begin{aligned}
		\bar{g}^{N*}(x) = \underset{v_1\dots v_N}{\text{maximize}} &\quad\frac{1}{N} \sum_{i=1} ^N g(v_i, x) \\
		\mbox{subject to} & \quad \frac{1}{N} \sum_{i=1} ^N  \| v_i - d_i \|^p \le \epsilon^p,
		\end{aligned}
\end{equation}
where~\eqref{eq:nocluster_one} is the formulation of the constraint without clustering, akin to traditional Wasserstein \gls{DRO},~\eqref{eq:nocluster_inf} is the same, except we drop the support constraint, and~\eqref{eq:cluster_one} is the formulation with $K$ clusters.
\reviewChanges{
From here on, when we mention that the support {\it affects the worst-case constraint value}, we refer to situations where at least one of the constraints $v_i \in \supp$ for $i = 1,\dots,N$ is binding.
Formally, the definition is $\bar{g}^N(x) \neq \bar{g}^{N*}(x)$ for any $x$ feasible for the \gls{DRO} problem. 
We note a sufficient but not necessary condition for the support to not affect the worst-case constraint value: the situation in which the support doesn't affect the uncertainty set, which is defined as 
$$\left \{u \in \reals^{N\times m}: (1/N)\sum_{i=1}^N  \| u_i - d_i \| \le \epsilon\right \} = \left \{ u \in \supp^{N\times m}: (1/N)\sum_{i=1}^N \| u_i - d_i \| \le \epsilon\right \}.$$
If the support satisfies this condition, then we can conclude that $\bar{g}^N(x) = \bar{g}^{N*}(x)$ for any $x$ feasible for the \gls{DRO} problem, and obtain improved bounds below. 
While the condition depends on the location of the datapoints, it is acceptable to have this dependency, as this is a condition we can check given data to potentially improve the following bounds, without having to solve the MRO problem.}

With these definitions, we can construct solutions for~\eqref{eq:nocluster_one},~\eqref{eq:cluster_one}, and~\eqref{eq:nocluster_inf} to prove the following relations.
\begin{theorem}
	\label{thm:1}
	\reviewChanges{
	With the same $x$ and $\epsilon$, and for any integer $p \geq 1$, we always have $$\bar{g}^{N}(x) \le \bar{g}^K(x).$$
Suppose that Assumption~\ref{assp:dom} holds, and $g$ satisfies an $L$-smooth condition on its domain with respect to the $\ell_2$-norm and for a given $x$,
	\begin{equation*}
		\left\|\nabla g(v,x) - \nabla g(u,x)\right\|_2 \leq L\| u - v\|_2.
	\end{equation*}
Then, with the same $x$ and $\epsilon$, and for any integer $p \geq 1$, we always have
	\begin{equation*}
	  \bar{g}^K(x) \leq \bar{g}^{N*}(x) + (L/2) D(K).
	\end{equation*}}
\end{theorem}
The proof is delayed to Appendix~\ref{app:thm1}. The results also hold for $p = \infty$, as we have shown in Remark~\ref{limit-proof} that the case $p = \infty$ is the limit of the case $p \geq 1$, and these results hold under the limit.

\reviewChanges{Let $\Delta$ be the maximum difference in constraint value resultant from relaxing the support constraint on the \gls{MRO} uncertainty sets, \ie, ${\Delta = \max_{x\in \mathcal{X}}\left(\bar{g}^{N*}(x) - \bar{g}^{N}(x)\right)}$, subject to $x$ being feasible for problem~\eqref{eq:robustopt}. 
As we assume Assumption~\ref{assp:dom} to hold, combined with the smoothness of $g$, we note that when solving for $\bar{g}^{N*}(x)$, the chosen $v_i$ values without the support constraint will still remain in the domain $\dom$ of $g$. 
Refer to a similar argument in Appendix~\ref{app:thm1} (ii) for details. 
The function $\bar{g}^{N*}(x) - \bar{g}^{N}(x)$ is then continuous in $x$ and everywhere defined for $x \in \mathcal{X}$, thus maximizing with respect to $\mathcal{X}$, a compact set, the value $\Delta$ is finite. }
Then, we observe that $ \bar{g}^K(x)  - \bar{g}^N(x)\leq \Delta +(L/2) D(K)$ for all such $x$, so the smaller the $D(K)$, (\ie, higher-quality clustering procedure), the smaller the increase in the worst-case constraint value. \reviewChanges{In addition, the value $\Delta$ is independent of $K$, as we calculate it with only $\bar{g}^{N*}(x)$ and $\bar{g}^{N}(x)$.}

\begin{remark} \label{rem:support}
While $\Delta$ could be constructed to be arbitrarily bad, in practice, we expect our relevant range of $\epsilon$ to be small enough such that the difference is insignificant.
We can then approximate $\Delta \approx 0$ and simply use the upper bound $(L/2)D(K)$, as this bound is often not tight. See Sections~\ref{sub:facility} and~\ref{sub:capitalbudget} for examples.
\end{remark}

\paragraph{Uncertain objective.}%
\label{par:uncertain_objective_}
When the uncertainty is in the objective, Theorem~\ref{thm:1} quantifies the difference in optimal values.
\begin{corollary}
	\label{cor:obj_unc}
Consider the problem where $g$ is itself the objective function we would like to minimize and $X \subseteq \reals^n$ represents the constraints, which are deterministic.
Then, $(L/2)D(K)$ + $\Delta$ upper bounds the difference in optimal values of the \gls{MRO} problem with $K$ and $N$ clusters.
\end{corollary}



\paragraph{Uncertain constraints.}%
\label{par:uncertain_constraints_}
When the uncertainty is in the constraints, the difference between $\bar{g}^K(x)$ and $\bar{g}^K(x)$ no longer directly reflects the difference in optimal values.
Instead, clustering creates a restriction on the feasible set for $x$ as follows.
For the same $\hat{x}$, $\bar{g}^K(\hat{x})$ takes a greater value than $\bar{g}^N(\hat{x})$.
Since both of them are constrained to be nonpositive from~\eqref{eq:robustopt}, the feasible region with $K$ clusters is smaller.


\paragraph{Affine dependence on uncertainty.}%
As a special case, when $g$ is affine in $u$, $L = 0$, so we observe the following corollary.
\begin{corollary}[Clustering with affine dependence on the uncertainty]
	\label{cor:affine}
	\reviewChanges{
If $g(u,x)$ is affine in $u$ and the worst-case constraint value is not affected by the support constraint, then clustering makes no difference to the optimal value and optimal solution to~\eqref{eq:robustopt}.}
\end{corollary}
\begin{proof}
In view of the primal problem and constraints, from Theorem~\ref{thm:1}, if $g(u,x)$ is affine in $u$ and the support does not affect the uncertainty set, $\bar{g}^N(x) = \bar{g}^K(x)$.
So for some fixed $\hat{x}$ we have ${\bar{g}^K(\hat{x}) \leq 0 \iff \bar{g}^N(\hat{x}) \leq 0}$. Therefore,
\[
    \hat{x} \text{ is feasible to~\eqref{eq:robustopt} for $K = N$} \iff \hat{x} \text{ is feasible to~\eqref{eq:robustopt} for $K < N.$}
\]
The feasible region of~\eqref{eq:robustopt} is identical for $K=N$ and $K < N$, and the optimal solutions will be identical so long as the optimal solution to~\eqref{eq:robustopt} is unique.
In view of the dual problem and constraints, if $g(u,x)$ is affine in $u$ following~\eqref{eq:affineg}, we observe from~\eqref{eq:example-affine-p} that the only term dependent on $K$ is  $(P^Tx)^T\sum_{k=1}^{K} w_k \bar{d}_k$, which is equivalent for all $K$.
\end{proof}

\reviewChanges{
\subsection{Maximum-of-concave functions}
\label{sub:maxofconcavebounds}
We now consider the more general case of a maximum-of-concave constraint function, $g(u,x) = \max_{j\leq J} g_j(u,x)$, subject to a polyhedral support, $\supp = \{u \mid Hu\leq h\}$. We define the new primal problems
\begin{equation}
	\tag{MRO-K-MAX}
	\bar{g}^K(x) = \begin{aligned}[t]
	\label{eq:nocluster_one_max}
	\underset{v_{11},\dots,v_{JK} \in \supp, \alpha \in \Gamma}{\text{maximize}} &\quad \sum_{k=1} ^K \sum_{j=1}^J \alpha_{jk} g_j(v_{jk}, x) \\
	\mbox{subject to} & \quad \sum_{k=1}^K \sum_{j=1}^J \alpha_{jk} \| v_{jk} - \bar{d}_k \|^p \le \epsilon^p\\
	& \quad v_i\in \supp \quad i = 1,\dots,K,
	\end{aligned}
\end{equation}
\begin{equation}
	\tag{MRO-N*-MAX}
	\label{eq:nocluster_inf_max}
	\bar{g}^{N*}(x) = \begin{aligned}[t]
		\underset{v_{11},\dots,v_{JN}, \alpha \in \Gamma}{\text{maximize}} &\quad \sum_{i=1} ^N \sum_{j=1}^J \alpha_{ji} g_j(v_{ji}, x) \\
		\mbox{subject to} & \quad \sum_{i=1}^N \sum_{j=1}^J \alpha_{ji} \| v_{ji} - d_i \|^p \le \epsilon^p\\
		\end{aligned}
\end{equation}
We also make use of the dual versions of the optimization problems, defined as follows. 
\begin{equation}
	\tag{MRO-N-Dual}
	\bar{g}^N(x) =\begin{aligned}[t]
	\label{eq:nocluster_dual}
	\underset{\lambda \geq 0, z_{ji},y_{ji}, s_i}{\text{minimize}} & \quad(1/N)\sum_i^N  s_i\\
	\mbox{subject to} & \quad[-g_j]^*(z_{ji} - H^T\gamma_{ji}) + \gamma_{ji}^T(h-Hd_i) - z_{ji}^Td_i  + \phi(q)\lambda \left\|z_{ji}/\lambda \right\|^q_* \\	
	& \hspace{3cm} + \lambda \epsilon^p \le s_i, \quad i = 1,\dots,N, \quad j = 1,\dots,J,
\end{aligned}
\end{equation}
where no clustering occurs, and 
\begin{equation}
	\tag{MRO-K-Dual}
	\label{eq:cluster_dual}
	\begin{aligned}
		\bar{g}^K(x) = \underset{\lambda \geq 0, z_{jk},y_{jk}, s_k}{\text{minimize}} &\quad \sum_k^K  (|C_k|/N) s_k\\
	\mbox{subject to} & \quad [-g_j]^*(z_{jk} - H^T\gamma_{jk})+ \gamma_{jk}^T(h-H\bar{d}_k) - z_{jk}^T\bar{d}_k + \lambda \epsilon^p   \\
	& \hspace{1cm} + \phi(q)\lambda \left\|z_{jk}/\lambda \right\|^q_*\le s_k, \quad k = 1,\dots,K, \quad j = 1,\dots,J,
\end{aligned}
\end{equation}
where we have $K$ clusters.
Given these definitions, we obtain bounds on the worst-case value of the constraint function for $K$ clusters.
\begin{theorem}
	\label{thm:maxconcave}
	When $g$ is the maximum of concave functions with domain ${\dom{g}= \reals^m}$ and polyhedral support $\supp = \{u \mid Hu\leq h\}$, and where each $g_j$ satisfies an $L$-smooth condition on its domain with respect to the $\ell_2$-norm at a given $x$, we have, for the same $x$ and $\epsilon$,
	\begin{equation*}
		\bar{g}^{N}(x) -\delta(K,z,\gamma) \leq \bar{g}^K(x) \leq  \bar{g}^{N*}(x) + \max_{j\leq J} (L_j/2)D(K),
	\end{equation*}
	where $\delta(K,z,\gamma) \define (1/N)\sum_{k=1}^K \sum_{i \in C_k}  \max_{j\leq J}((-z_{jk} -H^T\gamma_{jk})^T(d_i - \bar{d}_k))$, and $z$, $\gamma$ are the dual variable from~\eqref{eq:cluster_one}. The constants $L_j$ are the $L$-smoothness constants for the concave functions $g_j$.
	\end{theorem}
	The proof is delayed to Appendix~\ref{app:maxconcave}. We note that due to the nonconvex and nonconcave nature of maximum-of-concave functions, we can no longer directly fix the relationship between $\bar{g}^{N}(x)$ and $\bar{g}^{K}(x)$. Instead, we need to define the lower bound with the extra term $\delta(K,z,\gamma)$.
	However, in the special case where $g$ is a maximum-of-affine function, which is convex, we know $\bar{g}^{N*}(x)$ to be an upper bound on $\bar{g}^{K}(x)$.
	\begin{corollary}
		\label{cor:maxaffine}
		When $g$ is the maximum of affine functions with domain $\dom = \reals^m$ and polyhedral support $\supp = \{u \mid Hu\leq h\}$, for the same $x$ and $\epsilon$,
		\begin{equation*}
			\bar{g}^{N}(x) -\delta(K,z,\gamma)  \leq \bar{g}^K(x) \leq \bar{g}^{N*}(x) 
		\end{equation*}
	\end{corollary}
This follows from the fact that $L_j = 0$ for all affine functions $g_j$.}
\reviewChanges{
\paragraph{Uncertain objective.}%
\label{par:uncertain_objective_max}
When the uncertainty is in the objective, Theorem~\ref{thm:maxconcave} and Corollary~\ref{cor:maxaffine} quantifies the possible difference in optimal values between the \gls{MRO} problem with $K$ and $N$ clusters. We again define ${\Delta = \max_{x}\left(\bar{g}^{N*}(x) - \bar{g}^{N}(x)\right)}$, subject to $x$ being feasible for problem~\eqref{eq:robustopt}. 
Note, however, that this is only needed for the upper bound. The lower bound holds without needing to consider $\bar{g}^{N*}(x)$; we need not consider the effect of the support set on the problem.
\begin{corollary}
	\label{cor:obj_unc_max}
Consider the problem where $g$ is itself the objective function we would like to minimize and $X \subseteq \reals^n$ represents the constraints, which are deterministic.
Then, $\delta(K,z,\gamma)$ upper bounds the possible decrease in optimal values of the \gls{MRO} problem with $K$ clusters compared with that of $N$ clusters. 
Similarly, $\max_{j\leq J}(L_j/2)D(K) + \Delta$ upper bounds the possible increase in optimal values.
\end{corollary}}
\reviewChanges{
\paragraph{Uncertain constraints.}%
\label{par:uncertain_constraints_max}
When the uncertainty is in the constraints, the difference between $\bar{g}^N(x)$ and $\bar{g}^K(x)$ as given in Theorem~\ref{thm:maxconcave} no longer directly reflect the difference in optimal values.
Instead, clustering affects the feasible set for $x$ as follows.
In the case $\bar{g}^N(x) \geq \bar{g}^K(x)$, for any $\hat{x}$, $\bar{g}^K(\hat{x})$ can be at most $\delta(K,z,\gamma)$ lower in value than $\bar{g}^N(\hat{x})$.
Since both values are constrained to be nonpositive from~\eqref{eq:robustopt}, the feasible region of the \gls{MRO} problem with $K$ clusters may be less restricted than that of $N$ clusters. This indirectly allows \gls{MRO} with $K$ clusters to obtain a smaller optimal value.
On the other hand, in the case $\bar{g}^N(\hat{x}) \leq \bar{g}^K(\hat{x})$, 
$\bar{g}^K(\hat{x})$ can be at most $\max_{j\leq J}(L_j/2)D(K) + \Delta$ higher in value than $\bar{g}^N(\hat{x})$. 
Since both values are constrained to be nonpositive, the feasible region of the \gls{MRO} problem with $K$ clusters may be more restricted than that of $N$ clusters. This indirectly lets \gls{MRO} with $K$ clusters to obtain a larger optimal value.
}

\section{Parameter selection \reviewChanges{and outliers}}
\label{sec:parameter_selection}

\paragraph{Choosing $K$.}
When the uncertain constraint is affine and $\supp$ does not affect the worst-case constraint value, the number of clusters $K$ does not affect the final solution, so it is always best to choose $K = 1$.
We trivially cluster by averaging all data-points without using any clustering algorithm.
When $\supp$ affects the worst-case constraint value, there is a difference of at most $\Delta$ between setting $K=1$ and $K=N$, which can often be approximated $\approx 0$ for small $\epsilon$. Therefore, setting $K = 1$ remains the recommendation.
When the constraint is concave, we choose $K$ to obtain a reasonable upper bound on $\bar{g}^K(x)$, as described in Theorem~\ref{thm:1}.
This upper bound depends linearly on $D(K)$, the clustering value, so by choosing the {\it elbow} of the plot of $D(K)$, we choose a cluster number that, while being a reasonably low value, best conforms to the shape of the underlying distribution.
\reviewChanges{When the constraint is maximum-of-concave, the bounds given in Theorem~\ref{thm:maxconcave} also depends on $D(K)$. Notice that the value $\delta(K,z,\gamma)$ in the lower bound is also a weighted transformation of $D(K)$.}
The elbow method has been commonly used in machine learning problems pertaining the choice of hyper-parameters, especially for $K$-means, and can be traced back to Thorndike~\cite{elbow} in 1953.
Note that, by directly returning $D(K)$ and examining the elbow as an initial step, this procedure can be completed in the clustering step without having to solve the downstream optimization problem.
\reviewChanges{To further improve the choice of $K$, or if the elbow is unclear, cross-validation may be used for low $K$ values or $K$ values around the elbow. }
No matter if the uncertainty lies in the objective or the constraints, this bound will inform us of the potential difference between choosing different $K$.

\paragraph{Choosing $\epsilon$.}
While we have outlined theoretical results in Theorem~\ref{theorem_mro} for choosing $\epsilon$, in practice, we experimentally select $\epsilon$ through cross validation to arrive at the desired guarantee.
Therefore, while the theoretical bounds suggest to choose a larger $\epsilon$ when we cluster, this may not be the case experimentally.
\reviewChanges{In fact, for concave $g$, we may even choose a smaller $\epsilon$, due to the increase in the level of conservatism for small $K$. 
On the other hand, for maximum-of-concave $g$, we do need to choose larger $\epsilon$, as smaller $K$ leads to less conservative solutions. 
However, for both cases, we show a powerful result in the upcoming numerical examples: although for the same $\epsilon$, \gls{MRO} with $K$ clusters differs in conservatism from Wasserstein \gls{DRO} ($N$ clusters), there are cases where we can tune $\epsilon$ such that \gls{MRO} and \gls{DRO} provide almost identical tradeoffs between objective values and probabilistic guarantees, such that no loss in performance results from choosing a smaller cluster number $K$.}
\reviewChanges{
\paragraph{Data with outliers}
When the provided dataset contains outliers, one might imagine that the centoids created by the clustering algorithm will be biased towards the outliers. While this is true, the weights of the outliers will not increase through clustering, thus the effect of outliers on these clustered Wasserstein balls is not worse than their effect on the original Wasserstein balls, which include the Wasserstein ball around the outlier point. In fact, by clustering the outlier point with other points, \gls{MRO} offers protection against the outlier. We demonstrate this in on the numerical experiment in Section~\ref{sub:news}, where we compare three methods: \gls{MRO}, \gls{MRO} with outlier removal, and \gls{MRO} with the outlier considered as its own cluster.
}

\section{Numerical examples}
\label{sec:examples}

We now illustrate the computational performance and robustness of the proposed method on various numerical examples.
All the code to reproduce our experiments is available, \reviewChanges{in Python}, at
\begin{center}
	\url{https://github.com/stellatogrp/mro_experiments}.
\end{center}
We run the experiments on the Princeton Institute for Computational Science and Engineering (PICSciE) facility with 20 parallel 2.4 GHz Skylake cores.
We solve all optimization problems with MOSEK~\citep{mosek} optimizer with default settings.

\reviewChanges{
All numerical examples are solved through direct reformulations if not stated otherwise.
The calculated in-sample objective value and out-of-sample expected values, as well as the out-of-sample probability of constraint violation, are averaged over 50 independent runs of each experiment. For each run, we generate evaluation data of the same size $N$ as the training dataset.}

\reviewChanges{
For numerical examples with an uncertain objective, the probability of constraint violation is measured as the probability the average out-of-sample value is above the in-sample value.
For numerical examples with an uncertain constraint, the probability of constraint violation is measured as the probability the average constraint value is above zero.}

\reviewChanges{In Sections~\ref{sub:capitalbudget},~\ref{sub:quadraticconcave}, and~\ref{sub:logsumexp}, we demonstrate the performance of \gls{MRO} when the uncertain constraint is concave. In Section~\ref{sub:sparse_portfolio_optimization},~\ref{sub:facility}, and~\ref{sub:news}, we demonstrate the performance of \gls{MRO} for maximum-of-affine uncertainty. }

\subsection{Capital budgeting}
\label{sub:capitalbudget}
We consider the capital budgeting problem in ~\cite[Section 4.2]{ben-tal_deriving_2015}, where we select a portfolio of investment projects maximizing the total net present value (NPV) of the portfolio, while the weighted sum of the projects is less than a total budget $\theta$.
The NPV for all projects is $\eta(u) \in \reals^n$, where for each project $j$, $\eta_j(u)$ is the sum of discounted cash flows $F_{jt}$ over the years $t = 0,\dots, T$, \ie, $\eta_j(u) = \sum_{t = 0}^T F_{jt}/(1+u_j)^t$. Here, $u_j$ is the discount rate of project $j$.
We formulate the uncertain function to be minimized as
\begin{equation*}
	 g(u,x) = -\eta(u)^Tx,
\end{equation*}
where $x = (x_1, \dots, x_n) \in \{0,1\}^n$ is the indicator for selecting each project.
The discount rate $u_j$ is subject to uncertainty, as it depends on several factors, such as the interest rate of the country where project $j$ is located and the level of return the decision-maker wants to compensate the risk.
The function $g$ is concave and monotonically increasing in $u$, and we can define a domain $u \geq 0$ so that Assumption~\ref{assp:dom} and Theorem~\ref{thm:1} applies.
The robust problem can be written as
\begin{equation*}
	\begin{array}{ll}
	\underset{x, t}{\text{minimize}} & \tau\\
	\text{subject to} & \bar{g}(u,x) \leq \tau,\quad u \in \uncset(K,\epsilon)\\
	& h^Tx \leq \theta\\
& x\in \{0,1\},
\end{array}
\end{equation*}
where $h$ is the vector of project weights. We refer to~\eqref{eq:robustopt_p} and arrive at the convex reformulation for $p =2$
\begin{equation}
	\label{eq:capital}
	\begin{array}{ll}
	\text{minimize} & \tau\\
	\text{subject to} & \lambda \epsilon^2 + \sum^{K}_{k=1} w_k s_k \leq \tau\\
	& - F_{0}^Tx + \ones^T(\delta_{k}a - z_{k}) - z_k^T\bar{d}_k + \gamma_k^T(b - C\bar{d}_k) \\
	& \hspace{3cm} +1/(4\lambda)\| C^T\gamma_k - z_k \|^2_2 \leq s_k, \quad k = 1,\dots,K\\
	& \left(-(Y_k)_{jt},F_{jt}x_j,(\delta_k)_{jt}\right) \in  \mathcal{K}^{t/(t+1)}, \quad j = 1,\dots,n, \quad t= 1,\dots,T,\\
	& \hspace{7.4cm} \quad k= 1,\dots,K\\
	& Y_k \ones  = z_{k}, \quad \quad k= 1,\dots,K\\
	&  h^Tx \leq \theta\\
	& \lambda \geq 0, \quad \gamma_k \geq 0, \quad Y_k \leq 0,\quad \delta_{k} \leq 0,  \quad x \in \{0,1\},
	\end{array}
\end{equation}
where $a \in \reals^T$ with $a_t = t^{1/(t+1)} + t^{-t/(t+1)}$ for $t = 1,\dots,T$, and $(x,y,z) \in \mathcal{K}^{\alpha}$ is a power cone constraint given as $x^{\alpha} y^{1 - \alpha} \geq |z|$.
The vector $F_0$ indicates the first column of $F$, and matrix $C$ and vector $b$ encode the support of $u$, which we take to be $\{u\in \reals^m \mid 0 \leq u \leq \ones\}$, where $m = n$.
We have variables $x_j \in \reals$, $z_k \in \reals^{n}$, $Y_{k} \in \reals^{n\times T}$, $\delta_{k} \in \reals^{n\times T}$, $\tau \in \reals$, $\gamma_k \in \reals^{2n}$, $s_k \in \reals$, for $j = 1,\dots,n$, $k= 1,\dots,K$, and $t= 1,\dots,T$.
The derivation of reformulation~\eqref{eq:capital} is in Appendix~\ref{app:capital}.
Note that there are variables with total dimension $KnT$, which grows swiftly when any of the parameters are large.
For each cluster $k$, we introduce $nT$ new variables for $y$ and $\delta$, as well as $nT$ new power cone constraints, which greatly increases the compuational complexity of the problem.

\paragraph{Problem setup.}
We set $n = 20$, $N = 120$, $T= 5$. We generate $F_{jt}$ from a uniform distribution on ${[0.1,0.5+0.004t]}$ for $j = 1,\dots,n$, \quad $t = 0,\dots,T$.
For all $j$, $h_j$ is generated from a uniform distribution on $[1,3-0.5j]$, and the total budget $\theta$ is set to be 12.
We generate uncertain data from two slightly different uniform distributions, to simulate two different sets of predictions on the discount rates. The first half is generated on $[0.005j,0.02j]$, and the other half on $[0.01j,0.025j]$, for all $j$. We calculate an upper bound on the $L$-smooth parameter, ${L = \|\nabla^2 \sum_{j=1}^n \sum_{t=0}^T F_{jt}(\hat{x}_N)_j(1+u_j)^{-t} \|_{2,2} \leq  \|\sum_{j=1}^n \sum_{t=0}^T t(t+1)F_{jt}(\hat{x}_N)_j\|_{2,2}}$ for each data-driven solution $\hat{x}_N$.
\reviewChanges{
\paragraph{Choosing $K$}
Plotting the clustering value $D(K)$ over $K$, we note that the elbow occurs around $K=2$, which suggests using cross-validation for $K$ values around 2. }
\begin{figure}[h]%
	\centering
	\includegraphics[width=0.8\textwidth]{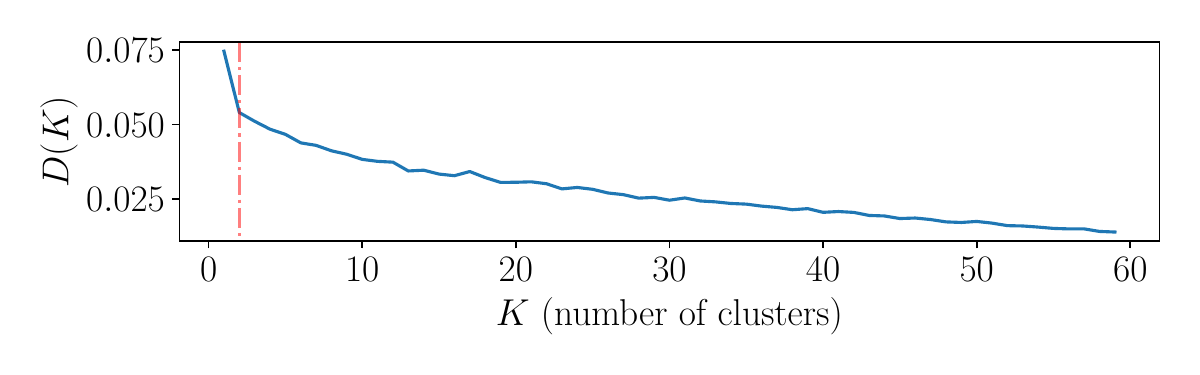}
	\caption{Capital budgeting. $D(K)$ vs. $K$. Dotted red line: $K = 2$.}%
	\label{fig:capitalk}%
\end{figure}
\begin{figure}[h]%
	\centering
	\includegraphics[width=\textwidth]{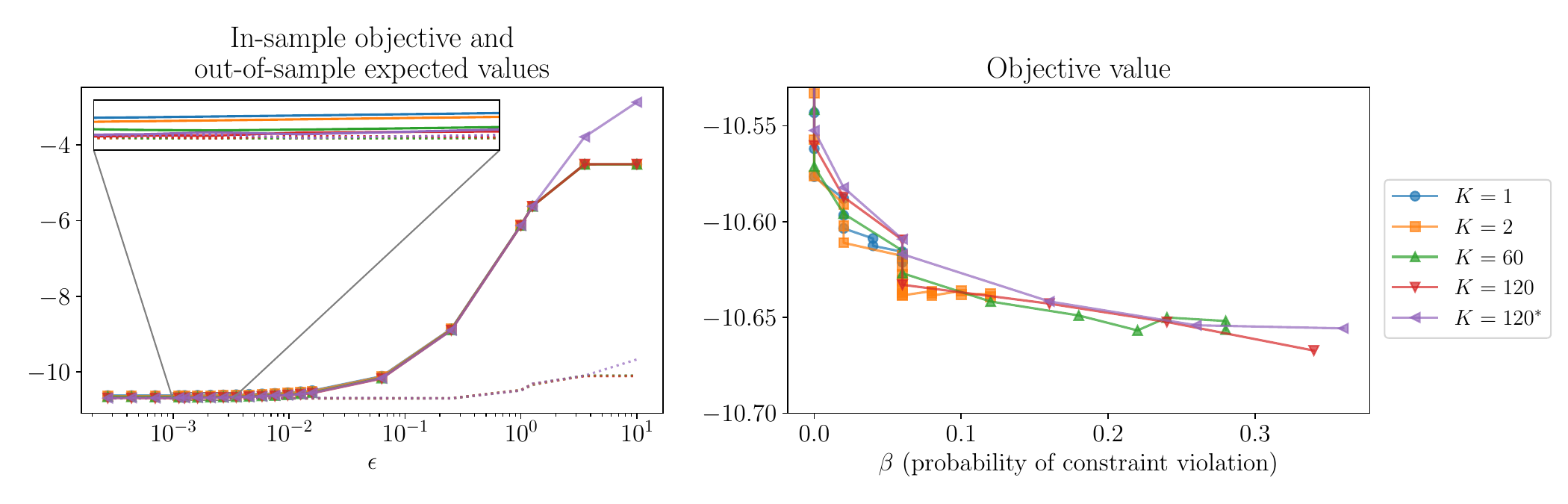}
	\caption{Capital budgeting. Left: in-sample objective values and out-of-sample expected values vs. $\epsilon$ for different $K$. Solid lines are the in-sample objective value, dotted lines are the out-of-sample expected value. Right: objective value vs. $\beta$ for different $K$; each point represents the solution for the $\epsilon$ achieving the smallest objective value. $K = 120^*$ is the formulation without the support constraint. }%
	\label{fig:capital1}%
\end{figure}
\begin{figure}[h]%
	\centering
	\includegraphics[width=\textwidth]{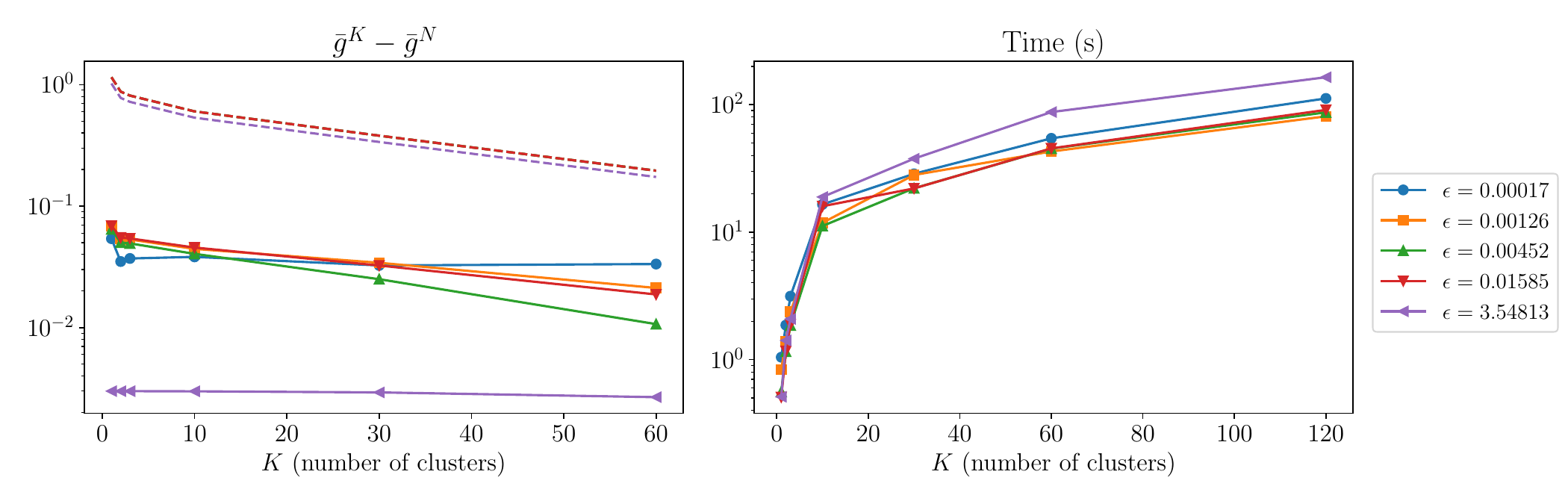}
	\caption{Capital budgeting. Left: the difference in the value of the uncertain objective between using $K$ and $N$ clusters, calculated as $\bar{g}^K(x) - \bar{g}^N(x)$, compared with the theoretical upper bound $(L/2)D(K)$ from Corollary~\ref{cor:obj_unc}. Solid lines are the difference, dotted lines are the upper bounds. Right: solve time. }%
	\label{fig:capital2}%
\end{figure}

\paragraph{Results.}  We observe in Figure~\ref{fig:capital1} that using two clusters is enough to achieve performance almost identical to that of using 120 clusters.
Although from the left image, we see that $K=2$ slightly upper bounds $K = 120$, from the right, their tradeoffs between the objective value and relevant constraint violation probability ($\beta \leq 0.2$) are largely the same, so we can always tune $\epsilon$ to achieve the same performance and guarantees.
Notice that the results for $K = 120$ and $K=120^*$ are near identical for small $\epsilon$, where $K=120^*$ is the formulation without the support constraint. Therefore, while $ \bar{g}^{N*}(x)$ slightly upper bounds $\bar{g}^N(x)$, we can approximate their difference $\Delta \approx 0$ for small enough $\epsilon$, for which the upper bound $(L/2)D(K)$ thus hold.
In fact in this example, even for larger $\epsilon$ where we observe $\Delta > 0$, the actual difference between $\bar{g}^K$ and $\bar{g}^N$ is bounded by $(L/2)D(K)$.
\reviewChanges{
In Figure~\ref{fig:capital2}, we see that the elbow of the upper bound is at $K = 2$, and the true difference follows the same trend, matching the suggestion from Figure~\ref{fig:capitalk}.}
Therefore, setting $K = 2$ is the optimal decision, with a time reduction of 2 orders of magnitude, and a complexity reduction from 26626 variables and 12000 power cones to 666 variables and 200 power cones.

\subsection{Quadratic concave uncertainty}
\label{sub:quadraticconcave}
We refer to the example from Ben-Tal et al.~\cite[Section 4.2]{ben-tal_deriving_2015} with concave uncertainty of the form $$g(u,x) = \sum_{i=1}^n h_i(u) x_i,$$ where $h_i(u) = -(1/2)u^TA_iu$, each $A_i\in \reals^{m\times m}$ a symmetric positive definite matrix, $u \in \reals^m$, and $x \in \reals^{n}_{+}$. For simplicity, we also require that $x$ sums to 1, $p = 2$, and the support of the uncertainty $\supp = \reals^m$. Assuming the uncertainty is in the objective, such that the uncertain constraint is created using epigraph form, we solve the problem
\begin{equation*}
	\begin{array}{ll}
	\text{minimize} & \tau\\
	\text{subject to} & \lambda \epsilon^2 + \sum^{K}_{k=1} w_k s_k \leq \tau\\
	&(1/2)\sum_{i=1}^n \left({(Y_k)_i}^TA_i^{-1}(Y_{k})_i\right)/x_i -z_k^T\bar{d}_k  + 1/(4\lambda) \left\| z_k \right\|^2_2 \leq s_k \\
	&\hspace{7cm} \quad k = 1,\dots,K\\
& Y_{k}\ones = z_k, \quad k = 1,\dots,K\\
	& \ones^T x = 1, \quad \lambda \geq 0,\quad x \geq 0.
	\end{array}
\end{equation*}
The $A_i^{-1}$ terms come from taking the conjugate of $g$, and the derivation can be found in~\cite[Example 24]{ben-tal_deriving_2015}. We have variables $x \in \reals^n$, $z_k \in \reals^{m}$, $Y_{k} \in \reals^{m\times n}$, $\tau \in \reals$, $s_k \in \reals$, for $k= 1,\dots,K$. We let $(Y_k)_i$ indicate the $i$th column of $Y_k$.

\paragraph{Problem setup.}
We set $n = m = 10, N = 90$, and generate synthetic uncertainty as a multi-modal normal distribution with 5 modes, where $\mu_i = \gamma_j0.03i$ for all $i = 1,\dots, n$ for mode $j$, with mode scales $\gamma = (1, 5,15,25,40)$. The variance is $\sigma_i = 0.02^2 + (0.025i)^2$ for all modes. We generate $A_i$ as random positive semi-definite matrices for all $i = 1,\dots, n$. For the upper bound, we calulate $L = \|\sum_{i=1}^n A_i(\hat{x}_{N})_i\|_{2,2}$ for each data-driven solution $\hat{x}_N$.
\reviewChanges{
\paragraph{Choosing $K$}
Plotting the clustering value $D(K)$ over $K$, we note that the elbow occurs at $K=5$, which suggests a choice of $K=5$. }
\begin{figure}[h]%
	\centering
	\includegraphics[width=0.8\textwidth]{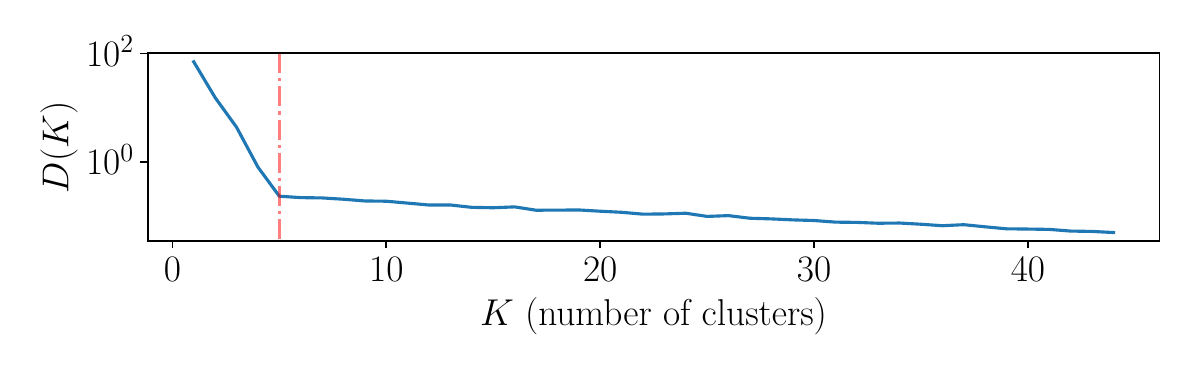}
	\caption{Quadratic concave uncertainty. $D(K)$ vs. $K$. Dotted red line: $K = 5$.}%
	\label{fig:concavek}%
\end{figure}
\begin{figure}[h]%
	\centering
	\includegraphics[width=\textwidth]{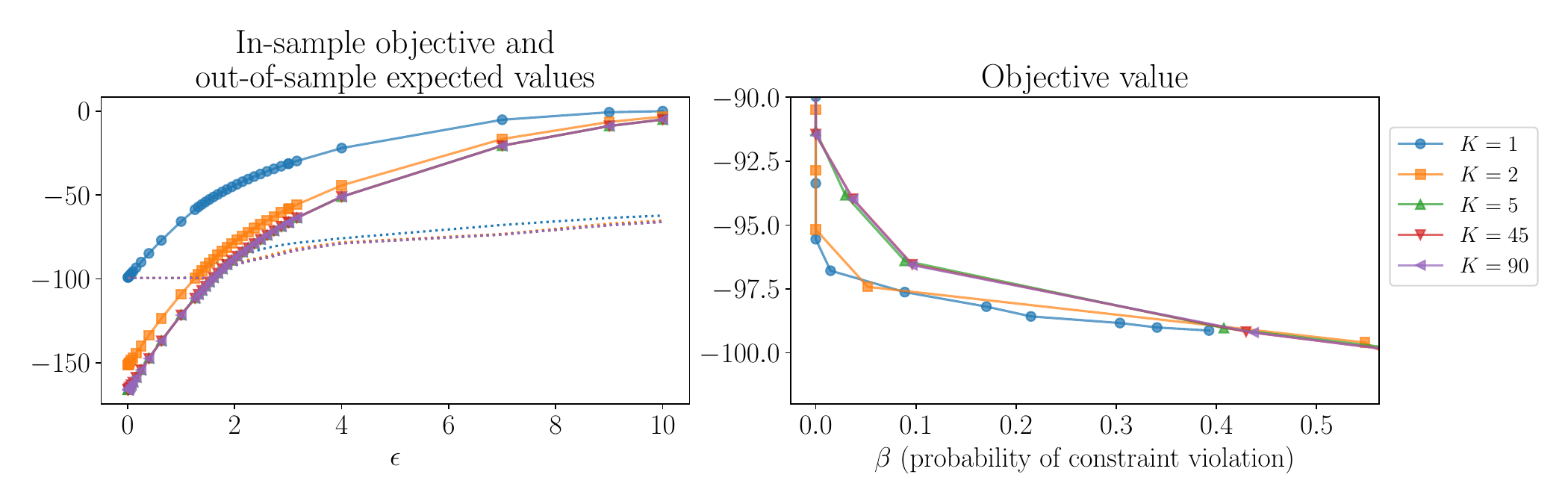}
	\caption{Quadratic concave uncertainty. Left: in-sample objective values and out-of-sample expected values of $g$ vs. $\epsilon$ for different $K$. Solid lines are the in-sample objective value, dotted lines are the out-of-sample expected value. Right: objective value vs. $\beta$ for different $K$; each point represents the solution for the $\epsilon$ achieving the smallest objective value. }%
	\label{fig:concave_syn1}%
\end{figure}
\begin{figure}[h]%
	\centering
	\includegraphics[width=\textwidth]{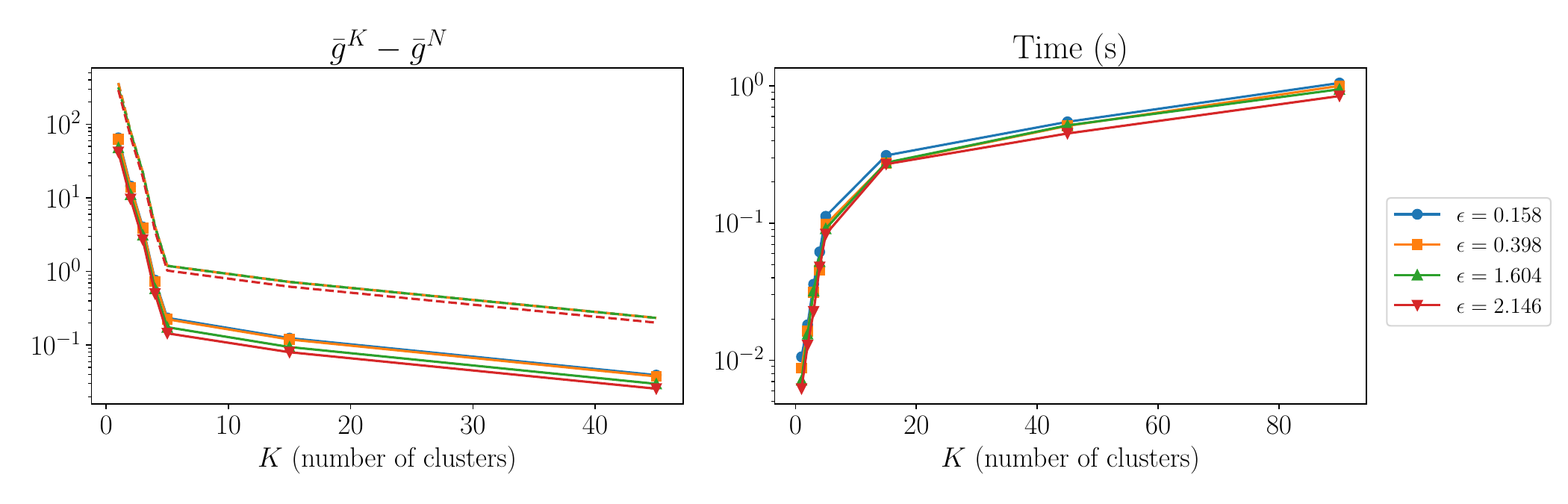}
	\caption{Quadratic concave uncertainty. Left: the difference in the value of the uncertain objective between using $K$ and $N$ clusters, calculated as $\bar{g}^K(x) - \bar{g}^N(x)$, compared with the theoretical upper bound $(L/2)D(K)$ from Corollary~\eqref{cor:obj_unc}. Solid lines are the difference, dotted lines are the upper bounds. Right: solve time. }%
	\label{fig:concave_syn2}%
\end{figure}

\paragraph{Results.}  We observe on the left of Figure~\ref{fig:concave_syn1} that using 5 clusters is enough to achieve performance almost identical to that of using 90 clusters.
\reviewChanges{
Indeed, in Figure~\ref{fig:concave_syn2}, the elbow of the upper bound (dotted lines) on the difference in objective values is at $K = 5$, and the true difference follows the same trend, corroborating with Figure~\ref{fig:concavek}.}
Furthermore, on the left plot of Figure~\ref{fig:concave_syn2}, we note for $K \geq 5$, the tradeoff between the objective value and constraint violation is the same, so we can tune $\epsilon$ to achieve the same performance and guarantees.
In fact, for this particular example, using a smaller $K$ such as 1 or 2 may allow us to tune $\epsilon$ to achieve an even better tradeoff. However, this result cannot be guaranteed in general, so the recommended action is still to choose $K=5$.

\subsection{Robust log-sum-exp optimization}
\label{sub:logsumexp}
We also consider uncertainty from Bertsimas and den Hertog~\cite[Chapter 14]{robustadaptopt} of the form
\begin{equation*}
	g(u,x) = \log\left(\sum_{i = 1}^n u_i e^{x_i}\right),
\end{equation*}
concave in $u$ and convex in $x$. This function $g$ is monotonically increasing in $u$, and we can define a domain $u \geq 0.01$ so that Assumption~\ref{assp:dom} and Theorem~\ref{thm:1} apply.
Assuming the simple case where the uncertainty is in the objective, we add some further restrictions on $x$ and use a cutting plane procedure to solve, for $p = 2$,
\begin{equation*}
	\begin{aligned}
	\underset{x}{\text{minimize}} \hspace{0.5em} \underset{u_1\dots u_k}{\text{maximize}} &\quad \sum_{k=1} ^K w_k  \log\left(\sum_{i = 1}^n u_i e^{x_i}\right) \\
	\mbox{subject to} & \quad \sum_{k=1} ^K w_k \| u_k - \bar{d}_k \|^2 \le \epsilon^2\\
	&\quad   \ones^T x \geq 10, \quad x \geq 0, \quad x \leq 10
	\end{aligned}
\end{equation*}
\paragraph{Problem setup.} We set $n = 30, N = 90$, and observe synthetic data from 3 sets of uniform distributions, scaled respectively by $\gamma = (1,3,7)$. Specifically, for each set $j$, each $d_i$ is generated uniformly on the intervals $0.01[\gamma_ji,\gamma_j(i+1)]$ for $i = 1,\dots, n$. For the upper bound, we calculate $L = \|\nabla^2 g\|_{2,2} \leq \exp(\hat{x}_N)^T\exp(\hat{x}_N)\min({d_k}^T\exp(\hat{x}_N))^{-2}$, for each data-driven solution $\hat{x}_N$.
\reviewChanges{
\paragraph{Choosing $K$}
Plotting the clustering value $D(K)$ over $K$, we note that the elbow occurs at $K=3$, which suggests using cross-validation for $K$ values around 3. }
\begin{figure}[h]%
	\centering
	\includegraphics[width=0.8\textwidth]{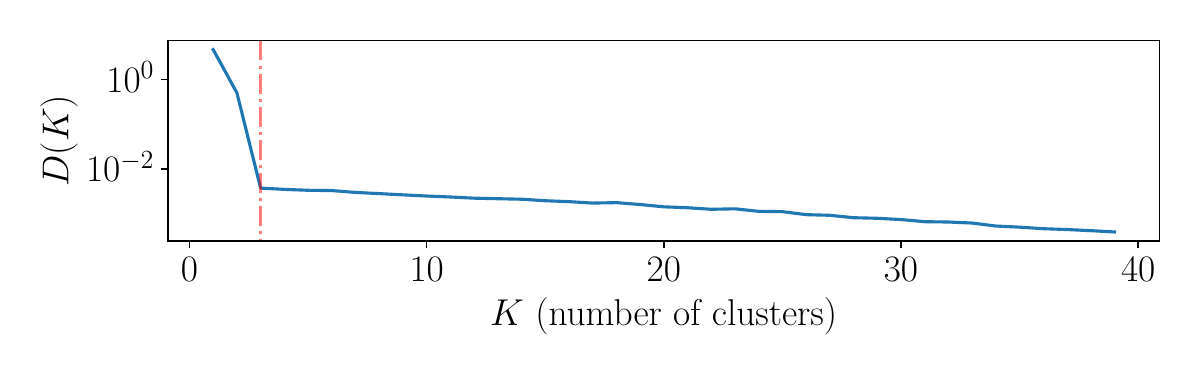}
	\caption{Log-sum-exp uncertainty. $D(K)$ vs. $K$. Dotted red line: $K = 3$.}%
	\label{fig:logk}%
\end{figure}
\begin{figure}[h]%
	\centering
	\includegraphics[width=\textwidth]{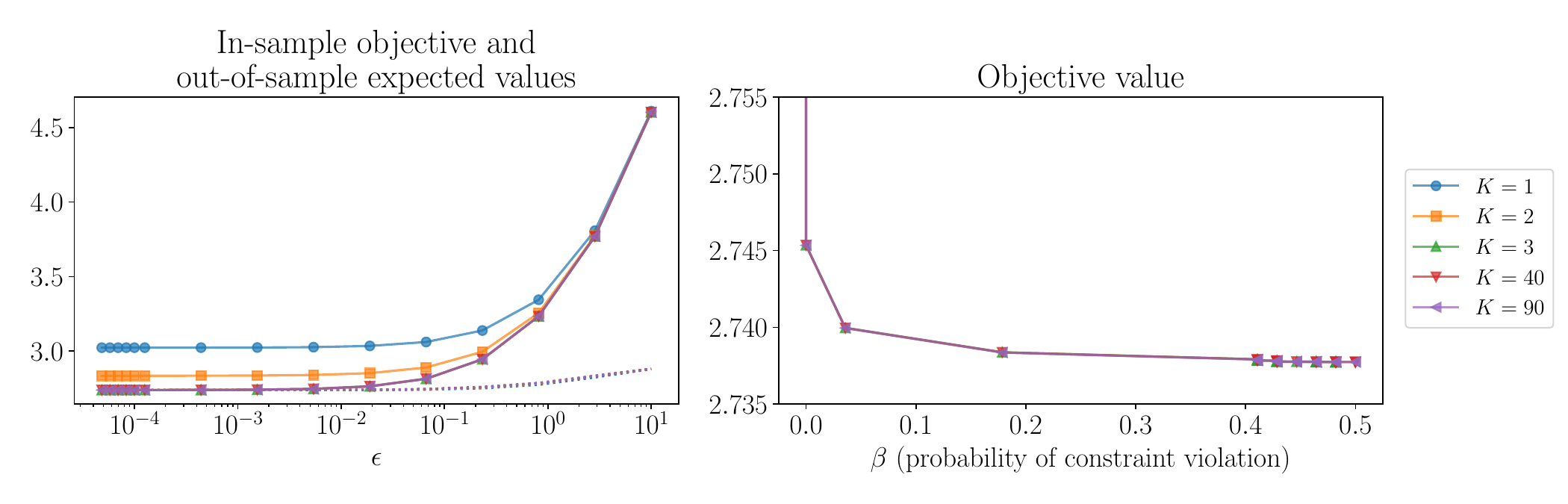}
	\caption{Log-sum-exp uncertainty. Left: in-sample objective values and out-of-sample expected values vs. $\epsilon$ for different $K$. Solid lines are the in-sample objective value, dotted lines are the out-of-sample expected value. Right: objective value vs. $\beta$ for different $K$. }%
	\label{fig:log_syn1}%
\end{figure}

\begin{figure}[h]%
	\centering
	\includegraphics[width=\textwidth]{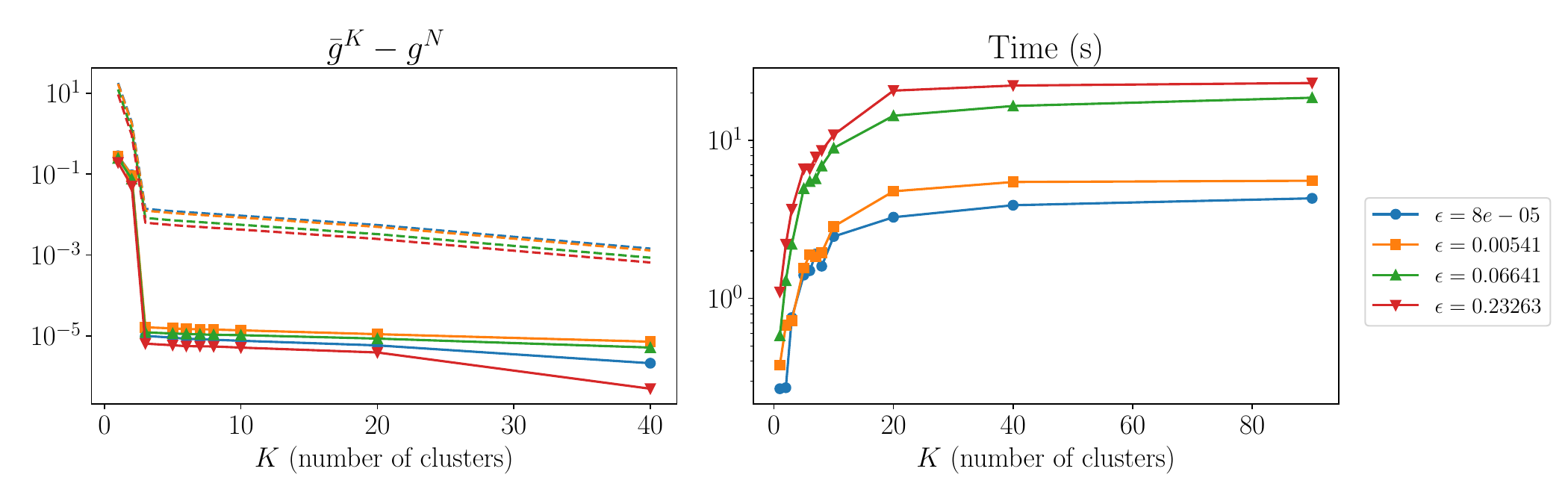}
	\caption{Log-sum-exp uncertainty. Left: the difference in the value of the uncertain objective between using $K$ and $N$ clusters. Solid lines are the difference, the dotted line is the upper bound. Right: solve time. }%
	\label{fig:log_syn2}%
\end{figure}

\paragraph{Results.} We observe on the left of Figure~\ref{fig:log_syn1} that while setting $K$ to smaller values increases the objective value, setting $K = 3$, the number of modes of the underlying distribution, already offers near identical performance to that of setting $K = 90$.
\reviewChanges{On the left of Figure~\ref{fig:log_syn2}, we see that $K = 3$ is at the elbow of upper bound and actual difference, corroborating with Figure~\ref{fig:logk}.}
Furthermore, we note that setting $K=3$ and above give identical tradeoff curves, therefore, choosing $K=3$ is the time-efficient solution.

\subsection{\reviewChanges{Sparse portfolio optimization}}%
\label{sub:sparse_portfolio_optimization}
\reviewChanges{
We consider a market that forbids short-selling and has $m$ assets as in~\citep{mohajerin_esfahani_data-driven_2018}.
Daily returns of these assets are given by the random vector $d = (d_1, \dots, d_m) \in \reals^m$.
The percentage weights (of the total capital) invested in each asset are given by the decision vector ${x = (x_1, \dots, x_n) \in \reals^n}$.
We restrict our selection to at most $\theta$ assets.
The underlying data-generating distribution $\prob$ is unknown, but we have observed a historical dataset  $\mathcal{D}_N$.
Our objective is to minimize the CVaR with respect to variable $x$,
\begin{equation*}
	\label{eq:portfolio}
	\begin{array}{ll}
	\text{minimize  } &\CVaR(- u^Tx,\alpha)\\
	\text{subject to} & \ones^T x = 1\\
			  & x \ge 0, \quad \|x\|_1 \leq \theta,
	\end{array}
\end{equation*}
which represents the average of the $\alpha$ largest portfolio losses that occur. In other words, the $\CVaR$ term seeks to ensure that the expected magnitude of portfolio losses, when they occur, is low.
The objective has an analytical form with an extra variable $\tau$ given as~\citep{cvaropt,mohajerin_esfahani_data-driven_2018}:
\begin{equation*}
	\text{minimize} \quad \Expect^{\prob}\left(\tau +\frac{1}{\alpha} \max\{-u^Tx -\tau, 0\}\right).
\end{equation*}
From this, we obtain $g$ as the maximum of affine functions,
\begin{equation*}
	g(u,x) \define \max\{(-1/\alpha)x^Tu + (1 - 1/\alpha)\tau, \tau \}.
\end{equation*}
Using the formulation~\eqref{eq:unif-dual-p} with $p= \infty$, we can write a convex reformulation of the form
\begin{equation*}
	\label{eq:portfolio_1}
	\begin{array}{ll}
	\text{minimize } & y\\
	\text{subject to } & \epsilon \left\|(1/\alpha)x\right\|_2 + \sum^{K}_{k=1} w_k s_k \leq y\\
	&  (1 - 1/\alpha)\tau - (1/\alpha)x^T\bar{d}_k    \leq s_k,\quad k = 1,\dots,K\\
	&  s_k \geq \tau, \quad k = 1,\dots,K\\
	& x \geq 0, \quad \mathbf{1}^T x = 1\\
	& z \in \{0,1\}, \quad z - x \leq 0, \quad \ones^Tz \leq \theta.
	\end{array}
\end{equation*}
with variables $x \in \reals^m$, $z \in \reals^m$, $y \in \reals$, $\tau \in \reals$, $s_k \in \reals$. The variables $z$ are introduced to replace the cardinality constraint using big-$M$ formulation~\citep{sparseportfolio}.}
\reviewChanges{
\paragraph{Problem setup.} We take stock data from the past 10 years of S\&P500, and generate synthetic data from their fitted general Pareto distributions.
We choose a generalized Pareto fit over a normal distribution as it better models the heavy tails of the returns~\citep{tails}.
See the Github repository for the code, which uses the "Rsafd" R package~\cite{Rsafd}.
We let $\alpha = 20\%$, $m = 50$ stocks, and generate a dataset size of $N = 1000$. Our portfolio can include at most $\theta = 5$ stocks. For the upper bound $\delta(K,z,\gamma)$ on $\bar{g}^N(x)-\bar{g}^K(x)$, we note the special structure of this problem, where one of the affine pieces is independent of $u$, to arrive at a bound $\delta(K,z,\gamma) = \max_k\{\max_i\{(\bar{d}_k - d_i)^Tx/\alpha\}\}$.}

\reviewChanges{
\paragraph{Choosing $K$}
Plotting the clustering value $D(K)$ over $K$, we note that the elbow occurs at $K=5$, which suggests using cross-validation for $K$ values around 5.}
\begin{figure}[h]%
	\centering
	\includegraphics[width=0.8\textwidth]{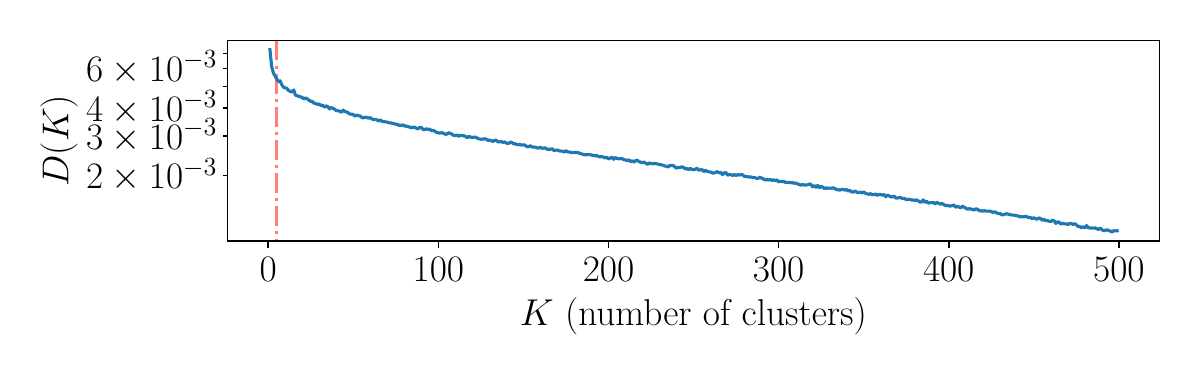}
	\caption{Sparse portfolio. $D(K)$ vs. $K$. Dotted red line: $K = 5$.}%
	\label{fig:portk}%
\end{figure}

\begin{figure}[h]%
	\centering
	\includegraphics[width=\textwidth]{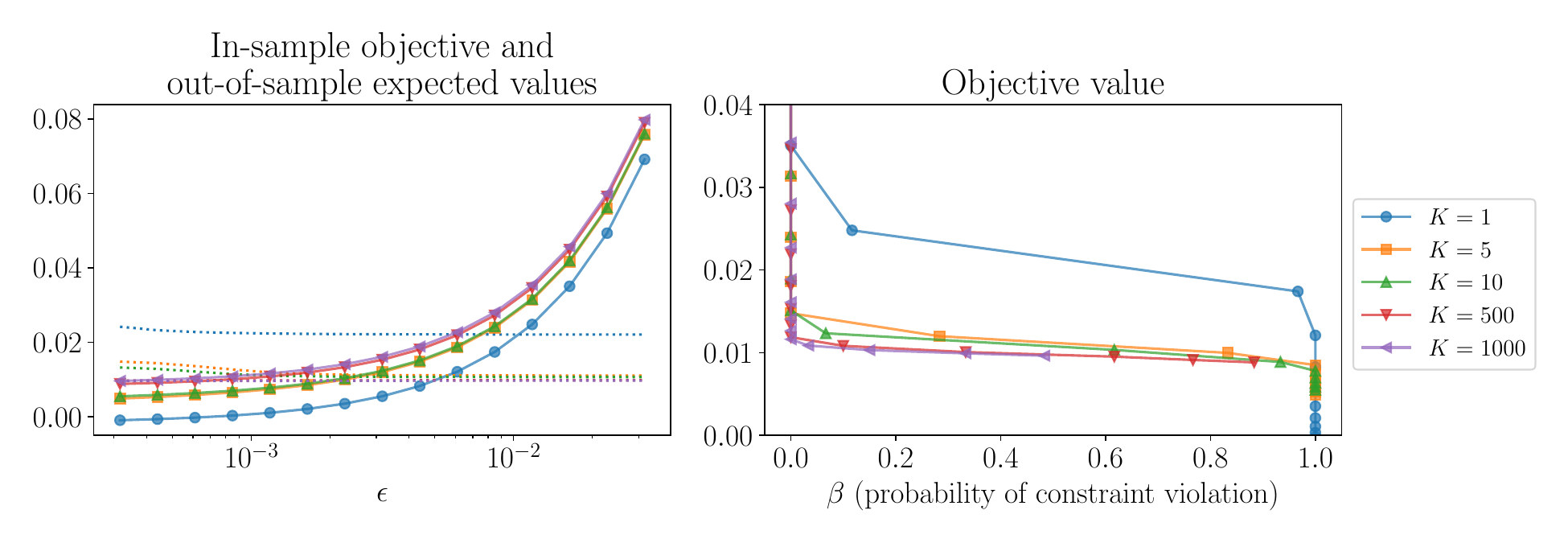}
	\caption{Sparse portfolio. Left: in-sample objective values and out-of-sample expected values vs $\epsilon$ for different $K$. Solid lines are the in-sample objective value, dotted lines are the out-of-sample expected value. Right: objective value vs $\beta$ for different $K$; each point represents the solution for the $\epsilon$ achieving the smallest objective value.}%
	\label{fig:portfolio1}%
\end{figure}

\begin{figure}[h]%
	\centering
	\includegraphics[width=\textwidth]{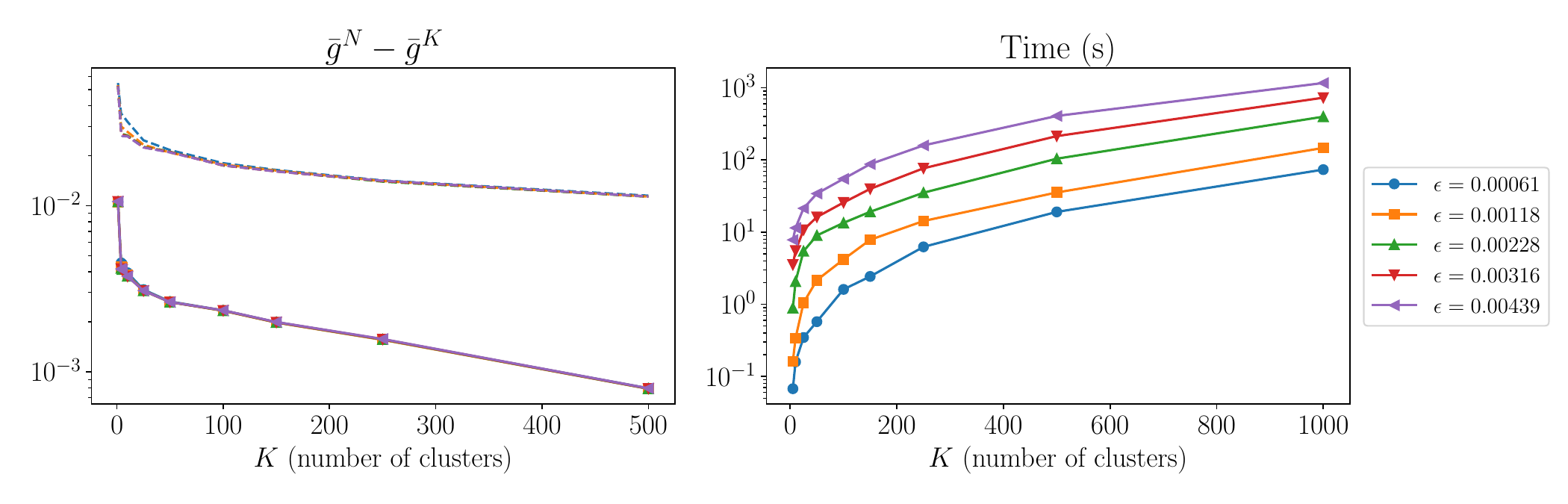}
	\caption{Sparse portfolio. Left: the difference in the value of the uncertain objective between using $N$ and $K$ clusters, calculated as $\bar{g}^N(x) - \bar{g}^K(x)$, compared with the theoretical upper bound $\delta(K,z,\gamma)$ from Corollary~\ref{cor:maxaffine}. Solid lines are the difference, dotted lines are the upper bounds. Right: solve time for $K \geq 5$. }%
	\label{fig:portfolio2}%
\end{figure}

\reviewChanges{
\paragraph{Results.}%
In Figure~\ref{fig:portfolio1}, while setting $K$ to smaller values lead to a decrease in the optimal value across $\epsilon$, we note that for $K=5$ and above, we can already achieve a tradeoff curve between the optimal value and probability of constraint satisfaction that is similar to that of $K = 1000$, and setting $K=10$ brings it slightly closer. In Figures~\ref{fig:portk} and~\ref{fig:portfolio2}, in the plots of $D(K)$ and of the upper bound on the difference, we also note that the elbow is around $K=5$. We thus recommend choosing $K$ through cross validation around 5, as tuning $\epsilon$ for these small $K$ gives 1 - 3 orders of magnitude time reduction.}


\subsection{Facility location}
\label{sub:facility}
We examine the classic facility location problem~\citep{bertsimas_probabilistic,HOLMBERG1999544}.
Consider a set of $n$ potential facilities, and $m$ customers.
Variable $x \in \{0, 1\}^n$ describes whether or not we construct each facility $i$ for $i=1, \dots, n$,  with cost $c_i$.
In addition, we would like to satisfy the uncertain demand $u\in \reals^m$ at minimal cost.
We define variable $X\in\mathbf{R}^{n\times m}$ where $X_{ij}$ corresponding to the portion of the demand of customer $j$ shipped from facility $i$ with corresponding cost $C_{ij}$.
Furthermore, $r \in \reals^n$ represents the production capacity for each facility, and $u\in \reals^m$ represents the uncertain demand from each customer.
For each customer $j$, $X_j$ represents the proportion of goods shipped from any facility to that customer, which sums to $1$.
For each facility $i$, $(X^T)_i$ represents the proportion of goods shipped to any customer.
Putting this all together, we obtain multiple affine uncertain capacity constraints,
\begin{equation*}
	g_i(u,x) \define (X^T)_i u - r_i x_i \leq 0 \quad i = 1,\dots,n,
\end{equation*}
\reviewChanges{which we combine to create a single maximum-of-affine constraint,
\begin{equation*}
	g(u,x) \define \max_{i \leq n}((X^T)_i u - r_i x_i) \leq 0.
\end{equation*}
Now, to ensure a high probability of constraint satisfaction, we use the $\CVaR$ reformulation, 
\begin{equation*}
	g(u,x,\tau) \define \tau + (1/\alpha)\max \left (\max_{i \leq n}((X^T)_i u - r_i x_i - \tau), 0\right ) \leq 0,
\end{equation*}
where we add the auxiliary variable $\tau$. We assume a polyhedral support $\supp= \{u \mid Hu\leq b\}$ for the demand, and solve the problem, for $p = \infty$,
\begin{equation*}
	\label{eq:mro_facility_mip_ref}
	\begin{array}{ll}
	\text{minimize} & c^Tx  +\Tr(C^TX)\\
	\text{subject to}  &\ones^TX_j = 1,\quad j = 1,\dots, m \\
	&   \tau + \sum_{k=1}^K w_k s_{k}  \leq 0,\\
	 & -(1/\alpha)\tau + \lambda_{k} \epsilon + (1/\alpha)((X^T)_i\bar{d}_k -r_ix_i) + \gamma_{ik}(b - H\bar{d}_k) \leq s_{k},\\
	 &\hspace{4cm} \quad i = 1,\dots,n, \quad k = 1,\dots, K\\
	 & \lambda_{k} \epsilon  \leq s_{k}, \quad k = 1,\dots, K\\
	 & \left\| H^T\gamma_{ik} + (1/\alpha)(X^T)_i\right\|_2 \leq \lambda_{k}, \quad i = 1,\dots,n, \quad k = 1,\dots, K\\
	 & \gamma_{ik} \geq 0, \quad i = 1,\dots,n, \quad k = 1,\dots, K\\
	 & x\in \{0,1\}^n, \quad X\in \reals^{n\times m}.
	\end{array}
\end{equation*}
We have variables $x \in \{0,1\}^{n}$, $X \in \reals^{n\times m}$, $s_{k} \in \reals$, $\tau \in \reals$, $\lambda_k \in \reals$, $\gamma_{ik} \in \reals^{m}$, for $i = 1,\dots,n$ and $k= 1,\dots,K$.
The $\gamma$ variables arise from enforcing the support constraints.}


\paragraph{Problem setup.}%
To generate data, we set $n = 5$ facilities, $m = 25$ customers, and $N = 50$ data samples.  \reviewChanges{For the $\CVaR$ reformulation, we set $\alpha = 20\%$. We set costs $c = (46.68, 58.81, 30, 42.09, 35.87)$, and generate the two coordinates of each customer's location from a uniform distribution on $[0,15]$. We then calculate $C$ as the $\ell_2$ distance between each pair of customers. We set production capacities $r = (33, 26, 41, 26, 22)$.
We assume the demand $d$ is supported between 1 and 6, which we write as $Hu\leq b$, where $H = [-I\; I]^T$ and $b$ is the concatenation of two vectors: a vector of $-1$'s of length $m$, and a vector of $6$'s of length $m$. 
We generate demands as the combination of two normal distributions. Half of the data is generated from the normal distribution with mean $\mu_1 = 3$ and variance $\sigma_1 = 0.9$, the second half has mean $\mu_1 = 4$ and variance $\sigma_1 = 0.8$.
We then project the demands onto $(1,6)$.}
\reviewChanges{For the upper bound $\delta(K,z,\gamma)$ on  $\bar{g}^N(x) - \bar{g}^K(x)$ from Corollary~\ref{cor:maxaffine}, we have $(1/N)\sum_{k=1}^K \sum_{i \in C_k}  \max(\max_{j\leq J}((X[i] -H^T\gamma_{jk}),0)^T(d_i - \bar{d}_k))$. Note that this upper bounds the difference in constraint values, and only indirectly affects the objective values through restrictions on the feasible region. Therefore, it is not an upper bound on the difference in objective value, merely a rough estimate.
We cannot directly compare this upper bound against the change in constraint values, as at the optimal chonsen $x, X$ for each $K$, which differ due to differences in the feasible regions, the constraint value will always be near 0 for optimality. We thus compare it against the change in objective values. }
\reviewChanges{
\paragraph{Choosing $K$}
Plotting the clustering value $D(K)$ over $K$, we note that the elbow occurs at $K=2$, which suggests using cross-validation for $K$ values around 2. }
\begin{figure}[h]%
	\centering
	\includegraphics[width=0.8\textwidth]{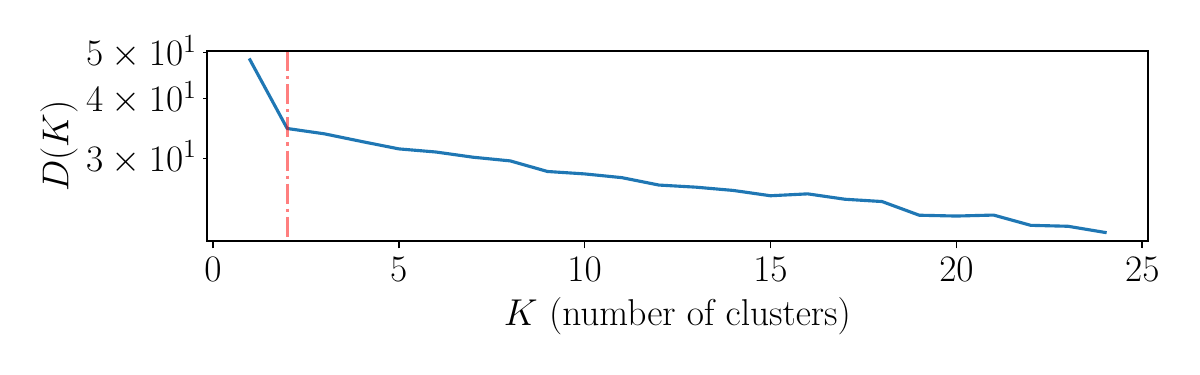}
	\caption{Facility location. $D(K)$ vs. $K$. Dotted red line: $K = 2$.}%
	\label{fig:facilityk}%
\end{figure}

\begin{figure}[h]%
	\centering
	\includegraphics[width=\textwidth]{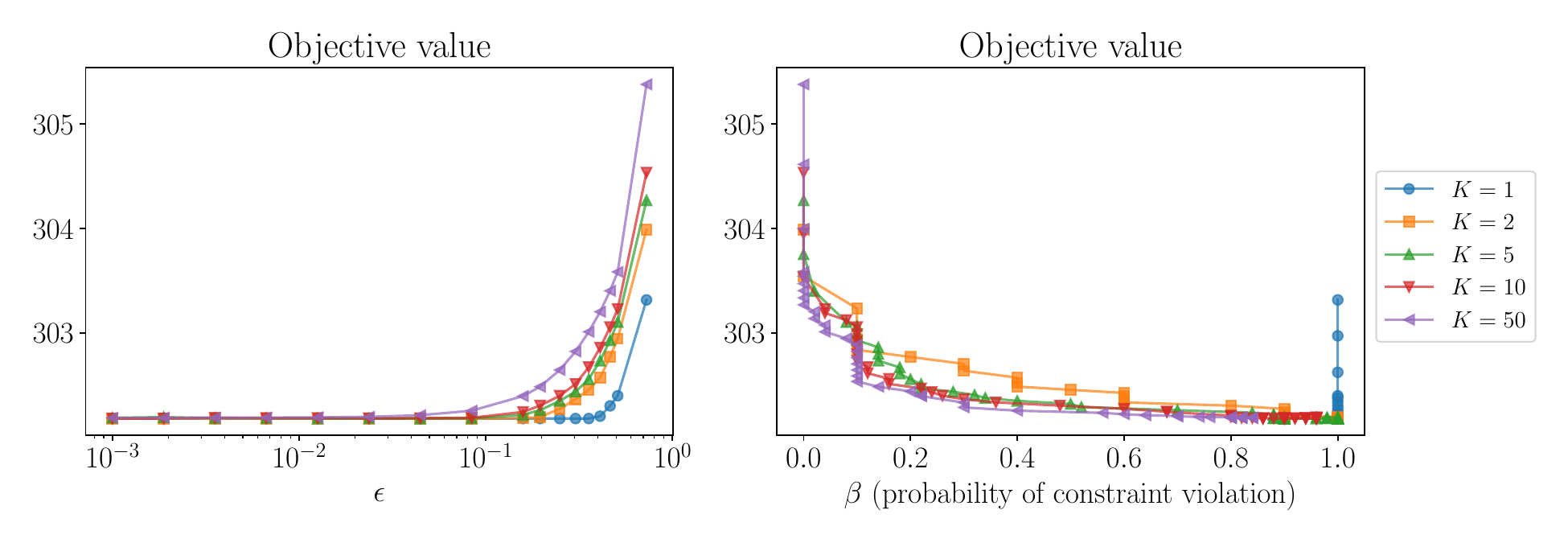}
	\caption{Facility location. Left: in-sample objective values vs $\epsilon$ for different $K$. Right: objective value vs $\beta$ for different $K$; each point represents the solution for the $\epsilon$ achieving the smallest objective value. }%
	\label{fig:facility}%
\end{figure}

\paragraph{Results.}%
\reviewChanges{
As expected of maximum-of-affine $g$, we note in Figure~\ref{fig:facility} that setting $K$ to smaller values lead to a decrease in the optimal value across different $\epsilon$ values. 
While $K=1$ yields poor performance in terms of the probability of constraint violation, we observe that $K=2$ already yields a tradeoff between the objective and probability of constraint violation close to that of $K=50$. Through cross-validation with different $K$, we select $K=5$, which provides a tradeoff curve closer to optimality. 
As this problem has uncertainty in the constraints and not the objective, the bounds given in Corollary~\ref{cor:maxaffine} do not directly reflect the difference in the objective values. 
However, they do give a reference value and inform us of the general trend of the difference.
In this case, they still upper bound the actual difference, as shown in Figure~\ref{fig:facility}.
We note that the bounds we use do not depend on $\bar{g}^{N*}(x)$, so it is irrelevant whether or not the support has an affect on the worst-case constraint value.
Overall, choosing $K = 5$ leads to a time reduction of an order of magnitude while achieving near-optimal performance. 
}
\begin{figure}[h]%
	\centering
	\includegraphics[width=\textwidth]{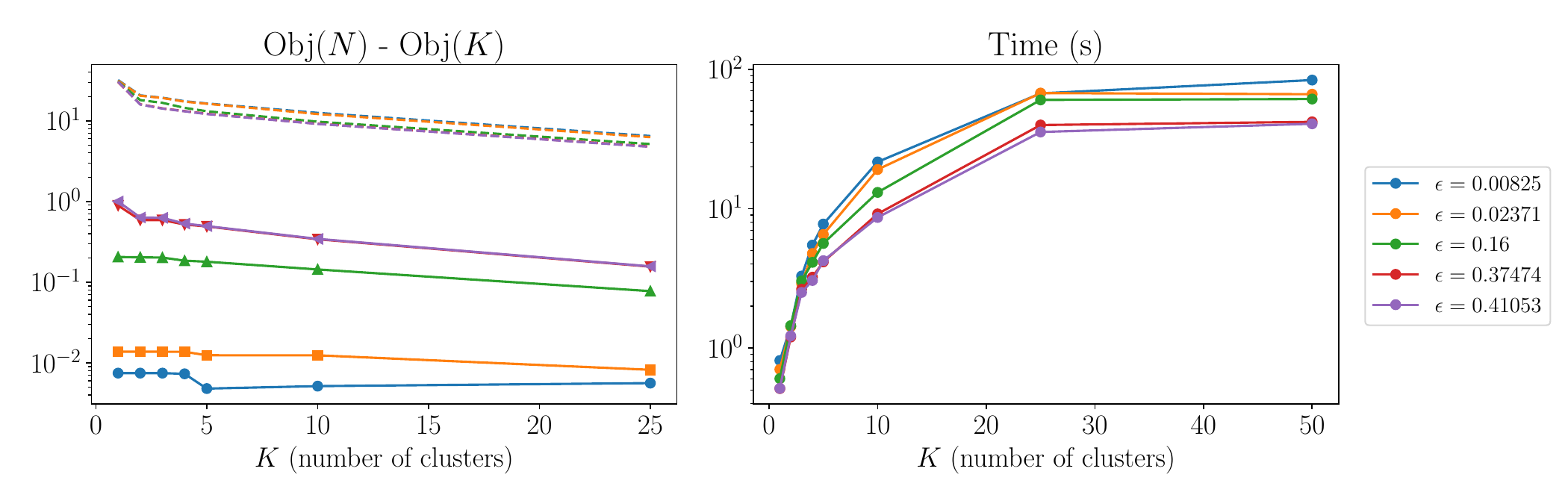}
	\caption{Facility location. Left: the difference in the value of the uncertain objective between using $N$ and $K$ clusters, calculated as Obj($N$) - Obj($K$), compared with the theoretical upper bound $\delta(K,z,\gamma)$ on the worst-case constraint value $\bar{g}^N(x) - \bar{g}^K(x)$,  from Corollary~\ref{cor:maxaffine}. Solid lines are the difference, dotted lines are the upper bounds. Right: solve time. }%
	\label{fig:facility1}%
\end{figure}

\subsection{\reviewChanges{Newsvendor problem}}
\label{sub:news}
\reviewChanges{
We consider a 2-item newsvendor problem where, at the beginning of each day, the vendor orders $x\in \reals^2$ products at price $h=(4,5)$. These products will be sold at the prices $c=(5,6.5)$, until either the uncertain demand $u$ or inventory $x$ is exhausted.
The objective function to minimize is the sum of the ordering cost minus the revenue:
\begin{equation*}
	h^Tx - c^T\min\{x,u\},
\end{equation*}
from which we obtain the maximum-of-affine uncertain function $g$ to minimize,
$$g(u,x) = h^Tx + \max(-c_1x_1 - c_2 x_2, -c_1 x_1 -c_2 u_2, -c_1u_1 -c_2x_2, -c_1u_1 - c_2u_2 ).$$}
\reviewChanges{
We assume a polyhedral support $\supp= \{u \mid Cu\leq b\}$, and solve, with $p = 1$,
\begin{equation*}
	\label{eq:news}
	\begin{array}{ll}
	\text{minimize} & h^Tx  + \lambda\epsilon + \sum_{k=1}^K w_k s_{k} \\
	\text{subject to}  & -c^Tx + \gamma_{1k}^T(b - C\bar{d}_k)  \leq s_k, \quad k = 1,\dots, K\\
	& -c_1 x_1 - \bar{d}_k^T(c_2e_2) + \gamma_{2k}^T(b - C\bar{d}_k) \leq s_k , \quad k = 1,\dots, K\\
& -c_2 x_2 - \bar{d}_k^T(c_1e_1) + \gamma_{3k}^T(b - C\bar{d}_k) \leq s_k, \quad k = 1,\dots, K\\
& -\bar{d}_k^T(c) + \gamma_{4k}^T(b - C\bar{d}_k) \leq s_k, \quad k = 1,\dots, K\\
	 & \left\| -C^T\gamma_{1k} \right\|_2 \leq \lambda_{k}, \quad k = 1,\dots, K\\
	 & \left\| -C^T\gamma_{2k}+ c_1 e_1 \right\|_2 \leq \lambda_{k}, \quad k = 1,\dots, K\\
	 & \left\| -C^T\gamma_{3k}+ c_2 e_2\right\|_2 \leq \lambda_{k}, \quad k = 1,\dots, K\\
	 & \left\| -C^T\gamma_{4k}+ c \right\|_2 \leq \lambda_{k}, \quad k = 1,\dots, K\\
	 & \gamma_{jk} \geq 0, \quad j = 1,\dots,4, \quad k = 1,\dots, K\\
	 & x \geq 0
	\end{array}
\end{equation*}
We have variables $x \in \reals^{n}$,$s_{k} \in \reals$, $\lambda \in \reals$, $\gamma_{jk} \in \reals^{m}$, for $j = 1,\dots,4$ and $k= 1,\dots,K$.
The $\gamma$ variables arise from enforcing the support. We denote $e_1 = (1,0)$ and $e_2 = (0,1)$.}

\begin{figure}[h]%
	\centering
	\includegraphics[width=0.6\textwidth]{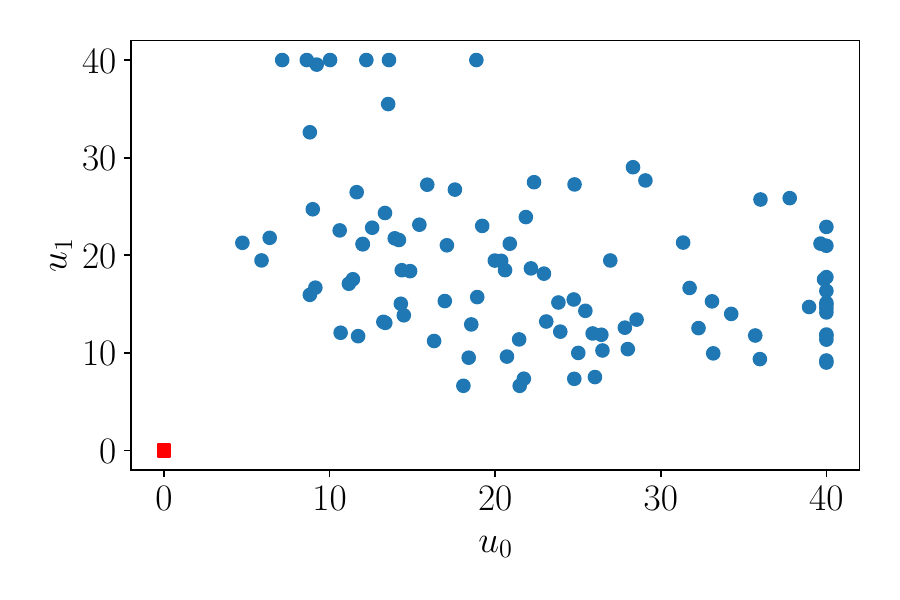}
	\caption{Newsvendor. Datapoints and the outlier at (0,0). }%
	\label{fig:news}%
\end{figure}

\reviewChanges{
For this problem, we consider the effects of outliers on the performance of \gls{MRO}. Therefore, we consider the data to have an outlier at $(0,0)$, the worst-case value of the support set. In Figure~\ref{fig:news}, we show a set of generated data along with this outlier point. }

\reviewChanges{
We consider three ways to solve the problem, decribed as follows. 
\begin{enumerate}
	\item \gls{MRO}, where we directly apply \gls{MRO} to the dataset with the outlier.
	\item ROB-\gls{MRO}, where we perform preliminary analysis on the dataset to remove the outlier point, then apply \gls{MRO} to the cleaned dataset. 
	\item AUG-\gls{MRO}, where we perform the clustering step on data without the outlier, then define an augmented distribution supported on $K+1$ points, where the extra point is the outlier point (0,0), with weight $1/N$. The weights of the other clusters are ajusted accordingly. 
\end{enumerate}}

\reviewChanges{
\paragraph{Problem setup.}%
To generate data, we set $N = 100$ data samples. We assume demand is supported between 0 and 40, which we write as $Cu\leq b$, where $C = [-I\; I]^T$ and $b = (0,0,0,0,40,40,40,40)$. We allow non-integer demand to allow for more variance in the data.
We generate the demand from a log-normal distribution, where the underlying normal distribution has parameters
\begin{equation*}
	\mu =
\left[ \begin{array}{l} 3.0 \\
2.8 \end{array} \right],\quad \Sigma =
\left[\begin{array}{ll} \phantom{-}0.3 & -0.1\\
-0.1 & \phantom{-}0.2 \end{array} \right],
\end{equation*}
and take the minimum between the generated values and 40. For the upper bound $\delta(K,z,\gamma)$ on  $\bar{g}^N(x) - \bar{g}^K(x)$ from Corollary~\ref{cor:maxaffine}, we have $(1/N)\sum_{k=1}^K \sum_{i \in C_k}  \max_{j\leq 4}((-\tilde{c}_j -C^T\gamma_{jk})^T(d_i - \bar{d}_k))$, where $\tilde{c}_1 = 0,~ \tilde{c}_2 = c_1e_1,~\tilde{c}_3 = c_2e_2,~\tilde{c}_4 = c$.}

\reviewChanges{
\paragraph{Choosing $K$}
Plotting the clustering value $D(K)$ over $K$, for the dataset both with and without the outlier, we note an elbow at around $K=5$, though not very prominent. The recommendation is setting $K$ around 5, to be fine tuned through cross-validation.}

\begin{figure}[h]%
	\centering
	\includegraphics[width=0.8\textwidth]{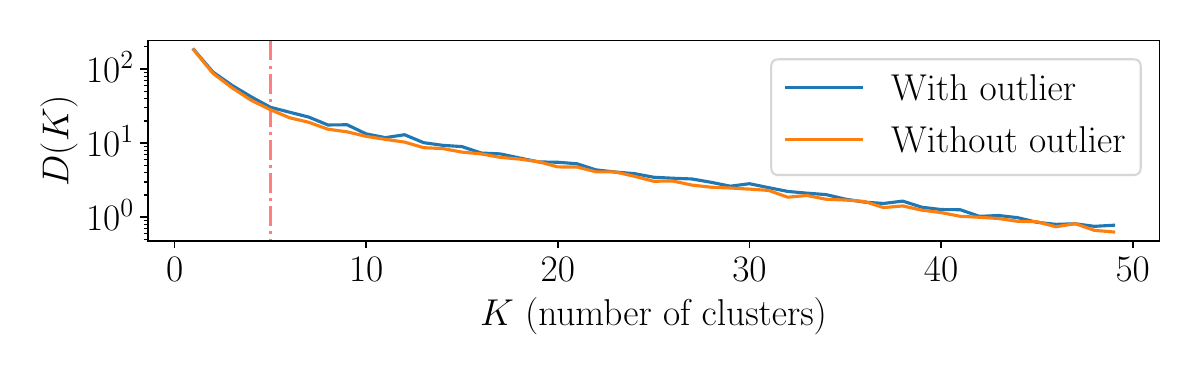}
	\caption{Newsvendor. $D(K)$ vs. $K$. Dotted red line: $K = 5$.}%
	\label{fig:newsk}%
\end{figure}
\begin{figure}[H]%
	\centering
	\includegraphics[width=\textwidth]{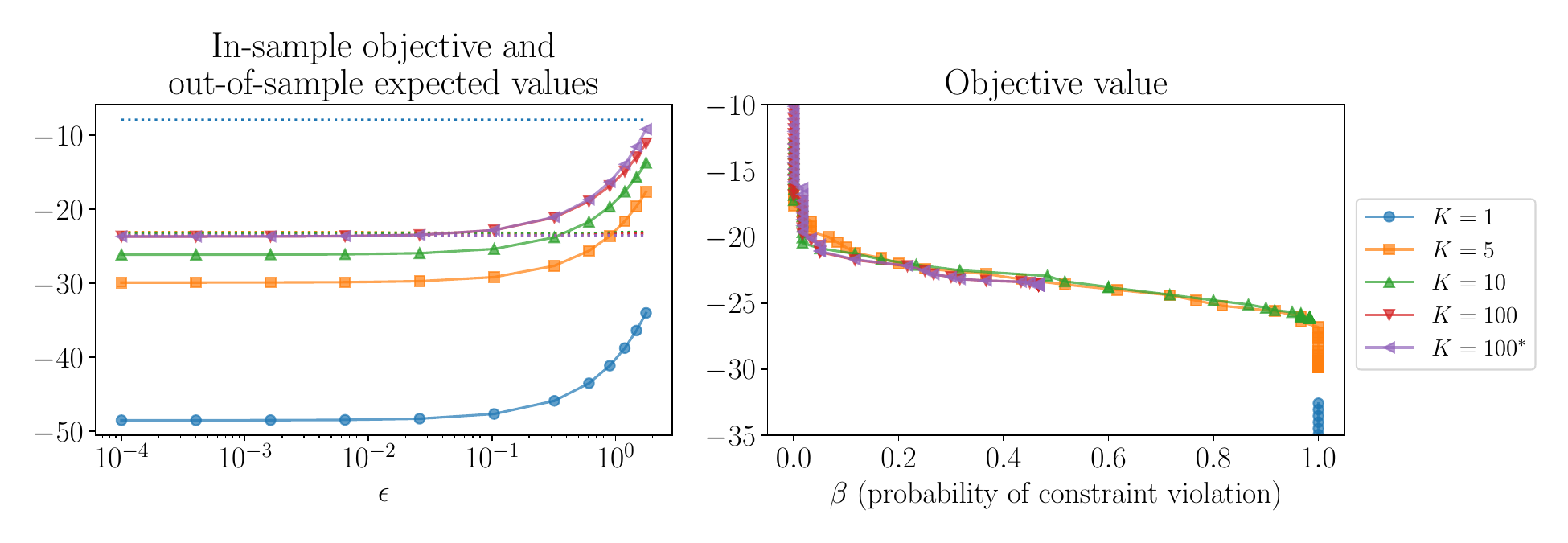}
	\includegraphics[width=\textwidth]{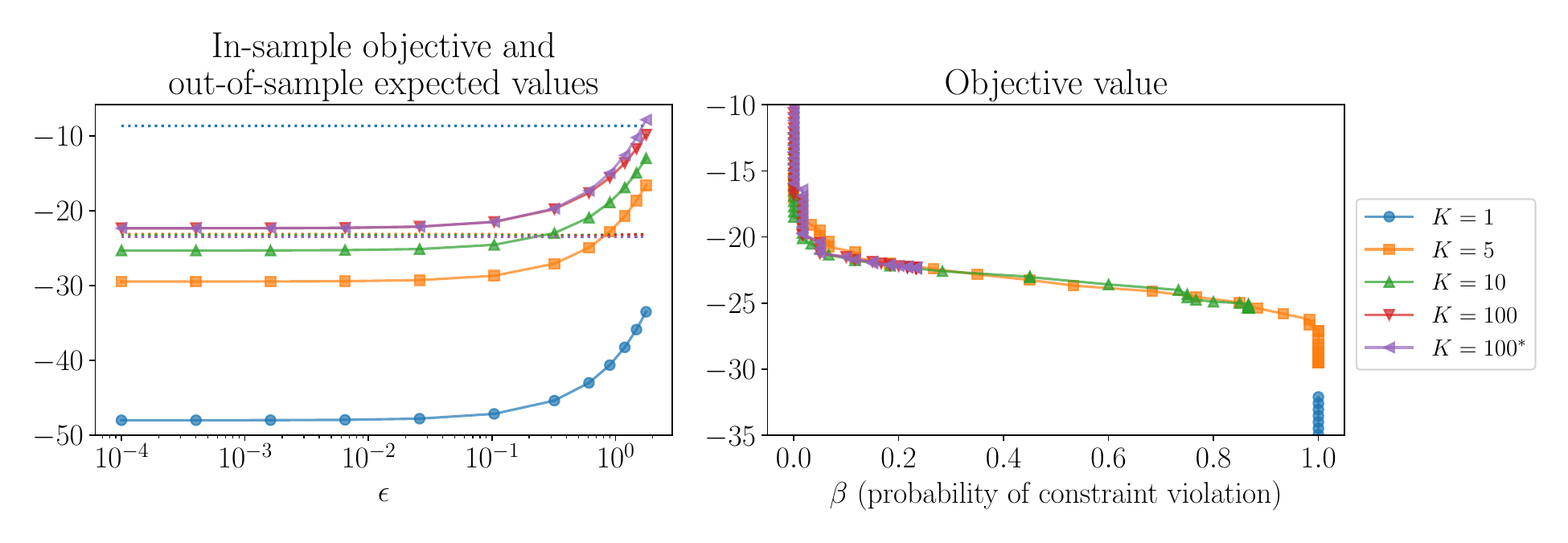}
	\includegraphics[width=\textwidth]{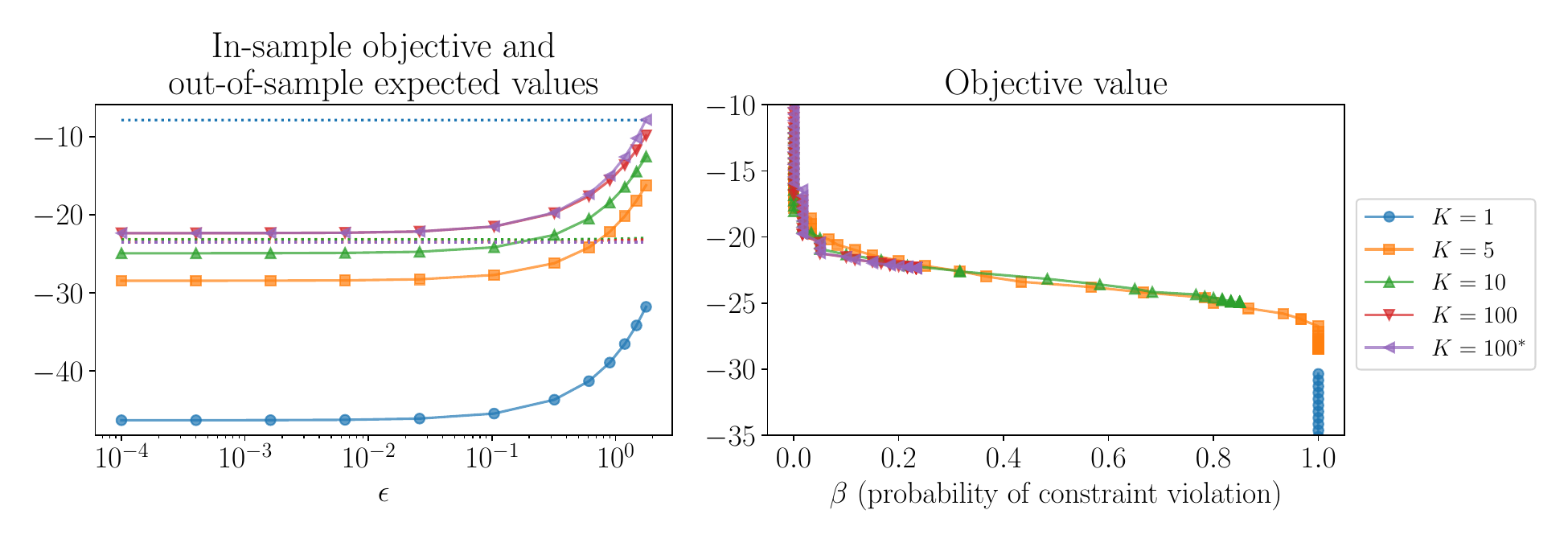}
	\caption{Newsvendor. Left: in-sample objective values vs $\epsilon$ for different $K$. Right: objective value vs $\beta$ for different $K$; each point represents the solution for the $\epsilon$ achieving the smallest objective value.  $K = 100^*$ is the formulation without the support constraint. Top: \gls{MRO}. Middle: ROB-\gls{MRO}. Bottom: AUG-\gls{MRO}.}%
	\label{fig:news2}%
\end{figure}
\begin{figure}[h]%
	\centering
	\includegraphics[width=\textwidth]{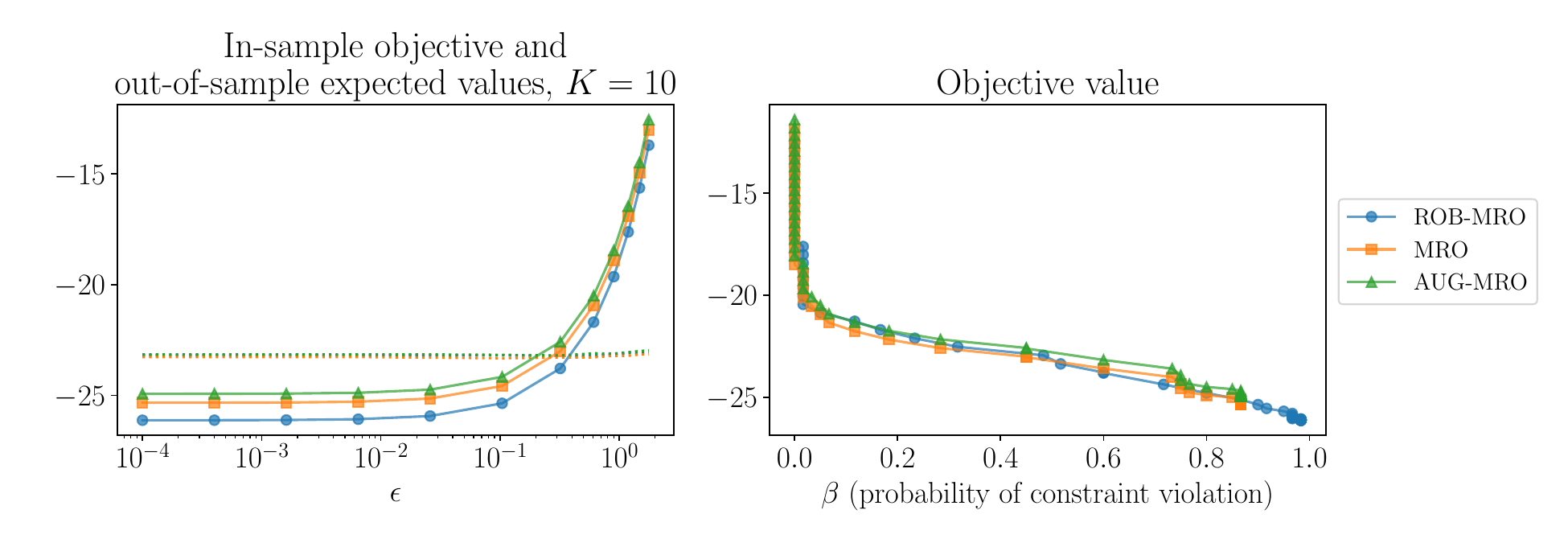}
	\caption{Newsvendor. Left: in-sample objective values vs $\epsilon$ for  $K=10$. Right: objective value vs $\beta$ for $K=10$.}%
	\label{fig:newscomp}%
\end{figure}
\begin{figure}[h]%
	\centering
	\includegraphics[width=\textwidth]{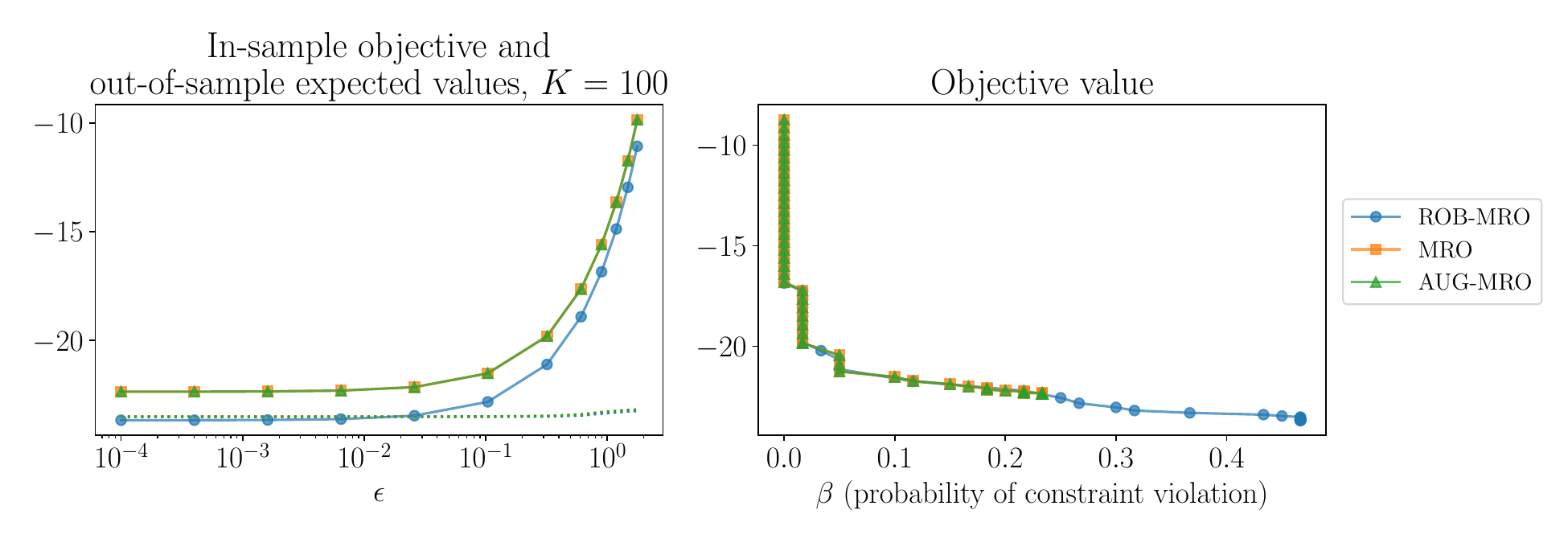}
	\caption{Newsvendor. Left: in-sample objective values vs $\epsilon$ for  $K=100$. Right: objective value vs $\beta$ for $K=100$.}%
	\label{fig:newscomp2}%
\end{figure}

\reviewChanges{
\paragraph{Results.}
In Figure~\ref{fig:news2}, we compare the in and out-of-sample objective values of the three methods, and note their similar performance. While setting $K=1$ yields suboptimal results, we note that for $K=5$ and above, we can achieve similar performance as setting $K=100$. To examine the effect of the outlier more closely, in Figures~\ref{fig:newscomp} and~\ref{fig:newscomp2}, we compare, for $K=10$ and $K=100$, the objectives and tradeoff curves for the three methods. We note that, when the outlier is averaged with other datapoints, the final in-sample objective may be improved, as the centroid moves closer to the non-outlier points. We observe this in Figure~\ref{fig:newscomp}, where \gls{MRO}, in which the outlier may be clustered with other points, offers a lower in-sample objective than AUG-\gls{MRO}, in which the outlier is considered its own cluster. \gls{MRO} has in fact offered protection against the outlier.  Lastly, as expected, ROB-\gls{MRO}, where the outlier point is removed, yields the best in-sample results. Regardless of the method, we note that the final out-of-sample tradeoff curves are near-identical. Comparing Figures~\ref{fig:newscomp} and~\ref{fig:newscomp2}, we note that the difference between \gls{MRO} and ROB-\gls{MRO} for $K=10$ is not larger than the difference for $K=N=100$, which shows that, while removing outliers before solving the problem may be helpful, the effect of outliers will not be worse for \gls{MRO} compared to classic Wasserstein \gls{DRO}.} 

\begin{figure}[h]%
	\centering
	\includegraphics[width=\textwidth]{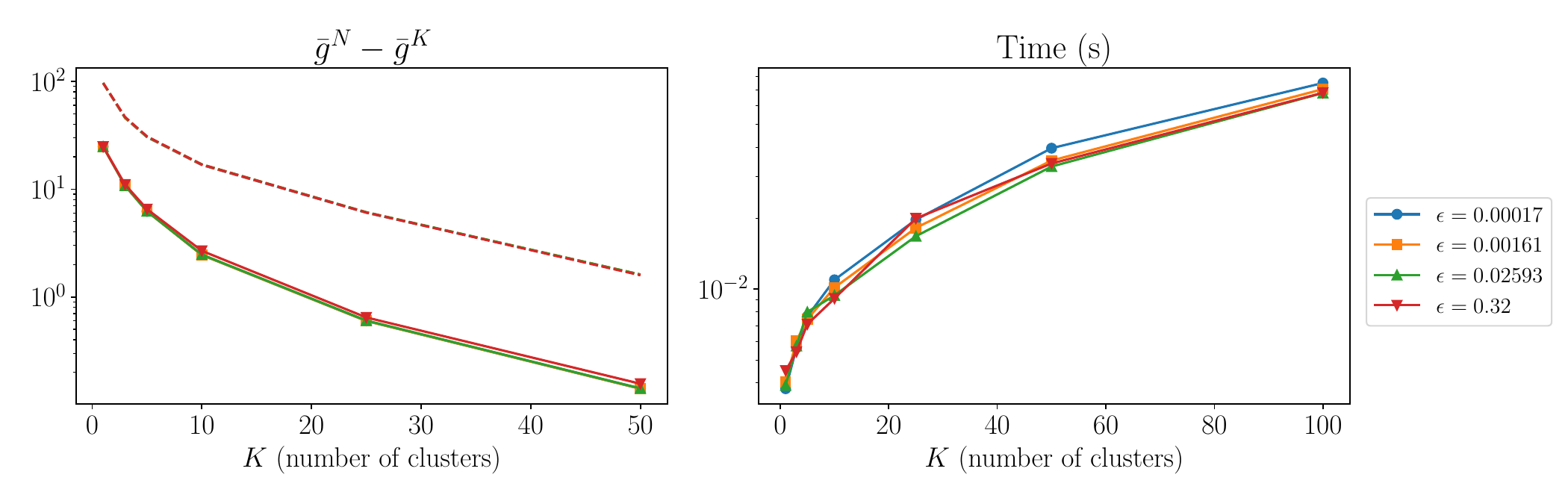}
	\caption{Newsvendor. Left: the difference in the value of the uncertain objective between using $N$ and $K$ clusters, calculated as $\bar{g}^N(x) - \bar{g}^K(x)$, compared with the theoretical upper bound $\delta(K,z,\gamma)$ from Corollary~\ref{cor:maxaffine}. Solid lines are the difference, dotted lines are the upper bounds. Right: solve time. }%
	\label{fig:news3}%
\end{figure}

\reviewChanges{We note in Figure~\ref{fig:news3} that the upper bound on $\bar{g}^N(x) - \bar{g}^K(x)$, given in Corollary~\ref{cor:maxaffine}, holds for \gls{MRO}.
We again note that bounds we observe do not depend on $\bar{g}^{N*}(x)$, so it is irrelevant whether or not the support has an affect on the worst-case constraint value.
Regardless, we see that the support only minimally affects the worst-case constraint value, at only at higher values of $\epsilon$.
Overall, choosing $K=5$, we obtain an order of magnitude computational speed-up.}

\section{Conclusions}
\label{sec:conclusions}
We have presented \glsentryfull{MRO}, a new data-driven methodology for decision-making under uncertainty that bridges robust and distributionally robust optimization while preserving rigorous probabilistic guarantees.
By clustering the dataset before performing \gls{MRO}, we solve an efficient and computationally tractable formulation with limited performance degradation.
In particular, we showed that when the constraints are affine in the uncertainty, clustering does not affect the optimal value of the objective.
\reviewChanges{When the constraint is concave or maximum-of-concave in the uncertainty, we directly quantified the change in worst-case constraint value that is caused by clustering. For problems with objective uncertainty, this directly bounds the change in the optimal value caused by clustering. }
We demonstrated this result through a set of numerical examples, where we observed the possibility of tuning the size of the uncertainty set such that using a small number of clusters achieves near-identical performance of traditional \gls{DRO}, with much higher computational efficiency.
\reviewChanges{In the final example, we also demonstrated that \gls{MRO} offers protection against outliers compared to Wasserstein \gls{DRO}.}


\ifpreprint
\section*{Acknowledgements}
\fi
\reviewChanges{Irina Wang and Bartolomeo Stellato are supported by the NSF CAREER Award ECCS 2239771.} The simulations presented in this article were performed on computational resources managed and supported by Princeton Research Computing, a consortium of groups including the Princeton Institute for Computational Science and Engineering (PICSciE) and the Office of Information Technology's High Performance Computing Center and Visualization Laboratory at Princeton University.

We would like to thank Daniel Kuhn for the useful feedback and for pointing us to the literature on scenario reduction techniques.

\bibliography{bibliography}

\ifpreprint
\appendix
\section{Appendices}
\else
\begin{appendices}
\fi

\ifpreprint
\subsection{Proof of the constraint reformulation in~\eqref{eq:robustopt_p}}
\fi

\label{app:dual_form}
To simplify notation, we define $c_k(v_k) \define   \| v_k - \bar{d}_k \|^p - \epsilon^p$. Then, starting from the inner optimization problem of~\eqref{eq:unif-const}:
\begin{equation*}
	\begin{aligned}
& \begin{cases}
	\underset{v_1 \dots v_K \in \supp}{\sup} &\sum_{k=1} ^K w_k g(v_k, x) \\
	\mbox{subject to} & \sum_{k=1} ^K w_k c_k(v_k) \le 0
\end{cases}\\
&= \begin{cases}
	\underset{v_1 \dots v_K\in \supp}{\sup}  & \underset{\lambda \geq 0}{\inf}   \sum_{k=1} ^K  w_k g(v_k, x)  - \lambda \sum_{k=1} ^K w_k c_k (v_k)\\
	\end{cases}\\
	&= \begin{cases}
		\underset{\lambda \geq 0}{\inf} & \underset{v_1 \dots v_K\in \supp}{\sup}    \sum_{k=1} ^K  w_k  g(v_k, x)  - \lambda \sum_{k=1} ^K w_k c_k (v_k)\\
		\end{cases}\\
&= \begin{cases}
	\underset{\lambda \geq 0}{\inf} &  \sum_{k=1} ^K  w_k s_k \\
	\mbox{subject to} &  \underset{v_k\in \supp}{\sup}\hspace{1mm} g(v_k, x) -\lambda c_k (v_k)  \le s_k \quad k = 1,\dots,K
\end{cases}\\
&= \begin{cases}
	\underset{\lambda \geq 0}{\inf} &  \sum_{k=1} ^K  w_k s_k \\
	\mbox{subject to} &  [-g +\chi_{\supp} + \lambda c_k]^*(0)  \le s_k \quad k = 1,\dots,K,
\end{cases}\\
\end{aligned}
\end{equation*}
where we used the Lagrangian in the first equality and Lagrangian strong duality in the second inequality, where we applied the assumption that $g$ is upper-semicontinuous and concave in $v$, the convexity of the supports, and the satisfaction of Slater's condition for all $\epsilon > 0$. For the last equality, we used the definition of conjugate functions. Now, borrowing results from Esfahani et al.~\cite[Theorem 4.2]{mohajerin_esfahani_data-driven_2018}, Rockafellar and Wets~\cite[Theorem 11.23(a), p. 493]{rockafellar_wets_1998}, and Zhen et al.~\cite[Lemma B.8]{zhen2021mathematical}, with regards to the conjugate functions of infimal convolutions and $p$-norm balls, we note that:
\begin{equation*}
	[-g + \chi_{\supp} + \lambda c_k]^* (0) = \underset{y_{k},z_{k}}{\inf} ([-g]^*(z_{k} - y_{k},x) + \sigma_{\supp}(y_{k}) + [\lambda c_k]^* (-z_{k})]),
\end{equation*}
and
\begin{equation*}
	[\lambda c_k]^* (-z_{k}) = \underset{v_k}{\sup}(-z_{k}^T v_k - \lambda \| v_k - \bar{d}_k\|^p + \lambda \epsilon^p) =-z_{k}^T\bar{d}_k + \phi(q)\lambda {\|}z_{k}/\lambda{\|}^q_* + \lambda \epsilon^p.
\end{equation*}
Substituting these in, we arrive at:
\begin{equation*}
	\begin{aligned}
\begin{cases}
	\underset{\lambda \geq 0, z_{k},y_{k}, s_k}{\inf} & \sum_k^K w_k s_k\\
	\mbox{subject to} & [-g]^*(z_{k} - y_{k},x) + \sigma_{\supp}(y_{k}) - z_{k}^T\bar{d}_k  + \phi(q)\lambda \left\|z_{k}/\lambda \right\|^q_* + \lambda \epsilon^p \le s_k \\
	& \hspace{7cm} \quad k = 1,\dots,K.
\end{cases}
\end{aligned}
\end{equation*}
\ifpreprint
\subsection{Proof of the constraint reformulation in~\eqref{eq:unif-dual-p}}
\fi
\label{app:dual_form-p}
Starting from the inner optimization problem of~\eqref{eq:unif-const-p}:
\begin{equation*}
	\begin{aligned}
& \begin{cases}
	\underset{v_1 \dots v_K \in \supp}{\sup} &\sum_{k=1} ^K w_k g(v_k, x) \\
	\mbox{subject to} &  \|v_k - \bar{d}_k\| \le \epsilon, \quad k = 1,\dots,K
\end{cases}\\
&= \begin{cases}
	\underset{v_1 \dots v_K\in \supp}{\sup}  & \underset{\lambda_k \geq 0}{\inf}   \sum_{k=1} ^K  w_k (g(v_k, x)  +  (1/w_k)\lambda_k (\epsilon  - \|v_k - \bar{d}_k\| ))\\
	\end{cases}\\
	&= \begin{cases}
		\underset{\lambda_k \geq 0}{\inf} & \underset{v_1 \dots v_K\in \supp}{\sup}    \sum_{k=1} ^K  w_k  (g(v_k, x)  +  (1/w_k)\lambda_k (\epsilon  - \|v_k - \bar{d}_k\| ))\\
		\end{cases}\\
&= \begin{cases}
	\underset{\lambda_k \geq 0}{\inf} &  \sum_{k=1} ^K  w_k s_k \\
	\mbox{subject to} &   (\lambda_k/w_k)\epsilon + \underset{v_k\in \supp}{\sup}\hspace{1mm} g(v_k, x) -(\lambda_k/w_k)\|v_k - \bar{d}_k\| \le s_k \quad k = 1,\dots,K
\end{cases}\\
&= \begin{cases}
	\underset{\lambda_k \geq 0}{\inf} &  \sum_{k=1} ^K  w_k s_k \\
	\mbox{subject to} & (\lambda_k/w_k)\epsilon + \underset{v_k\in \supp}{\sup}\hspace{1mm} g(v_k, x) -  \underset{\|z_k\|_* \leq \lambda_k/w_k}{\max} z_k^T(v_k - \bar{d}_k)\le s_k\\
	 &\hspace{7cm}\quad k = 1,\dots,K
\end{cases}\\
&= \begin{cases}
	\underset{\lambda_k \geq 0}{\inf} &  \sum_{k=1} ^K  w_k s_k \\
	\mbox{subject to} & (\lambda_k/w_k)\epsilon +  \underset{\|z_k\|_* \leq \lambda_k/w_k}{\min} \underset{v_k\in \supp}{\sup}\hspace{1mm} g(v_k, x) -   z_k^T(v_k - \bar{d}_k)\le s_k\\
	&\hspace{7cm}\quad k = 1,\dots,K
\end{cases}\\
&= \begin{cases}
	\underset{\lambda_k \geq 0, z_k}{\inf} &  \sum_{k=1} ^K  w_k s_k \\
	\mbox{subject to} & (\lambda_k/w_k)\epsilon + [-g +\chi_{\supp}]^*(-z_k) + z_k^T\bar{d}_k \le s_k \quad k = 1,\dots,K\\
	& \|z_k\|_* \leq \lambda_k/w_k,\quad k = 1,\dots,K
\end{cases}\\
&= \reviewChanges{\begin{cases}
	\underset{z_k}{\inf} &  \sum_{k=1} ^K  w_k s_k \\
	\mbox{subject to} &  [-g]^*(z_{k} - y_{k},x) + \sigma_{\supp}(y_{k})  - z_k^T\bar{d}_k  + \epsilon\|z_k\|_* \le s_k \quad k = 1,\dots,K.\\
\end{cases}}\\
\end{aligned}
\end{equation*}
We have again used the Lagragian, Lagrangian strong duality, and the definition of conjugate function. In particular, in the fourth equality, we refer to the proof of~\cite[Theorem 4.2]{mohajerin_esfahani_data-driven_2018} where we use the definition of the dual norm. \reviewChanges{In the final equality, we make the substition $z_k = -z_k$ and $ \lambda_k/w_k = \|z_k\|_* $, because the coefficient multiplying $\lambda_k$, $\epsilon/w_k$, is nonnegative and $\lambda_k/w_k$ achieves its lower bound, $\|z_k\|_*$, for all $k$.}

\ifpreprint
\subsection{Proof of the constraint reformulation in~\eqref{eq:robustopt_p_max}}
\fi

\label{app:dual_form_max}
\reviewChanges{
To simplify notation, we define $c_k(v_{jk}) \define   \| v_{jk} - \bar{d}_k \|^p - \epsilon^p$. Then, starting from the inner optimization problem of~\eqref{eq:unif-const_max}:
\begin{equation*}
	\begin{aligned}
& \begin{cases}
	\underset{v_{11}, \dots, v_{JK} \in \supp, \alpha \in \Gamma}{\sup} &\sum_{k=1}^K \sum_{j=1}^J \alpha_{jk} g_j(v_{jk}, x) \\
	\mbox{~~~~subject to} & \sum_{k=1} ^K \sum_{j=1}^J \alpha_{jk} c_k(v_{jk}) \le 0
\end{cases}\\
&=\begin{cases}
	\underset{v_{11}, \dots, v_{JK} \in \supp, \alpha \in \Gamma}{\sup} & \underset{\lambda \geq 0}{\inf}~\sum_{k=1}^K \sum_{j=1}^J \alpha_{jk} g_j(v_{jk}, x) - \lambda\sum_{k=1} ^K \sum_{j=1}^J \alpha_{jk} c_k(v_{jk})\\
\end{cases}\\
&=\begin{cases}
	\underset{\alpha \in \Gamma}{\sup}~\underset{\lambda \geq 0}{\inf}~\underset{v_{11}, \dots, v_{JK} \in \supp}{\sup} &\sum_{k=1}^K \sum_{j=1}^J \alpha_{jk} g_j(v_{jk}, x) - \lambda\sum_{k=1} ^K \sum_{j=1}^J \alpha_{jk} c_k(v_{jk}),\\
\end{cases}\\
\end{aligned}
\end{equation*}
We applied the Lagrangian in the first equality. Then, as the summation is over upper-semicontinuous functions $g_j(v_{jk},z)$ concave in $v_{jk}$, we can interchange the inf and the sup by Lagrangian strong duality. Next, we rewrite the formulation using an epigraph trick, and note that the $\alpha_{jk}$'s can be pulled outside the inner supremum due to their nonnegativity and the separability of the $v_{jk}$'s. 
\begin{equation*}
	\begin{aligned}
&= \begin{cases}
	\underset{\alpha \in \Gamma}{\sup}~\underset{\lambda \geq 0}{\inf} &\sum_{k=1}^K s_k  \\
	\mbox{subject to}~\underset{v_{11}, \dots, v_{JK} \in \supp}{\sup}& \sum_{j=1}^J \alpha_{jk}(g_j(v_{jk}, x) - \lambda c_k(v_{jk})) \le s_k \quad k = 1,\dots,K,
\end{cases}\\
&= \begin{cases}
	\underset{\alpha \in \Gamma}{\sup}~\underset{\lambda \geq 0}{\inf} &\sum_{k=1}^K s_k  \\
	\mbox{subject to}& \sum_{j=1}^J [-g_j' +\chi_{\supp'} + \lambda c_k']^*(0) \le s_k \quad k = 1,\dots,K.
\end{cases}\\
\end{aligned}
\end{equation*}
\reviewChanges{
Now, using the conjugate form $f^*(y) = \alpha g^*(y)$ of a right-scalar-multiplied function $f(x) = \alpha g(x/\alpha)$, we rewrite the above as 
\begin{equation*}
	\begin{aligned}
&= \begin{cases}
	\underset{\alpha \in \Gamma}{\sup}~\underset{\lambda \geq 0}{\inf} &\sum_{k=1}^K s_k  \\
	\mbox{subject to}& \sum_{j=1}^J \alpha_{jk}[-g_j +\chi_{\supp} + \lambda c_k]^*(0) \le s_k \quad k = 1,\dots,K.
\end{cases}\\
\end{aligned}
\end{equation*}
 Now, again using properties of conjuate functions as given in Appendix~\ref{app:dual_form}, we note that:
\begin{equation*}
	\begin{aligned}
		\alpha_{jk}[(-g_j +\chi_{\supp} + \lambda c_k)]^* (0) &= \alpha_{jk}\underset{y_{jk},z_{jk}}{\inf} ([-g_j]^*(z_{jk} - y_{jk},x) + \sigma_{\supp}(y_{jk}) + [\lambda c_k]^* (-z_{jk})).
	\end{aligned}
\end{equation*}
Substituting in the conjugate functions derived in Appendix~\ref{app:dual_form}, we have 
\begin{equation*}
	\begin{aligned}
\begin{cases}
	\underset{\alpha \in \Gamma}{\sup}~\underset{\lambda \geq 0, z_{jk},y_{jk}, s_k}{\inf} & \sum_k^K s_k\\
	\mbox{subject to} & \sum_{j=1}^J \alpha_{jk}([-g_j]^*(z_{jk} - y_{jk},x)\\
	& + \sigma_{\supp}(y_{jk}) - z_{jk}^T\bar{d}_k + \phi(q)\lambda \left\|z_{jk}/\lambda \right\|^q_* + \lambda \epsilon^p ) \le s_k \\
	&\hspace{2cm} \quad k = 1,\dots,K, \quad j = 1,\dots,J.
\end{cases}
\end{aligned}
\end{equation*}
Taking the supremum over $\alpha$, noting that $\sum_{j=1}^J \alpha_{jk} = w_k$ for all $k$, we arrive at
\begin{equation*}
	\begin{aligned}
\begin{cases}
	~\underset{\lambda \geq 0, z_{jk},y_{jk}, s_k}{\inf} & \sum_k^K s_k\\
	\mbox{subject to} & w_k ([-g_j]^*(z_{jk}' - y_{jk}',x) + \sigma_{\supp}(y_{jk}') - z_{jk}'^T\bar{d}_k \\
	&~~ + \phi(q)\lambda \left\|z_{jk}'/\lambda \right\|^q_* + \lambda \epsilon^p ) \le s_k \quad k = 1,\dots,K, \quad j = 1,\dots,J,
\end{cases}
\end{aligned}
\end{equation*}
which is equivalent to~\eqref{eq:robustopt_p_max}.
}

\ifpreprint
\subsection{Reformulation of the maximum-of-concave case for ${p = \infty}$}
\fi
\label{app:dual_form-p_max}

\reviewChanges{In the case where ${p = \infty}$, we have
\begin{equation}
	\label{eq:unif-const-p_max}
	\begin{array}{ll}
		\mbox{minimize} & f(x)\\
		\mbox{subject to} & \begin{dcases}
			\begin{rcases}
		\underset{v_{11},\dots,v_{JK} \in \supp, \alpha \in \Gamma}{\text{maximize}} &\quad \sum_{k=1} ^K\sum_{j=1} ^J  \alpha_{jk} g(v_{jk}, x) \\
		\mbox{~~~~subject to} & \quad  \sum_{j=1}^J (\alpha_{jk}/w_k)| v_{jk} - \bar{d}_k \|^p  \le \epsilon, \quad k = 1,\dots, K
			\end{rcases}
			\end{dcases} \le 0,
	\end{array}
\end{equation}
which has a reformulation where the constraint above is dualized,}
\reviewChanges{
	\begin{equation}
	\label{eq:unif-dual-p_max}
	\begin{array}{ll}
		\text{minimize} & f(x)\\
 \text{subject to} \quad & \sum_{k=1}^{K} w_k s_k \leq 0\\
&[-g_j]^*(z_{jk} - y_{jk},x) + \sigma_{{\supp}}(y_{jk}) - z_{jk}^T\bar{d}_k  + \lambda_k \epsilon \leq s_k\\
& \hspace{3cm} \quad k = 1,\dots, K, \quad j = 1,\dots, J\\
&  \left\|z_{jk} \right\|_* \leq \lambda_k \quad k = 1,\dots, K, \quad j = 1,\dots, J,
\end{array}
\end{equation}
with new variables $s_k \in \reals$, $z_{jk} \in \reals^m$, and $y_{jk} \in \reals^m$. }
\reviewChanges{
We prove this by starting from the inner optimization problem of~\eqref{eq:unif-const-p_max}:
\begin{equation*}
	\begin{aligned}
		& \begin{cases}
			\underset{v_{11}, \dots, v_{JK} \in \supp, \alpha \in \Gamma}{\sup} &\sum_{k=1}^K \sum_{j=1}^J \alpha_{jk} g_j(v_{jk}, x) \\
			\mbox{~~~~subject to} & \sum_{j=1}^J (\alpha_{jk}/w_k) \|v_{jk} - \bar{d}_k\| \le \epsilon, \quad k = 1,\dots,K
		\end{cases}\\
		&=\begin{cases}
			\underset{v_{11}, \dots, v_{JK} \in \supp, \alpha \in \Gamma}{\sup} & \underset{\lambda \geq 0}{\inf}~\sum_{k=1}^K (\sum_{j=1}^J \alpha_{jk} g_j(v_{jk}, x) + \lambda_k (\epsilon - \sum_{j=1}^J (\alpha_{jk}/w_k) \|v_{jk} - \bar{d}_k\|))\\
		\end{cases}\\
		&=\begin{cases}
			\underset{\alpha \in \Gamma}{\sup}~\underset{\lambda \geq 0}{\inf}~\underset{v_{11}, \dots, v_{JK} \in \supp}{\sup}&\sum_{k=1}^K (\sum_{j=1}^J \alpha_{jk} g_j(v_{jk}, x) + \lambda_k (\epsilon - \sum_{j=1}^J (\alpha_{jk}/w_k)  \|v_{jk} - \bar{d}_k\|))\\
		\end{cases}\\
	\end{aligned}
\end{equation*}
We have again formulated the Lagragian and noted that strong duality holds for the sum of concave functions in $v_{jk}$.
Now, using the same procedure as in Appendix~\ref{app:dual_form}, we can rewrite this in epigraph form, and use the dual norm definition,
		\begin{equation*}
			\begin{aligned}
				&= \begin{cases}
					\underset{\alpha \in \Gamma}{\sup}~\underset{\lambda \geq 0}{\inf} &\sum_{k=1}^K s_k  \\
					\mbox{subject to}&\underset{v_{11}, \dots, v_{JK} \in \supp}{\sup}~ \lambda_k \epsilon+ \sum_{j=1}^J \alpha_{jk}g_j(v_{jk}, x) - \lambda_k(\alpha_{jk}/w_k) \|v_{jk} - \bar{d}_k\|\le s_k\\
					& \hspace{4cm}  \quad k = 1,\dots,K,
				\end{cases}\\
				&= \begin{cases}
					\underset{\alpha \in \Gamma}{\sup}~\underset{\lambda \geq 0}{\inf} &\sum_{k=1}^K s_k  \\
					\mbox{subject to}&\underset{v_{11}, \dots, v_{JK} \in \supp}{\sup}~ \lambda_k \epsilon+ \sum_{j=1}^J \underset{\|z_{jk}\|_* \leq \lambda_k/w_k}{\min} \alpha_{jk}g_j(v_{jk}, x) \\
					& \hspace{4cm} - \alpha_{jk}z_{jk}^T(v_{jk} - \bar{d}_k)\le s_k \quad k = 1,\dots,K.
				\end{cases}\\
&= \begin{cases}
	\underset{\alpha \in \Gamma}{\sup}~\underset{\lambda_k \geq 0, z_{jk}}{\inf} &  \sum_{k=1} ^K s_k \\
	\mbox{subject to}&\lambda_k \epsilon+ \alpha_{jk}[-g_j +\chi_{\supp}]^*(-z_{jk},x) +\alpha_{jk}z_{jk}^T\bar{d}_k \le s_k \\
	&\hfill \quad k = 1,\dots,K, \quad j = 1,\dots,J\\
	& \|z_{jk}\|_* \leq \lambda_k/w_k \quad k = 1,\dots,K, \quad j = 1,\dots,J.
\end{cases}\\
\end{aligned}
\end{equation*}
Now, substituting $\lambda_k = \lambda_k w_k$, $z_{jk} = -z_{jk}$, and substituting in the conjugate functions derived in Appendix~\ref{app:dual_form_max}, we have 
\begin{equation*}
\begin{aligned}
&= \begin{cases}
	\underset{\alpha \in \Gamma}{\sup}~\underset{\lambda_k \geq 0, z_{jk}}{\inf} &  \sum_{k=1} ^K s_k \\
	\mbox{subject to} &  \lambda_k w_k\epsilon + \sum_{j=1}^J \alpha_{jk} ([-g_j]^*(z_{jk} - y_{jk},x) + \sigma_{\supp}(y_{jk}) \\
	&\hspace{2cm} - z_{jk}^T\bar{d}_k) \le s_k \quad k = 1,\dots,K, \quad j = 1,\dots,J\\
	& \|z_{jk}\|_*\leq \lambda_k,\quad k = 1,\dots,K, \quad j = 1,\dots,J.
\end{cases}\\
\end{aligned}
\end{equation*}
Note that rescaling $\lambda_k$ did not affect value of the problem, as minimizing $\lambda_k$ is equivalent to minimizing $\lambda_k w_k$.
Lastly, taking the supremum over $\alpha$, we arrive at
\begin{equation*}
	\begin{aligned}
	&= \begin{cases}
		\underset{\lambda_k \geq 0, z_{jk}}{\inf} &  \sum_{k=1} ^K s_k \\
		\mbox{subject to} & w_k(\lambda_k\epsilon + [-g_j]^*(z_{jk} - y_{jk},x) + \sigma_{\supp}(y_{jk})  - z_{jk}^T\bar{d}_k) \le s_k\\
		&\hfill  \quad k = 1,\dots,K, \quad j = 1,\dots,J\\
		& \|-z_{jk}\|_*\leq \lambda_k,\quad k = 1,\dots,K, \quad j = 1,\dots,J.
	\end{cases}\\
	\end{aligned}
	\end{equation*}
}

\ifpreprint
\subsection{Proof of the primal problem reformulation as $\mathbf{p \rightarrow \infty}$}
\fi
\label{app:primal_limit}

Consider again the function $\bar g^K$ discussed in Section \ref{sec:conservatism} and defined as
\[
\bar g^K(x; p)
\defn
\left\{
\begin{array}{l@{~~~~}l}
	\mbox{maximize} & \sum_{k=1}^K w_k g(v_k, x) \\[0.5em]
     \mbox{subject to} & \sum_{k=1}^K w_k \norm{v_k-d_k}^p\leq \epsilon^p\\[0.5em]
          & v_k\in \supp \quad k = 1,\dots,K
\end{array}
\right.
\]
where we make its dependence on $p$ explicit.
We have that $1/M<w_k=|C_k|/N<M$ for all $k= 1,\dots, K$ for some large $M\geq 1$.

\begin{theorem}
  \label{thm:extreme-case-infinity}
  Let the functions $\epsilon \mapsto g^\epsilon(d_k, x)$ be continuous for all $k= 1,\dots, K$ where
\(
  g^{\epsilon}(d, x)= \max \{g(v, x) \mid v\in \supp,\,\norm{v-d}\leq \epsilon\}.
\)
We have that
\(
  \bar g^K(x;\infty)\defn\lim_{p\to\infty} \bar g^K(x;p) = \sum_{k=1}^K w_k g^{\epsilon}(d_k, x).
  \)
\end{theorem}
\begin{proof}
  Using the auxiliary variables $t_k\geq 0$ for $k= 1,\dots, K$ we have that
  \[
    \bar g^K(x; p)
    =
    \left\{
      \begin{array}{l@{~~~~}l}
	      \mbox{maximize} & \sum_{k=1}^K w_k g(v_k, x) \\[0.5em]
        \mbox{subject to} & v_k\in \supp, ~ t_k\geq 0\quad k = 1,\dots,K\\[0.5em]
             &  \sum_{k=1}^K t_k^p \leq \epsilon^p\\[0.5em]
             &  t_k\geq \norm{v_k-d_k} w_k^{1/p} \quad k  =1,\dots,K.
      \end{array}
    \right.
  \]
  The function $\bar g^K(x;p)$ is hard to study directly. Hence, let us first introduce two auxiliary functions
  \[
    \bar{g}^K_u(x; p)
    =
    \left\{
      \begin{array}{l@{~~~~}l}
	      \mbox{maximize} & \sum_{k=1}^K w_k g(v_k, x) \\[0.5em]
        \mbox{subject to} & v_k\in \supp, \, t_k\geq 0 \quad k = 1,\dots,K\\[0.5em]
             &  \sum_{k=1}^K t_k^p \leq \epsilon^p\\[0.5em]
             &  t_k\geq \norm{v_k-d_k} M^{-1/p} \quad k = 1,\dots,K
      \end{array}
    \right.
  \]
  and
    \[
    \bar{g}^K_l(x; p)
    =
    \left\{
      \begin{array}{l@{~~~~}l}
	      \mbox{maximize} & \sum_{k=1}^K w_k g(v_k, x) \\[0.5em]
        \mbox{subject to} & v_k\in \supp, \, t_k\geq 0 \quad k = 1,\dots,K\\[0.5em]
             &  \sum_{k=1}^K t_k^p \leq \epsilon^p,\\[0.5em]
             &  t_k\geq \norm{v_k-d_k} M^{1/p} \quad k = 1,\dots,K.
      \end{array}
    \right.
  \]
  Observe that for $p\geq 1$ we have $1/M<w_k<M \implies M^{-1/p} < w^{1/p}_k < M^{1/p}$ for any $k\in 1,\dots, K$.
  As we hence have for all $k= 1,\dots, K$ that
  \(
    M^{-1/p} \norm{v_k-d_k} \leq w_k^{1/p} \norm{v_k-d_k}  \leq M^{1/p} \norm{v_k-d_k}
  \)
  we obtain the sandwich inequality $\bar{g}^K_l(x; p) \leq \bar{g}^K(x; p) \leq \bar{g}^K_u(x; p)$ for any $p\geq 1$.

  Furthermore, observe that when $t_k\geq 0$ for all $k=1,\dots, K$ then we have the implication $\sum_{k=1}^K t^p_k\leq \epsilon^p \implies \max_{k=1}^K t_k \leq \epsilon$. Hence, we have that
  \begin{align*}
    \bar{g}_u(x; p)
    \leq &
    \left\{
      \begin{array}{l@{~~~~}l}
				\mbox{maximize} & \sum_{k=1}^K w_k g(v_k, x) \\[0.5em]
        \mbox{subject to} & v_k\in \supp, \, t_k\geq 0 \quad k = 1,\dots,K,\\[0.5em]
             &  \max_{k=1}^K t_k \leq \epsilon,\\[0.5em]
             &  t_k\geq \norm{v_k-d_k} M^{-1/p} \quad k = 1,\dots,K
      \end{array}
           \right.\\
    = &    \left\{
      \begin{array}{l@{~~~~}l}
				\mbox{maximize} & \sum_{k=1}^K w_k g(v_k, x) \\[0.5em]
        \mbox{subject to} & v_k\in \supp\quad k = 1,\dots,K\\[0.5em]
             &  \max_{k=1}^K \norm{v_k-d_k} M^{-1/p} \leq \epsilon\\[0.5em]
      \end{array}
        \right.\\
    = &  \sum_{k=1}^K w_k \left[\max_{v\in \supp,\,\norm{v-d_k}\leq \epsilon M^{1/p}} g(v, x)\right].
  \end{align*}

  Similarly, observe that when $t_k\geq 0$ for all $k = 1,\dots, K$ we also have the inequality $\sum_{k=1}^K t_k^p \leq K (\max_{k=1}^K t_k)^p$. Hence, we have that
  \begin{align*}
    \bar{g}^K_l(x; p)
    \geq &
    \left\{
      \begin{array}{l@{~~~~}l}
				\mbox{maximize} & \sum_{k=1}^K w_k g(v_k, x) \\[0.5em]
        \mbox{subject to} & v_k\in \supp, \, t_k\geq 0 \quad k = 1,\dots,K\\[0.5em]
             &  K (\max_{k=1}^K t_k)^p \leq \epsilon^p\\[0.5em]
             &  t_k\geq \norm{v_k-d_k} M^{1/p} \quad k = 1,\dots,K
      \end{array}
           \right.\\
    = &
    \left\{
      \begin{array}{l@{~~~~}l}
				\mbox{maximize} & \sum_{k=1}^K w_k g(v_k, x) \\[0.5em]
        \mbox{subject to} & v_k\in \supp, \, t_k\geq 0 \quad k = 1,\dots,K,\\[0.5em]
             &   \max_{k=1}^K t_k \leq K^{-1/p} \epsilon,\\[0.5em]
             &  t_k\geq \norm{v_k-d_k} M^{1/p} \quad k = 1,\dots,K
      \end{array}
               \right.\\
    = &    \left\{
      \begin{array}{l@{~~~~}l}
				\mbox{maximize} & \sum_{k=1}^K w_k g(v_k, x) \\[0.5em]
        \mbox{subject to} & v_k\in \supp~~\forall k \in [1,\dots,K],\\[0.5em]
             &  \max_{k=1}^K \norm{v_k-d_k} M^{1/p} \leq \epsilon K^{-1/p}\\[0.5em]
      \end{array}
        \right.\\
    = &  \sum_{k=1}^K w_k \left[\max_{v\in \supp,\,\norm{v-d_k}\leq \epsilon (MK)^{-1/p}} g(v, x)\right].
  \end{align*}

  Finally, chaining all the inequalities together we obtain
  \[
    \sum_{k=1}^K w_k \left[\max_{v\in \supp,\,\norm{v-d_k}\leq \epsilon (MK)^{-1/p}} g(v, x)\right] \leq \bar{g}^K(x;p)  \leq \sum_{k=1}^K w_k \left[\max_{v\in \supp,\,\norm{v-d_k}\leq \epsilon M^{1/p}} g(v, x)\right]
  \]
  for any $p\geq 1$.
  Considering now the limit for $p$ to infinity
  \begin{align*}
    & \lim_{p\to\infty}\sum_{k=1}^K w_k  g^{\epsilon (MK)^{-1/p}}(d_k, x)\leq \lim_{p\to\infty} \bar{g}^K(x;p) \leq \lim_{p\to\infty} \sum_{k=1}^K w_k g^{\epsilon M^{1/p}}(d_k, x) \\
    \implies &\sum_{k=1}^K w_k \lim_{p\to\infty} g^{\epsilon (MK)^{-1/p}}(d_k, x)\leq \lim_{p\to\infty} \bar{g}^K(x;p) \leq \sum_{k=1}^K w_k  \lim_{p\to\infty}  g^{\epsilon M^{1/p}}(d_k, x) \\
    \implies &  \sum_{k=1}^K w_k g^{\epsilon}(d_k, x) \leq \lim_{p\to\infty} \bar{g}^K(x; p) \leq \sum_{k=1}^K w_k g^{\epsilon}(d_k, x)
  \end{align*}
  establishes the claim.
  The first implication follows from the fact that the finite sums and limits can be exchanged. The final implication follows from
  $\lim_{p\to \infty} (MK)^{-1/p} = \lim_{p\to \infty} M^{1/p}=1$ and the fact that the functions $\epsilon \mapsto g^\epsilon(d_k, x)$ are continuous  for all $k=1,\dots, K$.
\end{proof}

\ifpreprint
\subsection{Proof of the dual problem reformulation as $\mathbf{p \rightarrow \infty}$}
\fi
\label{app:dual_limit}

\begin{theorem}
  Let $S$ be a bounded set.
  Define here
  \begin{equation}
    \label{eq:def:barf}
    \bar g^K(x;\infty)\defn
    \left\{
      \begin{array}{l@{~~~~}l}
				\mbox{\rm{minimize}} & \sum_{k=1}^K w_k s_k\\
        \mbox{\rm{subject to}} &  \lambda \geq 0,~ z_{k}\in \Re^m, ~y_{k}\in \Re^m, ~s_k\in \Re^m \quad k=1,\dots, K,\\
             &  [-g]^*(z_{k} - y_{k}) + \sigma_{\supp}(y_{k}) - z_{k}^T{d}_k  + \epsilon \left\|z_{k} \right\|_* \le s_k\\
						 &\hspace{5.5cm} \quad k=1,\dots, K.
      \end{array}
    \right.
  \end{equation}
  Then,
  \(
    \lim_{p\to \infty} \bar g^K(x;p)= \bar g^K(x;\infty)
  \)
  for any $x\in X$.
\end{theorem}
\begin{proof}
  First, from Equation \eqref{eq:robustopt_p} we have for any $p> 1$ that
  \begin{align*}
    & \bar g^K(x;p)\\
    =&\reviewChanges{
                       \left\{
                       \begin{array}{l@{~~~~}l}
												\mbox{minimize} & \sum_{k=1}^K w_k s_k\\
                         \mbox{subject to} &  \lambda \geq 0,~ z_{k}\in \Re^m, ~y_{k}\in \Re^m, ~s_k\in \Re^m \quad k = 1,\dots,K,\\
                              &  [-g]^*(z_{k} - y_{k}) + \sigma_{\supp}(y_{k}) - z_{k}^T{d}_k  + \phi(q)\lambda \left\|z_{k}/\lambda \right\|^q_* + \lambda \epsilon^p\le s_k\\
															&\hspace{7cm} \quad k=1,\dots, K
                       \end{array}
                       \right.}\\
    \geq &
                       \left\{
                       \begin{array}{l@{~~~~}l}
												\mbox{minimize} & \sum_{k=1}^K w_k s_k\\
                         \mbox{subject to} &  \lambda_k \geq 0,~ z_{k}\in \Re^m, ~y_{k}\in \Re^m, ~s_k\in \Re^m \quad k = 1,\dots,K,\\
                              &  [-g]^*(z_{k} - y_{k}) + \sigma_{\supp}(y_{k}) - z_{k}^T{d}_k  + \phi(q)\lambda_k \left\|z_{k}/\lambda_k \right\|^q_* + \lambda_k \epsilon^p\le s_k\\
															&\hspace{7cm} \quad k=1,\dots, K
                       \end{array}
                       \right.\\
    = &
           \left\{
           \begin{array}{l@{~~~~}l}
						\mbox{minimize} & \sum_{k=1}^K w_k s_k\\
             \mbox{subject to} & z_{k}\in \Re^m, ~y_{k}\in \Re^m, ~s_k\in \Re^m \quad k = 1,\dots, K,\\
                  &  [-g]^*(z_{k} - y_{k}) + \sigma_{\supp}(y_{k}) - z_{k}^T{d}_k + \epsilon \left\|z_{k} \right\|_*\le s_k \quad k =1,\dots, K
           \end{array}
           \right.
  \end{align*}
  \reviewChanges{where the first equality is established in Appendix \ref{app:dual_form}} and the second equality follows from Lemma \ref{lemma:dual-helper}. Remark that the inequality in the second step simply follows as we introduce $\lambda_k$ and do not impose that $\lambda_k=\lambda$ for all $k=1,\dots, K$.
  Hence, considering the limit for $p$ tending to infinity gives us now
  \(
  \reviewChanges{\liminf_{p\to \infty}} \bar g^K(x;p)\geq \bar g^K(x;\infty).
  \)
  It remains to prove the reverse
  \(
  \reviewChanges{\limsup_{p\to \infty}} \bar g^K(x;p) \leq \bar g^K(x;\infty).
  \)

  Second, we have for any $p>1$ with $1/p+1/q=1$ that
  \begin{align*}
    & \bar g^K(x;p)\\
    \leq&
                       \left\{
                       \begin{array}{l@{~~~~}l}
												\mbox{minimize} & \sum_{k=1}^K w_k s_k\\
                         \mbox{subject to} &  z_{k}\in \Re^m, ~y_{k}\in \Re^m, ~s_k\in \Re^m \quad  k=1,\dots,K,\\
                              &  [-g]^*(z_{k} - y_{k}) + \sigma_{\supp}(y_{k}) - z_{k}^T{d}_k  \\
															& \qquad \qquad + \phi(q)\left[\frac{q-1}{q} \epsilon^{\frac{1}{1-q}} \max_{k'=1}^K\norm{z_{k'}}_*\right]^{1-q} \left\|z_{k} \right\|^q_* \\
                              & \qquad \qquad+ \left[\frac{q-1}{q} \epsilon^{\frac{1}{1-q}} \max_{k'=1}^K\norm{z_{k'}}_*\right] \epsilon^p\le s_k \quad  k=1,\dots, K\\
                              & (q-1)^{1/4} \leq \norm{z_k}_* \leq (q-1)^{-1/4}\quad  k =1,\dots, K
                       \end{array}
                     \right.\\
    =&
                       \left\{
                       \begin{array}{l@{~~~~}l}
												\mbox{minimize} & \sum_{k=1}^K w_k s_k\\
                         \mbox{subject to} &  z_{k}\in \Re^m, ~y_{k}\in \Re^m, ~s_k\in \Re^m \quad k= 1,\dots, K,\\
                              &  [-g]^*(z_{k} - y_{k}) + \sigma_{\supp}(y_{k}) - z_{k}^T{d}_k  + \frac{1}{q}\epsilon \left[\max_{k'=1}^K\norm{z_{k'}}_*\right]^{1-q} \left\|z_{k} \right\|^q_* \\
                              & \qquad + \frac{q-1}{q} \epsilon \max_{k'=1}^K\norm{z_{k'}}_*\le s_k \quad  k =1,\dots, K\\
			      & (q-1)^{1/4} \leq \norm{z_k}_* \leq (q-1)^{-1/4}\quad  k =1,\dots, K
                       \end{array}
                       \right.\\
    =&
                       \left\{
                       \begin{array}{l@{~~~~}l}
												\mbox{minimize} & \sum_{k=1}^K w_k s_k\\
                         \mbox{subject to} &  z_{k}\in \Re^m, ~y_{k}\in \Re^m, ~s_k\in \Re^m \quad k= 1,\dots, K,\\
                              &  [-g]^*(z_{k} - y_{k}) + \sigma_{\supp}(y_{k}) - z_{k}^T{d}_k  \\
                              & \qquad +\epsilon \left\|z_{k} \right\|_* \left[\frac{1}{q} \left[\max_{k'=1}^K\frac{\norm{z_{k'}}_*}{\left\|z_{k} \right\|_*}\right]^{1-q} +\frac{q-1}{q} \max_{k'=1}^K\frac{\norm{z_{k'}}_*}{\left\|z_{k} \right\|_*} \right] \le s_k\\
															&\hspace{7cm} \quad k =1,\dots, K\\
                         & (q-1)^{1/4} \leq \norm{z_k}_* \leq (q-1)^{-1/4}\quad k=1,\dots, K
                       \end{array}
                       \right.\\
    \leq&
                       \left\{
                       \begin{array}{l@{~~~~}l}
												\mbox{minimize} & \sum_{k=1}^K w_k s_k\\
                         \mbox{subject to} &  z_{k}\in \Re^m, ~y_{k}\in \Re^m, ~s_k\in \Re^m \quad k=1,\dots, K,\\
                              &  [-g]^*(z_{k} - y_{k}) + \sigma_{\supp}(y_{k}) - z_{k}^T{d}_k  +\epsilon \left\|z_{k} \right\|_* D(q) \le s_k \quad k =1,\dots, K\\
                         & (q-1)^{1/4} \leq \norm{z_k}_* \leq (q-1)^{-1/4}\quad  k=1,\dots, K.
                       \end{array}
                       \right.
  \end{align*}
  \reviewChanges{
    The first inequality follows from the choice $\lambda_k = \left[\frac{q-1}{q} \epsilon^{\frac{1}{1-q}} \max_{k'=1}^K\norm{z_{k'}}_*\right]$ and by imposing the restrictions $(q-1)^{1/4} \leq \norm{z_k}_* \leq (q-1)^{-1/4}$ for all $k =1,\dots, K$.
    The first equality follows from the identities
    \[
      \phi(q)\left[\frac{q-1}{q}\right]^{1-q} = (q-1)^{(q-1)}/q^q \left[\frac{q-1}{q}\right]^{1-q} = 1/q
    \]
    and
    \[
      p+\frac{1}{1-q} = \frac{1}{\frac{1}{p}} + \frac{1}{1-q} = \frac{1}{1-\frac{1}{q}}+\frac{1}{1-q} = \frac{-q}{-q+1}+\frac{1}{1-q}=\frac{1-q}{1-q}=1.
    \]
    The second equality follows by pulling $\epsilon \left\|z_{k} \right\|_*>0$ out of the last two terms in the constraints.
  }
  To establish the last inequality we note that $\max_{k'=1}^K \tfrac{\norm{z_{k'}}_*}{\left\|z_{k} \right\|_*} \geq \tfrac{\norm{z_{k}}_*}{\left\|z_{k} \right\|_*} = 1$ and $\max_{k'=1}^K \tfrac{\norm{z_{k'}}_*}{\left\|z_{k} \right\|_*} \leq (q-1)^{-1/2}$ and hence we can apply Lemma \ref{lemma:dual-helper-II}. Let
  \begin{equation}
    \label{eq:def:barf2}
    \bar g_u^K(x;p)\defn \left\{
      \begin{array}{l@{~}l}
        \mbox{minimize} & \sum_{k=1}^K w_k s_k\\
        \mbox{subject to} &  z_{k}\in \Re^m, ~y_{k}\in \Re^m, ~s_k\in \Re^m \quad k =1,\dots, K,\\
             &  [-g]^*(z_{k} - y_{k}) + \sigma_{\supp}(y_{k}) - z_{k}^T{d}_k  +\epsilon \left\|z_{k} \right\|_* D\left(\frac{p}{p-1}\right) \le s_k\\
						 &\hspace{6cm} \quad k =1,\dots, K\\
             & (p-1)^{-1/4} \leq \norm{z_k}_* \leq (p-1)^{1/4}\quad  k =1,\dots, K.
      \end{array}
    \right.
  \end{equation}
  Hence, as $q=p/(p-1)$ and $q-1=1/(p-1)$  we have
  \(
    \bar g^K(x;p) \leq \bar g_u^K(x;p)
  \)
  for all $p>1$.
  Hence, taking the limit $p\to\infty$ we have $\bar g^K(x;\infty)=\lim_{p\to\infty} \bar g^K(x;p) \leq \liminf_{p\to\infty} \bar g_u^K(x;p)$.
  \reviewChanges{In fact, as the function $D\left(\frac{p}{p-1}\right)$ defined in Lemma \ref{lemma:dual-helper-II} is nonincreasing for all $p$ sufficiently large this implies that $\bar g^K(x;p)$ is nonincreasing for $p$ sufficiently large and hence we have $\bar g^K(x;\infty)=\lim_{p\to\infty} \bar g^K(x;p) \leq \liminf_{p\to\infty} \bar g_u^K(x;p) = \lim_{p\to\infty} \bar g_u^K(x;p)$. }
    We now prove here that $\lim_{p\to\infty} \bar g_u^K(x;p) = \liminf_{p\to\infty} \bar g_u^K(x;p)\leq \bar g^K(x;\infty)$. \reviewChanges{Consider any feasible sequence $\{(z^t_{k}, \,y^t_{k}, \,s^t_k=[-g]^*(z^t_{k} - y^t_{k}) + \sigma_{\supp}(y^t_{k}) - (z^t_{k})^T{d}_k  +\epsilon \left\|z^t_{k} \right\|_* )\}_{t\geq 1}$ in the optimization problem characterizing $\bar g^K(x;\infty)$ in Equation \eqref{eq:def:barf} so that $\lim_{t\to\infty}\sum_{k=1}^K w_k s^t_k=\bar g^K(x;\infty)$.}
  Let $\tilde z^t_{k}\in \arg\max \set{\norm{z}_*}{z\in \Re^m, \, \norm{z-z_k^t}_*\leq 1/t}$ for all $t\geq 1$ and $k=1,\dots, K$ and \reviewChanges{observe that $\norm{\tilde z^t_k}_*= 1/t+\norm{z_k^t}_*\geq 1/t$.}
  Consider now an increasing sequence $\{p_t\}_{t\geq 1}$ so that $(p_t-1)^{1/4}\geq \max^K_{k=1} \norm{\tilde z^t_k}_*$ and $(p_t-1)^{-1/4}\leq 1/t$.
  Finally observe that the auxiliary sequence $\{(\tilde z^t_{k}, \,\tilde y^t_{k}=y^t_k + (\tilde z^t_{k}-z^t_{k}), \,\tilde s^t_k= [-g]^*(\tilde z^t_{k} - \tilde y^t_{k}) + \sigma_{\supp}(\tilde y^t_{k}) - (\tilde z^t_{k})^T{d}_k  +\epsilon \left\|\tilde z^t_{k} \right\|_* D\left(\tfrac{p_t}{(p_t-1)}\right))\}_{t\geq 1}$ is by construction feasible in the minimization problem characterizing the function $\bar g_u^K(x;p_t)$ in Equation \eqref{eq:def:barf2}. Hence, finally, we have
  \begin{align*}
    \lim_{p\to\infty} &g_u^K(x;p) = \lim_{t\to\infty} g_u^K(x;p_t)\\
    =  & \lim_{t\to\infty} \sum_{k=1}^K w_k \tilde s_k^t\\
    \reviewChanges{=} & \lim_{t\to\infty} \sum_{k=1}^Kw_k\left([-g]^*(\tilde z^t_{k} - \tilde y^t_{k}) + \sigma_{\supp}(\tilde y^t_{k}) - (\tilde z^t_{k})^T{d}_k  +\epsilon \left\|\tilde z^t_{k} \right\|_* D\left(\tfrac{p_t}{(p_t-1)}\right)\right)\\
    \leq &\lim_{t\to\infty} \sum_{k=1}^Kw_k\left([-g]^*(z^t_{k} - y^t_{k}) + \sigma_{\supp}(y^t_{k}) - (z^t_{k})^T{d}_k  +\epsilon \left\|z^t_{k} \right\|_* D\left(\tfrac{p_t}{(p_t-1)}\right)\right)\\
                                            & \qquad  + \sum_{k=1}^Kw_k\left(\max_{s\in S} \norm{s}+\norm{d_k}+\epsilon D\left(\tfrac{p_t}{(p_t-1)}\right) \right)/t\\
    \leq &\lim_{t\to\infty} \sum_{k=1}^Kw_k\left([-g]^*(z^t_{k} - y^t_{k}) + \sigma_{\supp}(y^t_{k}) - (z^t_{k})^T{d}_k  +\epsilon \left\|z^t_{k} \right\|_* D\left(\tfrac{p_t}{(p_t-1)}\right)\right)\\
    \reviewChanges{=} &\lim_{t\to\infty} \sum_{k=1}^Kw_k([-g]^*(z^t_{k} - y^t_{k}) + \sigma_{\supp}(y^t_{k}) - (z^t_{k})^T{d}_k  \\
		& \qquad +\epsilon \left\|z^t_{k} \right\|_* + \epsilon \left\|z^t_{k} \right\|_* \left(D\left(\tfrac{p_t}{(p_t-1)}\right)-1\right))\\
    \leq &\lim_{t\to\infty} \sum_{k=1}^Kw_k([-g]^*(z^t_{k} - y^t_{k}) + \sigma_{\supp}(y^t_{k}) - (z^t_{k})^T{d}_k \\
		& \qquad +\epsilon \left\|z^t_{k} \right\|_* + \epsilon  (p_t-1)^{1/4} \left(D\left(\tfrac{p_t}{(p_t-1)}\right)-1\right))\\
    \leq & \lim_{t\to\infty} \sum_{k=1}^Kw_ks_k = \bar g^K(x;\infty).
  \end{align*}  
  To establish the third inequality observe first that $-(\tilde z^t_{k})^T{d}_k = -(z^t_{k})^T{d}_k - (\tilde z^t_k-z^t_k)^Td_k\leq -(z^t_{k})^T{d}_k + \norm{\tilde z^t_k-z^t_k}_* \norm{d_k}\leq -(z^t_{k})^T{d}_k + \norm{d_k}/t$. \reviewChanges{Second, remark that we have 
    \begin{align*}
      \sigma_{\supp}(\tilde y^t_{k}) = & \sigma_{\supp}(y^t_{k}+(\tilde z^t_k-z_k^t)) \\
      \leq & \max_{s\in S} s^T(y^t_{k}+(\tilde z^t_k-z_k^t)) \\
      \leq & \max_{s\in S} s^Ty^t_{k}+\max_{s\in S} s^T(\tilde z^t_k-z_k^t)\leq \norm{s}\norm{\tilde z^t_k-z_k^t}_* \leq 1/t \max_{s\in S} \norm{s} .
    \end{align*}
    as $\norm{\tilde z_t-z_t}\leq 1/t$.
  }
  Lemma \ref{lemma:dual-helper-II} guarantees that $\lim_{t\to\infty} D\left(\tfrac{p_t}{(p_t-1)}\right)=1$.
  Finally, $\left\|z^t_{k} \right\|_*\leq \left\|\tilde z^t_{k} \right\|_*\leq (p_t-1)^{1/4}$ and
  \begin{align*}
    \lim_{t\to\infty} (p_t-1)^{1/4} \left(D\left(\tfrac{p_t}{(p_t-1)}\right)-1\right)=& \lim_{p\to\infty} (p-1)^{1/4} \left(D\left(\tfrac{p}{(p-1)}\right)-1\right)\\
    = & \lim_{q\to 1} (q-1)^{-1/4} \left( D(q)-1 \right)=0
  \end{align*}
  with $1/p+1/q=1$ using again Lemma \ref{lemma:dual-helper-II}.
\end{proof}

\begin{lemma}
  \label{lemma:dual-helper}
  We have
  \[
    \min_{\lambda \geq 0} ~\phi(q)\lambda \left\|z/\lambda \right\|^q_* + \lambda \epsilon^p = \norm{z}_*\epsilon
  \]
  for any $p>1$ and $q>1$ for which $1/p+1/q=1$, $\phi(q)=(q-1)^{q-1}/q^q$ and $\epsilon>0$.
\end{lemma}
\begin{proof}
  Remark that as the objective function $\lambda\mapsto \phi(q)\lambda \left\|z/\lambda \right\|^q_* + \lambda \epsilon^p$ is continuous and we have $\lim_{\lambda\to 0} \phi(q)\lambda \left\|z/\lambda \right\|^q_* + \lambda \epsilon^p = \lim_{\lambda\to\infty} \phi(q)\lambda \left\|z/\lambda \right\|^q_* + \lambda \epsilon^p=\infty$ as $\epsilon>0$ there must exist a minimizer $\lambda^\star\in \min_{\lambda \geq 0}\, \phi(q)\lambda \left\|z/\lambda \right\|^q_* + \lambda \epsilon^p$ with $\lambda_\star>0$.
The necessary and sufficient first-order convex optimality conditions of the minimization problem guarantee
\begin{align*}
  & \lambda^\star\in \min_{\lambda \geq 0} ~\phi(q)\lambda \left\|z/\lambda \right\|^q_* + \lambda \epsilon^p\\
  \iff & (1-q)\phi(q) \lambda_\star^{-q} \norm{z}_\star^q + \epsilon^p = 0\\
  \iff & \epsilon^p = (q-1) \phi(q) \lambda_\star^{-q} \norm{z}_\star^q \\
  \iff & \lambda_\star = \left[(q-1) \phi(q)\right]^{1/q} \norm{z}_\star\epsilon^{-p/q}\\
  \iff & \lambda_\star = \frac{q-1}{q} \epsilon^{\frac{1}{1-q}} \norm{z}_\star
\end{align*}
where we exploit that $1/p+1/q=1$ and $\phi(q)=(q-1)^{q-1}/q^q$. Indeed, we have
\[
  \left[(q-1) \phi(q)\right]^{1/q} = \left[(q-1)^{q}/q^q\right]^{1/q}= (q-1)/q
\]
and
\[
  -\frac{p}{q} = -\frac{1}{\frac{1}{p} q} = -\frac{1}{(1-1/q) q}=-\frac{1}{q-1}=\frac{1}{1-q}.
\]
Hence, we have
\begin{align*}
  & \min_{\lambda\geq 0} ~\phi(q) \lambda^{1-q} \norm{z}_\star^q+\lambda\epsilon^p\\
  = ~&\phi(q) \lambda_\star^{1-q} \norm{z}_\star^q+\lambda_\star\epsilon^p\\
  = ~&\phi(q) \left[ \frac{(q-1)^{1-q}}{q^{1-q}} \epsilon \norm{z}^{1-q}_\star \right] \norm{z}_\star^q +\left[\frac{q-1}{q} \epsilon^{\frac{1}{1-q}} \norm{z}_\star\right]\epsilon^p\\
  = ~&\phi(q) \frac{(q-1)^{1-q}}{q^{1-q}} \epsilon \norm{z}_\star +\frac{q-1}{q} \epsilon^{p+\frac{1}{1-q}} \norm{z}_\star\\
  = ~&\phi(q) \frac{(q-1)^{1-q}}{q^{1-q}} \epsilon \norm{z}_\star +\frac{q-1}{q} \epsilon \norm{z}_\star\\
  = ~&\frac{(q-1)^{q-1}}{q^q} \frac{(q-1)^{1-q}}{q^{1-q}} \epsilon \norm{z}_\star +\frac{q-1}{q} \epsilon \norm{z}_\star\\
  = ~&\frac{1}{q} \epsilon \norm{z}_\star +\frac{q-1}{q} \epsilon \norm{z}_\star\\
  = ~&\left[\frac{1}{q}  +\frac{q-1}{q} \right] \epsilon \norm{z}_\star\\
  = ~ & \epsilon \norm{z}_\star
\end{align*}
where we exploit that $1/p+1/q=1$ and $\phi(q)=(q-1)^{q-1}/q^q$. Indeed, we have
\[
  p+\frac{1}{1-q} = \frac{1}{\frac{1}{p}} + \frac{1}{1-q} = \frac{1}{1-\frac{1}{q}}+\frac{1}{1-q} = \frac{-q}{-q+1}+\frac{1}{1-q}=\frac{1-q}{1-q}=1
\]
establishing the claim.
\end{proof}

\begin{lemma}
  \label{lemma:dual-helper-II}
  Let $q>1$ then
  \[
    \max_{t\in [1, \tfrac{1}{\sqrt{q-1}}]} ~\frac{1}{q} t^{1-q} +\frac{q-1}{q} t = D(q)\defn \max\left(1, \frac{1}{q} \frac{1}{(q-1)^{(1-q)/2}} +\frac{\sqrt{q-1}}{q}\right)
  \]
  with $\lim_{q\to 1} D(q)=1$ and $\lim_{q\to 1} (q-1)^{1/4} \left(D(q)-1\right)=0$.
\end{lemma}
\begin{proof}
  Observe that the objective function is convex in $t$. Convex functions attain their maximum on the extreme points of their domain. The limits can be verified using standard manipulations.
\end{proof}

\ifpreprint
\subsection{Proof of the equivalence between $\mathbf{p =1}$ and  $\mathbf{p = \infty}$ for affine uncertainty}
\fi
\label{app:equivalence}
We again consider the single affine constraint~\eqref{eq:affineg}. In formulation~\eqref{eq:example-affine-p}, when $p=1$, we observe from Kuhn et al.~\citep[Section 2.2 Remark 1]{kuhn2019wasserstein} that
\begin{equation*}
	 \lim_{q \rightarrow \infty}\phi(q)\lambda \left\|z_{k}/\lambda \right\|^q_* =
	\begin{dcases} 0 & \text{if}\; \|z_{k} \| \leq \lambda \\
		\infty & \text{otherwise},
	\end{dcases}
	\end{equation*}
so when the support is $\supp = \reals^m$, the reformulation becomes
\begin{equation}
	\begin{array}{ll}
		\mbox{minimize} & f(x)\\
		\mbox{subject to} & a^Tx - b +\lambda \epsilon + (P^Tx)^T\sum_{k=1}^{K} w_k \bar{d}_k \le 0\\
		& \|P^Tx \|_* \leq \lambda \\
		& \lambda \ge 0,\\
	\end{array}
\end{equation}
where we can make the substitution $\lambda = \| P^Tx \|_*$, in which case this becomes equivalent to~\eqref{eq:example-affine-inf}.

\ifpreprint
\subsection{Proof of the convex reduction of the worst-case problem~\eqref{eq:dro-reduct}}
\fi
\label{app:convex_red}
We assume all preconditions given in Section 2.2. Referencing the proof for the case where $p=1$ in~\cite{mohajerin_esfahani_data-driven_2018}, we first expand-out the definition of the expected value and the Wasserstein-ball constraint. Then, we replace the joint distribution by a conditional one, since one of the distributions is the known empirical distribution, given by data. We use $K$ instead of $N$ and $w_i$ instead of $1/N$ such that this generalizes to ambiguitity sets defined as the Wasserstein-ball around the weighted empirical distribution $\prob^K$ of the clustered and averaged dataset.
\begin{equation*}
\begin{aligned}
  \sup_{\mathbf{Q}\in \mathbf{B}_\epsilon^p(\hat{\prob}^K)} \Expect^\mathbf{Q}[g(u,x)] &=
    \begin{cases}
      \underset{\Pi, \mathbf{Q}}{\sup} & \int_\reviewChanges{\supp} g(u,x)\mathbf{Q}(\dif u)\\
			\text{s.t.} & \int_\reviewChanges{\supp} \|u - u'\|^p \Pi (\dif u,\dif u') \leq \epsilon^p
    \end{cases} \\
		&=   \begin{cases}
      \underset{\mathbf{Q_k} \in \mathcal{M}(\supp)}{\sup} & \sum_{k = 1}^K w_k \int_\reviewChanges{\supp} g(u,x)\mathbf{Q_k}(\dif u)\\
			\text{s.t.} &  \sum_{k = 1}^K w_k \int_\reviewChanges{\supp} \|u - \bar{d_k}\|^p \mathbf{Q_k}(\dif u) \leq \epsilon^p.
    \end{cases} \\
\end{aligned}
\end{equation*}
Next, we take the Lagrangian and utilize the definition of conjugacy.
\begin{equation*}
\begin{aligned}
		&= \begin{cases}
			\underset{\mathbf{Q_k} \in \mathcal{M}(\supp)}{\sup}\underset{\lambda \geq 0}{\inf} &  \sum_{k = 1}^K w_k \int_\reviewChanges{\supp} g(u,x)\mathbf{Q_k}(\dif u) + \lambda(\epsilon^p - \sum_{k = 1}^K w_k \int_\reviewChanges{\supp} \|u - \bar{d_k}\|^p \mathbf{Q_k}(\dif u) )\\
	\end{cases}\\
	&= \begin{cases}
		\underset{\lambda \geq 0}{\inf} \underset{\mathbf{Q_k} \in \mathcal{M}(\supp)}{\sup} \quad \lambda\epsilon^p +  \sum_{k = 1}^K w_k \int_\reviewChanges{\supp} g(u,x)-\lambda\|u - \bar{d_k}\|^p \mathbf{Q_k}(\dif u)\\
\end{cases}\\
&= \begin{cases}
	\underset{\lambda \geq 0}{\inf} \quad \lambda \epsilon^p +  \underset{u = (v_1,\dots,v_K) \in \supp}{\sup} \sum_{k = 1}^K w_k (g(v_k,x)-\lambda\|v_k - \bar{d_k}\|^p),\\
\end{cases}
\end{aligned}
\end{equation*}
where the second equality is due to a well-known strong duality result for moment problems~\citep{mohajerin_esfahani_data-driven_2018}. \reviewChanges{Now, separating $g$ into its constituent functions following~\citep[Theorem 4.2]{mohajerin_esfahani_data-driven_2018},
\begin{equation*}
	\begin{aligned}
&= \begin{cases}
	\underset{\lambda \geq 0}{\inf} & \lambda\epsilon^p + \sum_{k=1} ^K  w_k s_k \\
	\mbox{subject to} &  \underset{v_k\in \supp}{\sup}\hspace{1mm} g_j(v_k, x) -\lambda \|v_k - \bar{d_k}\|^p  \le s_k \quad k = 1,\dots,K,\quad j = 1,\dots,J
\end{cases}\\
&= \begin{cases}
	\underset{\lambda \geq 0}{\inf} & \sum_{k=1} ^K  w_k s_k \\
	\mbox{subject to} &  [-g_j +\indicator_{\supp} + \lambda c_k]^*(0)  \le s_k \quad k = 1,\dots,K,\quad j = 1,\dots,J,
\end{cases}
\end{aligned}
\end{equation*}
where $c_k(v_k) \define   \| v_k - \bar{d}_k \|^p$. The last expression is identical to the form in Appendix~\ref{app:dual_form}, except now with multiple $j$, so the final dual is equivalent to the dualized constraint in~\eqref{eq:robustopt_p_max}.}

\ifpreprint
\subsection{Proof of Theorem~\ref{thm:1}}
\fi
\label{app:thm1}
We prove (i) $\bar{g}^N(x) \leq \bar{g}^K(x)$, (ii) $\bar{g}^K(x) \leq \bar{g}^{N*}(x) + (L/2) D(K)$, and (iii) when the support constraint does not affect the uncertainty set, ${\bar{g}^K(x) \leq \bar{g}^{N}(x) + (L/2) D(K)}$.

\paragraph{Proof of (i).} We begin with a feasible solution $v_1, \dots, v_{N}$ of~\eqref{eq:nocluster_one}, and set $u_k = \sum_{i \in C_k}{v_i}/|C_k|$ for each of the $K$ clusters. We see $u_k$ with $k = 1,\dots,K$ satisfies the constraints of~\eqref{eq:cluster_one}, as
\reviewChanges{
\begin{align*}
	\sum_{k = 1}^K \frac{|C_k|}{N}\left \| u_k - \bar{d}_k \right \|^p &= 	\sum_{k = 1}^K \frac{|C_k|}{N} \left \|  \frac{\sum_{i \in C_k}{v_i}}{|C_k|} -  \frac{\sum_{i \in C_k}{d_i}}{|C_k|} \right \|^p \\
	&\leq \sum_{k = 1}^K \frac{|C_k|}{N}   \sum_{i \in C_k} \frac{1}{|C_k|} \|v_i - d_i\|^p\\
	&= \sum_{k = 1}^K \frac{1}{N}  \sum_{i \in C_k} \|v_i - d_i\|^p\\
	&\leq \epsilon^p,
\end{align*}
where we have applied triangle inequality, Jensen's inequality for the convex function $f(x) = \|x\|^p$, and the constraint of~\eqref{eq:nocluster_one}. In addition, since the support $\supp$ is convex, for every $k$ our constructed $u_k$, as the average of select points $v_i$ $\in \supp$, must also be within $\supp$. The same applies with respect to the domain of $g$.
}

Since we have shown that the $u_k$'s satisfies the constraints for~\eqref{eq:cluster_one}, it is a feasible solution. We now show that for this pair of feasible solutions, in terms of the objective value, $\bar{g}^K(x) \geq \bar{g}^N(x)$.
By assumption, $g$ is concave in the uncertain parameter, so by Jensen's inequality,
\begin{align*}
\sum_{k = 1}^K \frac{|C_k|}{N} g \left(\frac{1}{|C_k|}\sum_{i \in C_k}v_i, x\right) &\geq  \sum_{k = 1}^K \frac{|C_k|}{N} \frac{1}{|C_k|}\sum_{i \in C_k} g(v_i)\\
\sum_{k = 1}^K \frac{|C_k|}{N}g(u_k,x) &\geq \frac{1}{N}\sum_{i \in N} g(v_i).
\end{align*}
Since this holds true for $u_k$'s constructed from any feasible solution $v_i, \dots, v_N$, we must have $\bar{g}^K(x) \geq \bar{g}^N(x)$.

\paragraph{Proof of (ii).}
Next, we prove $\bar{g}^K(x) \leq \bar{g}^{N*}(x) + (L/2) D(K)$ by making use of the $L$-smooth condition on $-g$. We first solve~\eqref{eq:cluster_one} to obtain a feasible solution $u_1, \dots, u_k$. We then set $\Delta_k \define u_k - \bar{d}_k$ for each $k \leq K$, and set $v_i = d_i + \Delta_k \quad \forall i \in C_k, k = 1,\dots,K$. These satisfy the constraint of~\eqref{eq:nocluster_inf}, as
\begin{align*}
	\frac{1}{N} \sum_{i=1}^N  \| v_i - d_i \|^p &= \frac{1}{N} \sum_{k=1}^K \sum_{i \in C_k} \|\Delta_k\|^p\\
	&= \sum_{k=1}^K \frac{|C_k|}{N} \|u_k - \bar{d}_k\|^p\\
	&\leq \epsilon^p,
\end{align*}
where the inequality makes use of the constraint of~\eqref{eq:cluster_one}. Since the constraints are satisfied, the constructed $v_i \dots v_N$ are a valid solution for~\eqref{eq:nocluster_inf}.
We note that these $v_i$'s are also in the domain of $g$, given that the uncertain data $\mathcal{D}_N$ is in the domain of $g$.
\reviewChanges{For monotonically increasing functions $g$, (\eg, $\log(u)$, $1/(1+u)$), we must have $\Delta_k = u_k - \bar{d_k} \geq 0$ in the solution of~\eqref{eq:cluster_one}, as the maximization of $g$ over $u_k$ wil lead to $u_k \geq \bar{d_k}$.
Therefore, $v_i = d_i + \Delta_k$ is also in the domain, as the $L$-smooth and concave function $g$ with only a potential lower bound will not have holes in its domain above the lower bound. }
For monotonically decreasing functions $g$, the same logic applies with a nonpositive $\Delta_k$.
We now make use of the convex and $L$-smooth conditions~\citep[Theorem 5.8]{amirbeck} on $-g: \forall v_1, v_2 \in \supp, \lambda \in [0,1],$
\begin{align*}
	g(\lambda v_1 + (1-\lambda) v_2) \leq \lambda g(v_1) + (1 - \lambda)g(v_2) + \frac{L}{2}\lambda(1 - \lambda)\| v_1 - v_2 \|^2_2,
\end{align*}
\reviewChanges{which, we can apply iteratively, with the first iteration being
$$ g\left (\frac{1}{|C_k|} v_1 + \frac{|C_k|-1}{|C_k|} \bar{v}_2\right ) \leq \frac{1}{|C_k|}g(v_1) + \frac{|C_k|-1}{|C_k|}g(\bar{v}_2) + \frac{L}{2}\frac{1}{|C_k|}\frac{|C_k|-1}{|C_k|}\| v_1 - \bar{v}_2 \|^2_2, $$
where $\bar{v}_2 = \frac{1}{|C_k|-1}\sum_{i \in C_k, i\neq 1} v_i$. Note that $v_1 - \bar{v}_2 = d_1 - \frac{1}{|C_k|-1}\sum_{i \in C_k, i\neq 1} d_i$, as they share the same $\delta_k$. The next iteration will be applied to $g(\bar{v}_2)$, and so on.} 
\reviewChanges{
For each cluster $k$, this results in:
\begin{align*}
	g \left(\frac{1}{|C_k|}\sum_{i \in C_k} v_i, x\right) &\leq  \frac{1}{|C_k|}\sum_{i \in C_k} g(v_i,x) + \frac{L}{2|C_k|} \sum_{i=2} ^{|C_k|} \frac{i-1}{i} \left \|d_i - \frac{\sum_{j=1}^{i-1}d_j}{i-1} \right \|^2_2\\
	 g(\bar{d}_k + \Delta_k,x) &\leq \frac{1}{|C_k|}\sum_{i \in C_k} g(v_i,x) + \frac{L}{2|C_k|} \sum_{i \in C_k} \|d_i - \bar{d}_k\|^2_2\\
	g(u_k,x) &\leq  \frac{1}{|C_k|}\sum_{i \in C_k} g(v_i,x) + \frac{L}{2|C_k|} \sum_{i \in C_k} \|d_i - \bar{d}_k\|^2_2,
\end{align*}
where we used the equivalence $$ \sum_{i=2} ^{|C_k|} \frac{i-1}{i} \left \|d_i - \frac{\sum_{j=1}^{i-1}d_j}{i-1} \right \|^2_2 =  \sum_{i \in C_k} \|d_i - \bar{d}_k\|^2_2. $$
Now, summing over all clusters, we have 
\begin{align*}
	\sum_{k = 1}^K \frac{|C_k|}{N}  g(u_k,x) &\leq \frac{1}{N} \sum_{i = 1}^N g(v_i,x) + (L/2) D(K).
\end{align*}
Since this holds for any feasible solution of~\eqref{eq:cluster_one}, we must have  $ {\bar{g}^K(x) \leq \bar{g}^{N*}(x) + (L/2) D(K)}$.
}

\ifpreprint
\subsection{Proof of Theorem~\ref{thm:maxconcave}}
\fi
\label{app:maxconcave}
\reviewChanges{
\paragraph{Proof of the lower bound.}
We use the dual formulations of the \gls{MRO} constraints. We first solve~\eqref{eq:cluster_dual} to obtain dual variables $z_{jk}$, $\gamma_{jk}$ for $k = 1,\dots,K$, $j = 1,\dots,J$. For each data label $i$ in cluster $C_k$, for all clusters $k = 1,\dots,K$, and for all pieces $j = 1,\dots, J$, if we set
\begin{equation*}
	\begin{aligned}
		z_{ji} &= z_{jk}, \quad \gamma_{ji} = \gamma_{jk}\\
		s_i &= s_k + \max_j\{(-z_{jk} - H^T\gamma_{jk})^T (d_i - \bar{d}_k)\},
	\end{aligned}
\end{equation*}
we have obtained a valid solution for~\eqref{eq:nocluster_dual}. The increase in the objective value from $\bar{g}^K(x)$ to that of the constructed solution of~\eqref{eq:nocluster_dual} is
\begin{equation*}
	\begin{aligned}
		\delta(K,z,\gamma) & \define (1/N)\sum_{i=1}^N s_i - \sum_{k=1}^K(|C_k|/N)s_k \\
		&= (1/N)\sum_{k = 1}^K\sum_{i\in C_k} \max_j\{(-z_{jk} - H^T\gamma_{jk})^T(d_i - \bar{d}_k)\}.
	\end{aligned}
\end{equation*}
We note that $\delta(K,z,\gamma) \geq 0$, as
\begin{equation*}
	\begin{aligned}
		\delta(K,z,\gamma) &= (1/N)\sum_{k = 1}^K\sum_{i\in C_k} \max_j\{(-z_{jk} - H^T\gamma_{jk})^T(d_i - \bar{d}_k)\}\\
		&\geq (1/N)\sum_{k = 1}^K\sum_{i\in C_k} (-z_{1k}^T - \gamma_{1k}^TH)(d_i - \bar{d}_k)\\
		& = (|C_k|/N)\sum_{k = 1}^K  (-z_{1k} - H^T\gamma_{1k})^T (1/|C_k|)\sum_{i\in C_k} (d_i - \bar{d}_k)\\
		&= (|C_k|/N) \sum_{k = 1}^K  (-z_{1k} - H^T\gamma_{1k})^T0\\
		&= 0
	\end{aligned}
\end{equation*}
The constructed feasible solution for~\eqref{eq:nocluster_dual} is an upper bound for its optimal solution, since we solve a minimization problem. As this holds true for any solution of~\eqref{eq:cluster_dual}, we must have $\bar{g}^K(x) + \delta(K,z,\gamma) \geq \bar{g}^N(x)$, which translates to $\bar{g}^N(x) - \delta(K,z,\gamma) \leq \bar{g}^K(x)$.
\paragraph{Proof of the upper bound.}
We use the primal formulations of the \gls{MRO} constraints. 
We first solve~\eqref{eq:cluster_one} to obtain a feasible solution $u_{11}, \dots, u_{JK}$. We then set $\Delta_{jk} \define u_{jk} - \bar{d}_k$ for each $k \leq K$, and set $v_{ji} = d_i + \Delta_{jk}$, $\alpha_{ji} = \alpha_{jk}/|C_k| \quad \forall i \in C_k, k = 1,\dots,K$.
These satisfy the constraint of~\eqref{eq:nocluster_inf}, as
\begin{align*}
	\sum_{i=1}^N \sum_{j=1}^J  \alpha_{ji} \| v_{ji} - d_i \|^p &= \sum_{k=1}^K \sum_{j=1}^J \sum_{i \in C_k}  \alpha_{ji} \|\Delta_{jk}\|^p\\
	&= \sum_{k=1}^K \sum_{j=1}^J \alpha_{jk} \|u_{jk} - \bar{d}_k\|^p\\
	&\leq \epsilon^p.
\end{align*}
The constructed solutions remain in $\dom$, following the arguments in the proof of (ii) in Appendix~\ref{app:thm1}. 
Now, for each cluster $k$ and function $g_j$, we note using the $L$-smooth condition on $-g_j$,
\begin{align*}
	\alpha_{jk} g_j\left(\frac{1}{|C_k|}\sum_{i \in C_k} v_{ji}, x\right) &\leq  \alpha_{jk} \left (\sum_{i \in C_k} \frac{1}{|C_k|} g_j(v_{ji},x) + \frac{L_j}{2|C_k|} \sum_{i=2} ^{|C_k|} \frac{i-1}{i} \left \|d_i - \frac{\sum_{j=1}^{i-1}d_j}{i-1} \right \|^2_2\right )\\
	\alpha_{jk} g_j(\bar{d}_k + \Delta_k,x) &\leq  \sum_{i \in C_k} \alpha_{ji} g_j(v_{ji},x)  + \frac{\alpha_{ji} L_j}{2} \sum_{i \in C_k} \|d_i - \bar{d}_k\|^2_2
\end{align*}
Then, summing over all the clusters and functions, we have 
\begin{align*}
\sum_{k=1}^K \sum_{j=1}^J \alpha_{jk}  g_j(u_k,x) &\leq \sum_{k=1}^K \sum_{j=1}^J \sum_{i \in C_k} \alpha_{ji} g_j(v_{ji},x)  + \sum_{k=1}^K \sum_{j=1}^J \frac{\alpha_{ji}  \max_{j\leq J} L_j}{2} \sum_{i \in C_k} \|d_i - \bar{d}_k\|^2_2.\\
&\leq \sum_{i=1}^N \sum_{j=1}^J \alpha_{ji} g(v_i,x) + \max_{j\leq J} (L_j/2N)\sum_{i=1}^N \|d_i - \bar{d}_k\|^2_2 .
\end{align*}
Since this holds for all feasible solutions of~\eqref{eq:cluster_one}, we conclude that 
$$ \bar{g}^K(x) \leq \bar{g}^{N*}(x) +  \max_{j\leq J} (L_j/2) D(K).$$
}

\ifpreprint
\subsection{Conjugate derivation for the Capital Budgeting problem~\eqref{eq:capital}}
\fi
\label{app:capital}
We begin with
\begin{equation*}
g(u,x) = -\eta(u)^Tx = -\sum_{j = 1}^n \sum_{t = 0}^T F_{jt}x_j/(1+u_j)^{t},
\end{equation*}
and take its conjugate $[-g]^*(z)$ in the uncertain parameter $u$. We use the theorem on infimal convolutions~\citep[Theorem 11.23(a), p. 493]{rockafellar_wets_1998} to arrive at
\begin{equation*}
	[-g]^*(z) = \sum_{t = 0}^T \sup_{u} \left( u^Ty_t -  \left[\sum_{j = 1}^n F_{jt}x_j(1+ u_j)^t\right]\right), \quad \sum_{t = 0}^T y_t = z.
	\end{equation*}
We calculate the inner conjugates using the the first order optimality condition,
\begin{align*}
\nabla \left(y_t^Tu - \sum_{j = 1}^n F_{jt}x_j(1+ u_j)^{-t}\right) = 0 \\
y_{jt} = -tF_{jt}x_j(1+u_j^\star)^{-(t+1)}, \quad j = 1,\dots,n \\
u_j^\star = (tF_{jt}x_j(-y_{jt})^{-1})^{1/(t+1)} - 1, \quad j = 1,\dots,n.
\end{align*}
Substituting this back into the expression for the conjugate, we have, for each $j$ and $t$,
\begin{align*}
	y_{jt}u_j^\star - F_{jt}x_j&(1+u_j)^{-t} = y_{jt}(tF_{jt}x_j(-y_{jt}^{-1})^{(1/(t+1)} - 1)) \\
	&\hspace{4cm} - F_{jt}(tF_{jt}x_j(-y{jt})^{-1})^{-t/(t+1)}\\
&= -y_{jt} - ((-y_{jt})^{t/(t+1)}(F_{jt}x_j)^{1/(t+1)})(t^{1/(t+1)} + t^{-1/(t+1)}),
	\end{align*}
after combining terms. Note that $-(-y_{jt})^{t/(t+1)}(F_{jt}x_j)^{1/(t+1)}$ can be replaced by an auxiliary variable $\delta_{jtk} \leq 0$, by introducing the power cone constraint $(-y_{jt})^{t/(t+1)}(F_{jt}x_j)^{1/(t+1)} \geq |\delta_{jtk} |$. By substituting these results into~\eqref{eq:robustopt_p} and further vectoring some constraints, we arrive at the desired formulation.

\ifpreprint \else
\end{appendices}
\fi
\end{document}